\documentclass{amsart}

\usepackage[foot]{amsaddr}
\usepackage{amsmath,amsfonts,amssymb,amsthm,amscd,comment,euscript,
latexsym,mathrsfs,amsfonts,amscd,euscript,graphicx,lastpage}
\usepackage[all]{xy}
\usepackage[colorlinks=true, pdfstartview=FitV, linkcolor=blue,
citecolor=blue,urlcolor=blue,breaklinks=true]{hyperref}
\usepackage{color}
\usepackage[top=.8in,bottom=1in,left=1in,right=1in]{geometry} 
\usepackage{fancyhdr}
\usepackage{paralist}

\definecolor{yellowtimes}{rgb}{1,1,0.65}

\newcommand{\comments}[1]{}

\newcommand\vfe{\medskip} %
\newcommand{\ignore}[1]{}  

\newcommand{\arxiv}[1]{\href{http://arxiv.org/abs/#1}{\tt arXiv:\nolinkurl{#1}}}

\newcommand\C{\mathbb{C}}
\newcommand\Z{\mathbb{Z}}

\newcommand\R{\mathbb{R}}
\newcommand\N{\mathbb{N}}
\newcommand\kk{{\Bbbk}}

\newcommand\id{\mathrm{id}}

\newcommand\reg{\mathrm{reg}}

\newcommand\g{\mathfrak{g}}
\newcommand\h{\mathfrak{h}}
\newcommand\n{\mathfrak{n}}

\newcommand\cI{\mathcal{I}}

\newcommand\bW{\mathbf{W}}

\newcommand\bfR{\mathbf{Res}}

\newcommand{\ts}{\textstyle}

\newcommand{\op}{^\mathrm{op}}
\DeclareMathOperator{\Hom}{Hom}

\DeclareMathOperator{\End}{End}

\DeclareMathOperator{\Span}{Span}
\DeclareMathOperator{\Spec}{Spec}

\DeclareMathOperator{\Id}{Id}
\DeclareMathOperator{\Ann}{Ann}
\DeclareMathOperator{\supp}{supp} 
\DeclareMathOperator{\Ob}{Ob}

\DeclareMathOperator{\ad}{ad}
\DeclareMathOperator{\wt}{wt}
\DeclareMathOperator{\Ker}{Ker}

\DeclareMathOperator{\U}{U}

\DeclareMathOperator{\diag}{diag}
\DeclareMathOperator{\Mat}{Mat}
\DeclareMathOperator{\can}{can}
\DeclareMathOperator{\tr}{tr}

\theoremstyle{plain}

\newtheorem*{theo*}{Theorem}

\newtheorem*{lem*}{Lemma}

\theoremstyle{definition}

\newtheorem*{rem*}{Remark}

\numberwithin{equation}{subsection}

\newcommand*{\defeq}{\mathrel{\vcenter{\baselineskip0.5ex \lineskiplimit0pt
                     \hbox{\scriptsize.}\hbox{\scriptsize.}}}%
                     =}
\allowdisplaybreaks
\newcommand{\sgn}{\mathrm{sgn}}
\newcommand{\sym}{\mathrm{sym}}

\newcommand{\bfT}{\mathbf{T}}

\DeclareSymbolFont{bbold}{U}{bbold}{m}{n}
\DeclareSymbolFontAlphabet{\mathbbold}{bbold}
\newcommand{\idu}{\mathbbold{1}_{\rm U}}
\newcommand{\TS}{\mathrm{TS}}

\newcommand\sm{\smallskip}
\newcommand\ms{\medskip}

\newcommand\ot{\otimes}

\newcommand\Iff{\Leftrightarrow}
\newcommand\lan{\langle} \newcommand\ran{\rangle}

\newcommand\tsum{\textstyle \sum}

\newcommand\ch{\sp{\scriptscriptstyle\vee}}

\newcommand\ideal{\triangleleft}
\newcommand\wid{\widetilde} 
\newcommand\co{\colon}

\newcommand{\simlgr}{\buildrel \sim \over \longrightarrow}

 \DeclareMathOperator{\rma}{A}
 \DeclareMathOperator{\rmb}{B}
 \DeclareMathOperator{\rmbc}{BC}
 \DeclareMathOperator{\rmc}{C}
 \DeclareMathOperator{\rmd}{D}
 \DeclareMathOperator{\rme}{E}

\newcommand\bfI{{\mathbf I}\mathbf{n}\mathbf{t}}

\newcommand\CC{\mathbb{C}}

\newcommand\NN{\mathbb{N}}

\newcommand\Se{\mathbb{S}}
\newcommand\ZZ{\mathbb{Z}}

\newcommand\scC{\mathcal{C}}

\newcommand\scI{\mathcal{I}}

\newcommand\scP{\mathcal{P}}

\newcommand{\euP}{\EuScript{P}}
\newcommand\euQ{\EuScript{Q}}

\newcommand\euZ{\EuScript{Z}}  

\newcommand\frC{\mathfrak{C}}

\newcommand\frn{\mathfrak{n}} 

\newcommand\frs{\mathfrak{s}}
\newcommand\frS{\mathfrak{S}} 

\newcommand\frh{\ensuremath{\mathfrak{h}}} 

\newcommand\lsl{\ensuremath{\mathfrak{sl}}}
\newcommand\gl{\ensuremath{\mathfrak{gl}}}

\newcommand\sfB{{\sf B}}
\newcommand\sfC{{\sf C}}

\newcommand\sfS{{\sf S}}

\newcommand\sfT{{\sf T}} \newcommand\TB{\sfT\sfB}

\newcommand\al{\alpha}
\newcommand\be{\beta}
\newcommand\ga{\gamma} 
\newcommand\de{\delta} \newcommand\De{\Delta}
\newcommand\eps{\epsilon} 
\newcommand\io{\iota}

\newcommand\la{\lambda} 
 \newcommand\vphi{\varphi}

\newcommand\si{\sigma} 

\newcommand\ta{\tau}
  \newcommand\Th{\Theta}
\newcommand\ze{\zeta}
\newcommand\om{\omega}

\newcommand{\kalg}{\kk\mathchar45\mathbf{alg}}
\newcommand{\kkalg}{{\Bbbk}\mathchar45\mathbf{alg}}
\newcommand{\MOD}{\text{-}\ensuremath\mathcal{M}\mathit{od}}
\makeatletter
\renewcommand{\email}[2][]{%
  \ifx\emails\@empty\relax\else{\g@addto@macro\emails{,\space}}\fi%
  \@ifnotempty{#1}{\g@addto@macro\emails{\textrm{(#1)}\space}}%
  \g@addto@macro\emails{#2}%
}
\makeatother

\begin{document}
%
\title[]{Integrable representations of root-graded Lie algebras}
\author{Nathan Manning}
\address[ Manning]{Department of Mathematics, University of Maryland}
\email{nmanning@math.umd.edu}
\author{Erhard Neher}
\address[Neher, Salmasian]{Department of Mathematics and Statistics, University of Ottawa}
\email{neher@uottawa.ca}
\author{Hadi Salmasian}
\email{hsalmasi@uottawa.ca}
\thanks{}


\date{\today}


\date{\today}


\begin{abstract}
In this paper we study the category of representations of a root-graded Lie
algebra $L$ which are integrable as representations of a finite-dimensional
semisimple subalgebra $\g$ and whose weights are bounded by some dominant
weight of $\g$. We link this category to the module category of an
associative algebra, whose structure we determine for map algebras and
$\lsl_n(A)$.

Our approach unifies recent work of Chari and her collaborators on map
algebras, of Fourier and Savage and their collaborators on equivariant map
algebras, as well as the classical work of Seligman on isotropic Lie
algebras.
\end{abstract}

 \maketitle \thispagestyle{empty}
\centerline{\it Dedicated to Efim Zelmanov on the occasion of his 60th
birthday}

\tableofcontents


\maketitle \thispagestyle{empty}
\tableofcontents

%
\section*{Introduction}
%

This paper unites two strands of research which up to now did not have any
interaction. Historically, the first strand is the work of Seligman on
representations of finite-dimensional isotropic central-simple Lie algebras
$L$ over fields $\kk$ of characteristic $0$ (\cite{Sel81, Sel88}). While
these could be studied in the spirit of Tits' approach to representation
theory \cite{tits-rep} using Galois descent, Seligman pursued ``rational
methods'', i.e., describing the representations over $\kk$
(rather than the algebraic closure of $\kk$).  The Lie algebra $L$ decomposes
with respect to a maximal split toral subalgebra $\h\subset L$ as $L=
\bigoplus_{\al \in \Th \cup \{0\}} L_\al$ where $\Th \subset \h^*$ is an
irreducible, possibly non-reduced root system and the $L_\al$ are the root
spaces of $L$ with respect to $\h$. Any finite-dimensional irreducible
representation of $L$ has weights bounded by a dominant integral weight $\la$
of $\Th$. Seligman introduced a unital associative algebra $\Se^\la$  (see
Definition \ref{def:adm})  and linked the finite-dimensional irreducible
representations of $L$ and of $\Se^\la$. \sm

The second strand is essentially due to Chari and her collaborators. It
concerns integrable representations of map algebras, i.e., Lie algebras
$L=\g\ot_\CC A$ where $\g$ is a finite-dimensional simple complex Lie algebra
and $A$ is a unital commutative associative $\CC$-algebra. (The name ``map
algebra'' comes from the interpretation of $L$ as regular maps from the
affine scheme $\Spec(A)$ to the affine variety $\g$.) While originally
$A=\CC[t^{\pm 1}]$ (loop algebras \cite{CP01}) or $A=\CC[t]$ (current
algebras \cite{CL06}), the algebra $A$ was soon taken to be arbitrary, see
for example \cite{FL04} for an earlier paper.
 Note that the Lie algebra
$L=\g \ot_\CC A$ decomposes in a similar way as Seligman's Lie algebras: Let
$\De$ be the root system of $\g$ with respect to a Cartan subalgebra $\h$ of
$\g$ so that $\g = \bigoplus_{\al \in \De \cup \{0\}} \g_\al$. Then $L=
\bigoplus_{\al \in \De} L_\al$ where $L_\al = \g_\al \ot_\CC A$ are the
weight spaces of $L$ under the canonical $\h$-action. Regarding the
representation theory of map algebras, a  big step forward was
undertaken in the paper \cite{CFK10} by Chari-Fourier-Khandai, which stressed
a categorical point of view. Denote by $\scI^\la(L,\g)$ the category of
integrable representations of $L=\g \ot_\CC A$ whose weights with respect to
a Cartan subalgebra $\h$ of $\g$ are bounded by a dominant weight $\la$ for
$(\g,\h)$. It was shown in loc.\ cit., among many other things, that
$\scI^\la(L, \g)$ is closely related to the module category of a commutative
associative $\CC$-algebra $\mathbf A_\la$ (see Definition \ref{def:cha}).
Moreover, the algebra $\mathbf A_\la$ was determined in the Noetherian case.
The paper also introduced global Weyl modules as the initial cyclic objects
in $\scI^\la(L,\g)$, generalizing the case $A=\CC[t^{\pm 1}]$ of \cite{CP01}.
It has been the blueprint for many sequels, in which the results of
\cite{CFK10} have been generalized to the setting of various equivariant map
algebras, culminating with the recent paper \cite{FMSb}. \sm

The starting point of this paper is the observation that map algebras $\g \ot
A$ and the equivariant map algebras studied in the sequels to
\cite{CFK10} like for example \cite{FMSb},  as well as the Lie algebras
considered in \cite{Sel81, Sel88}, are examples of Lie algebras graded by
finite (not necessarily reduced nor irreducible) root systems. For the sake
of simplicity, in this introduction we will only consider irreducible root
systems $\Th$. Given such a $\Th$ we let $\De=\Th$ if $\Th$ is reduced and
$\De = \rmb_n$ or $= \rmc_n$ if $\Th \cong \rmbc_n$. By definition, a Lie
algebra $L$ over a field $\kk$ of characteristic $0$ is {\em
$(\Th,\De)$-graded\/} if it contains a finite-dimensional split simple
subalgebra $\g$ whose root system with respect to a splitting Cartan
subalgebra $\h\subset \g$ is $\De$ such that
\[ L = \textstyle \bigoplus_{\al \in \Th \cup \{0\}} L_\al \quad \hbox{and} \quad
\quad L_0 = \sum_{\al \in \Th} [L_\al, L_{-\al}],
\]
where the $L_\al$ are the weight spaces of $L$ under the adjoint
representation of $\h$. There exists an elaborate structure theory of
root-graded Lie algebras \cite{abg,abg2,BeSm,beze, BM92,n:3g}. The reader
will recognize that the two types of Lie algebras mentioned above are indeed
$(\Th,\De)$-graded for appropriate choices of $(\Th,\De)$. For example, the
Lie algebras considered in \cite{CIK, FK12} are $(\Th,\De)$-graded with a
non-reduced $\Th$. But there exist many more examples of root-graded Lie
algebras (see Section \ref{rg-exa} for further examples). \sm

The main point of this paper is that root-graded Lie algebras provide a
natural setting, in which one can develop a satisfactory theory of integrable
representations. Indeed, we are able to define the main concepts developed
for representations of map algebras on the one hand and for
finite-dimensional central-simple Lie algebras on the other hand in the
setting of $(\Th,\De)$-graded Lie algebras $(L,\g)$:
\begin{enumerate}

 \item[-] The category $\scI(L,\g)$ (defined in Section \ref{cat-int-def})
 is the category of $L$-modules $V$ that are integrable as $\g$-modules (See Section \ref{int-mod}) and
thus have a weight space decomposition $V= \bigoplus_{\mu \in \euP} V_\mu$
with respect to $\h$ for $\euP$ the set of integral weights of $\De$.

  \item[-] $\scI^\la(L,\g)$ is the subcategory of $\scI(L,\g)$ consisting
      of modules $V$ for which all non-zero weight spaces have weights
      bounded above by $\la$ in a canonical order (which is the natural
      partial order of $\euP$ in the reduced case), see Section \ref{def:cat-I-bd}.

 \item[-] Global Weyl modules $W(\la) \in \scI(L,\g)^\la$ are defined by
     taking an appropriate quotient of $\U(L) \ot _{\U(\g)} V(\la)$
     (\ref{def:global-weyl}) and have the usual presentation
     (See Proposition \ref{prop:global-weyl-def-rel}).

 \item[-] In Definition \ref{def:cha}, we define the
     algebra $\mathbf{A}_\la = \U(L_0)/ \Ann_{U(L_0)} (w_\la)$,
where $w_\la \in W(\la)$ is the canonical generator of $W(\la)$.
     This
     extends  the definition originally given by Chari and Pressley.

 \item[-] In Definition \ref{def:adm}, we define the algebra $\Se^\la =
     \U(L_0)/J^\la$, where the ideal $J^\la$ is given in terms of
     generators  -- it turns out to be the maximal quotient of
  $\U(L_0)$ which acts on the weight space $V_\la$ for any module $V\in
  \scI^\la(L,g)$. The algebra $\Se^\la$ extends the definition of Seligman.
\end{enumerate}
The main results of this paper are: \sm

(a) Theorem~\ref{thm:chari=seli}: {\em $\mathbf
    A_\la$ and $\Se^\la$ are equal
associative algebras}, that is, $J^\la=\Ann_{U(L_0)} (w_\la)$.
Because of this equality,  this algebra will be called the
Seligman--Chari--Pressley algebra.
\sm

(b) Extending results from the previously considered cases, we show that we
have a restriction
    functor \[ \bfR \co \scI^\la(L,\g) \to \Se^\la\MOD,\]
     given by $V \mapsto
    V_\la$, and an ``integrable induction'' functor \[ \bfI \co \Se^\la\MOD
    \to \scI^\la(L,\g)\] that satisfy $\bfR \circ \bfI \xrightarrow{\simeq}
    \Id_{\Se^\la\MOD}$, see \eqref{func1}, and $\bfI \circ \bfR \Rightarrow
    \Id_{\scI^\la(L,\g)}$, see Proposition \ref{poth}\eqref{pothb}. \sm

(c) We determine the structure of $\Se^\la$ in two important examples.
    For $L=\g \ot A$, where $A$ is a commutative unital $\kk$-algebra,  and a dominant $\la = \sum_{i\in I} \ell_i \varpi_i$
    we have
\[ \Se^\la \cong \TS^{\ell_1}(A) \ot_\kk \cdots \ot_\kk \TS^{\ell_r}(A)
\]
(Theorem~\ref{thm:comm}), where $\TS^n(A)$ is the fixed point subalgebra of
$A^{\ot n}$ under the obvious action of the symmetric group $\frS_n$. This
generalizes \cite[Thm.~4]{CFK10}. For the $\rma_{n-1}$-graded Lie algebra
$\lsl_n(A)$, $n\ge 3$ and $A$ associative but not necessarily commutative,
and a dominant weight  $\la$ with totally disconnected support, we describe
$\Se^\la$ in Theorem~\ref{disc-th}. For example, taking $\la = \ell \varpi_i$
with $\ell\in \NN_+$ we have
\[
    \Se^{\ell \varpi_i} \cong
        \begin{cases} \TS^\ell(A), & i=1, \\
        \TS^\ell(A)/\scC & 1<i<n-1, \\
          \TS^\ell(A\op) & i=n-1. \end{cases}
\]
where $\scC$ is the ideal of $\TS^\ell(A)$ generated by the commutator space
$[\TS^\ell(A), \TS^\ell(A)]$. We point out that it may very well happen that
$\Se^\la = \{0\}$ in the middle case, which means that $V_\la = \{0\}$ for
every $V\in \scI^\la$.

\sm

 Besides the classes of root-graded Lie algebras explicitley mentioned
above, Tits-Kantor-Koecher algebras $L(J)$ of unital Jordan algebras $J$ and
their central coverings are another important class of root-graded Lie
algebras --- in this case $\Th=\De=\rma_1$. Integrable representations of the
universal central extension  $\widehat{L(J)}$   bounded by $\la = 2\varpi$
have been studied in the recent preprint \cite{KS} by Kashuba-Serganova.
Results for $\la=\varpi$ are contained in the sequel \cite{KS-II} of
\cite{KS}, while integrable representations in general are investigated in
on-going work of Lau and Mathieu. The cases $\la=\om$ and $\la=2\om$ are
closely related to associative specializations and bimodules of Jordan
algebras $J$, which for finite-dimensional simple $J$ have been classified in
the classical paper \cite{JJ} by F.~D.\ and N.~Jacobson (see also
\cite[Ch.~II and VII]{jake:srj}). Representations of $\rmc_n$-graded Lie
algebras have been investigated in \cite{Zel-85, Zel-94}. \sm

{\em Outlook.} There are many questions left open in this paper. For
instance,  what is structure of the Seligman-Chari-Pressley algebra for the
other types of root-graded Lie algebras? What are local Weyl modules in the
setting of root-graded Lie algebras?  Lie tori form a special
class of root-graded Lie algebras that enter in the construction of extended
affine Lie algebras, see e.g.\ \cite{ n:eala,   n:tori, Neh11}. It would be
of interest to work out the relation between the integrable representations
of Lie tori and the representations of the associated extended affine Lie
algebras. Last but not least, one can also consider these questions for
root-graded Lie superalgebras, for which one has a well-developed structure
theory (see \cite{BEM} and the references therein). \sm

The paper is organized as follows. In Section~\ref{sec:sta} we study
symmetric tensor algebras $\TS^\ell(A)$ for $A$ a unital associative
$\kk$-algebra. Generalizing a result proven in \cite[\S2]{Sel80} for
finite-dimensional $A$, we show in Theorem~\ref{sel-thm} that $\TS^\ell(A)$
has a universal property with respect to certain symmetric identities in the
sense of \cite{Sel80}. In light of the explicit descriptions of $\Se^\la$ given above, understanding the
structure of $\TS^\ell(A)$ is important. We discuss this in Section \ref{mopro}. For
example, $\TS^\ell\big(\Mat_d(\kk)\big)$ is described by the classical
Schur-Weyl duality.
In Section~\ref{sec:gen} we investigate integrable
representations of certain pairs $(L,\g)$ which in Section~\ref{sec:rg} are
specialized to root-graded Lie algebras. Section~\ref{sec:admiss} is
devoted to the Seligman algebra $\Se^\la=\Se^\la(L,\g)$, its module category and the
link to $\scI^\la(L,\g)$. For example, when $ L_0$ is finite-dimensional,  we show that $\dim \Se^\la$ is finite and indeed we give an explicit upper bound for
the  dimension of $\Se^\la$
(see Proposition~\ref{fd-Se}).
 In Sections \ref{sec:rep-gencur} and
\ref{sec:slna} we prove our results on $\Se^\la$ for map algebras and
$\lsl_n(A)$ respectively.
We finish by investigating the structure of $\Se^\la$ for $\lsl_4(A)$ and $\la=\varpi_1+\varpi_2$.
\sm

\paragraph{\textbf{Acknowledgements}} The authors thank Vyjayanthi Chari for encouragement and
hints to the literature. We are also grateful to Alistair Savage for
comments on a preliminary version of this paper.

The second and third authors gratefully acknowledge research funding from the
NSERC Discovery Grant program. During part of the work on this paper the
first and second author enjoyed the hospitality of the Centre de Recherches
Math\'ematiques (Universit\'e de Montr\'eal) as members of the thematic
program ``New Directions in Lie Theory''. They gratefully acknowledge support
from the CRM during that program.
%
%
\ms

\subsection*{Notation and conventions}\label{notation-section}
Throughout the paper, $\kk$ will denote a field of characteristic zero, $\Z$
the integers, $\N$ the nonnegative integers and $\N_+$ the positive integers.
All algebras will be defined over $\kk$. Unless stated otherwise we
abbreviate $\ot = \ot_\kk$. If $X$ is a set, its cardinality will be denoted
by $|X|$. We use $M \in \Ob \scI$ to denote an object $M$ of a category
$\scI$. \sm

Throughout $A$ denotes an associative unital but not necessarily commutative
$\kk$-algebra. Its identity element is denoted $1_A$ or sometimes just $1$.
Given another unital associative algebra $B$, a unital algebra homomorphism
$\vphi \co A \to B$ is a $\kk$-linear map satisfying $\vphi(a_1 a_2) =
\vphi(a_1) \vphi(a_2)$ for all $a_i \in A$ and $\vphi(1_A) = 1_B$. As usual,
$[a_1, a_2] = a_1a_2 - a_2 a_1$ is the commutator in $A$. We denote by
$\kkalg$ the category of associative commutative unital $\kk$-algebras with
unital algebra homomorphisms as morphisms. We abbreviate $A\in \kkalg$ for
$A\in \Ob\kkalg$.

Any associative algebra $B$ becomes a Lie algebra with respect to the
commutator product $[.,.]$, denoted $B^-$. Its derived algebra $[B,B] =
\Span_k\{ [b_1, b_2]: b_i \in B\}$ is an ideal of the Lie algebra $B^-$. A
linear map $A \to B$ will be called a {\em Lie homomorphism\/} if it is a
homomorphism $A^- \to B^-$ of the associated Lie algebras.

The universal enveloping algebra of a Lie algebra $L$ is denoted $\U(L)$. We
will say that $V$ is an $L$-module if there exists a Lie algebra homomorphism
$\rho \co L \to \gl(V)$. However, if we need to be more precise this
situation will be abbreviated by $(V,\rho)$. Unless explicitly stated
otherwise, $\g$ is a finite-dimensional split semisimple Lie algebra over
$\kk$. All other unexplained notation and terminology can be found in
\ref{g-not}.

\vfe

\section{Symmetric tensor algebras} \label{sec:sta}

In this section we review and generalize some of Seligman's results
\cite{Sel80} on symmetric tensor algebras. Throughout this section we fix
$\ell \in \NN_+$.

\subsection{Definition (Symmetric tensor algebras)}\label{subsn:symm-alg}
We denote by $A^{\ot \ell}= A \ot \cdots \ot A$ the $\ell^{\rm th}$-tensor
product of $A$. We will view $A^{\otimes \ell}$ as an associative algebra
whose product is given by ``coordinate-wise'' multiplication, i.e., by
extending linearly the assignment $(a_1\otimes\cdots\otimes a_\ell)\cdot
(b_1\otimes\cdots\otimes b_\ell)\defeq (a_1b_1)\otimes\cdots\otimes (a_\ell
b_\ell)$. The symmetric group $\frS_l$ acts by automorphisms on $A^{\ot
\ell}$ by linear extension of $\si \cdot (a_1\otimes\cdots\otimes
a_\ell)=a_{\sigma^{-1}(1)}\otimes\cdots\otimes a_{\sigma^{-1}(\ell)}$.

We denote by $\TS^\ell(A)$ the subspace of symmetric tensors in $A^{\otimes
\ell}$, i.e., the invariant subspace under the $\frS_\ell$-action ($\TS$
stands for ``tenseurs sym\'etriques''). Since $\frS_\ell$ acts on $A^{\ot
\ell}$ by  automorphisms, the subspace $\TS^\ell(A)$ is a subalgebra of
$A^{\otimes \ell}$. In the following, we will always view $\TS^\ell(A)$ with
this associative structure. We call $\TS^\ell(A)$ the {\em $\ell^{\rm
th}$-symmetric tensor algebra.}

Since $\TS^\ell(A)$ is an associative algebra, $\TS^\ell(A)^-$ is a Lie
algebra. For $1 \le i \le \ell$ and $a\in A$ we put
\begin{align*}
   a^{\ot \ell} &= a \ot a \ot \cdots \ot a \in \TS(A) \quad \hbox{($\ell$ factors), and}\\
  s_i(a) &=  1^{\otimes i-1}\otimes a\otimes 1^{\otimes \ell-i} \in A^{\ot \ell}.
\end{align*}
Since $[s_i(a), s_j(b)]=\delta_{ij} s_i([a,b])$ the \textit{symmetrization
map}
\begin{equation} \label{subsn:symm-alg1}
\sym_\ell: A^-\to \TS^\ell(A)^-,\quad a\mapsto \textstyle
   \sum_{i=1}^\ell  s_i(a)\end{equation}
is a homomorphism of Lie algebras.

\subsection{Lemma}\label{lem:symm-folklore} {\em 
\begin{inparaenum}[\upshape (a)]
\item\label{sym-folk-a} $\TS^\ell(A)$ is spanned as vector space by $a^{\ot
    \ell}$, $a\in A$.

\item\label{sym-folk-c} $\TS^\ell(A)$ is generated as an associative algebra
    by $\sym_\ell(A)$.
\end{inparaenum} }
\sm

\noindent  More properties of $\TS^\ell(A)$ are discussed in \ref{mopro}.

\begin{proof} \eqref{sym-folk-a} We denote by $\sfS(A)$ the symmetric algebra
of the vector space underlying $A$, by $\circ$ its product and by
$\sfS^\ell(A)$ its $\ell^{th}$ graded component. Since our base field $\kk$
has characteristic 0, the canonical map $\pi: A^{\otimes \ell}\to
\sfS^\ell(A)$, given by $a_1\otimes\cdots\otimes a_\ell\mapsto
a_1\circ\cdots\circ a_\ell$, becomes a linear isomorphism when restricted to
$\TS^\ell(A)$ \cite[III, \S6.3 Rem.]{Bou70}. We have
$\pi(a^{\otimes\ell})=a\circ a\cdots\circ a$. Since $\sfS^\ell(A)$ is spanned
by the powers $a\circ \cdots \circ a$ \cite[III, \S6.1, Rem.\,3]{Bou70},
\eqref{sym-folk-a} follows.

\eqref{sym-folk-c} Fix $a\in A$. By \eqref{sym-folk-a} it is enough to show
that $a^{\ot \ell}$ lies in the associative subalgebra of $\TS^\ell(A)$
generated by $\sym_\ell(A)$. Since $[s_i(a), s_j(a)] = 0$ for $1\le i, j\le
\ell$, we have a unique unital algebra homomorphism $\ze \co \kk[x_1, \ldots
, x_\ell] \to A^{\ot \ell}$ such that $\ze(x_i) = s_i(a)$. Recall that
$\frS_l$ acts on the polynomial algebra $\kk[x_1, \ldots, x_\ell]$ by $(\si
\cdot p)(x_1, \ldots, x_\ell) = p(x_{\si(1)}, \ldots, x_{\si(\ell)})$ for
$\si \in \frS_l$. It follows from $\ze(\si \cdot x_i) = \ze(x_{\si(i)} )=
s_{\si(i)}(a) = \si \cdot s_i(a) = \si \cdot \ze(x_i)$ that $\ze$ is an
$\frS_\ell$-module map. It therefore maps the symmetric polynomial algebra
$k[x_1, \ldots, x_\ell]^{\frS_l}$ to $\TS^\ell(A)$. Let $p_i =
\sum_{j=1}^\ell x_j^i$ be the $i^{\rm th}$-power sum. Since $\ze(p_i) =
\sym_\ell(a^i)$ and the $p_1, \ldots, p_\ell$ generate the algebra $\kk[x_1,
\ldots, x_\ell]^{\frS_\ell}$, see e.g.\  \cite[I.2, (2.12)]{Mac95}, our claim
follows. \end{proof}
\sm

We will characterize $\TS^\ell(A)$ by a universal property in
Theorem~\ref{sel-thm}. To do so, we will use the following technical results,
the first of which is standard vector space theory.

\subsection{Lemma}\label{lem:basis-filtration}{\it 
Let $\sfC$ be a basis of a vector space $Y$ and assume that $\sfC=
\bigcup_{i=0}^\ell \sfC_i$ is a partition of $\sfC$. Put $Y_i = \Span (
\sfC_0 \cup \cdots \cup \sfC_i)$ and $Y_{-1} = \{0\}$. Furthermore, suppose
that $\sfB\subseteq Y$ has a partition $\sfB = \bigcup_{i=0}^\ell \sfB_i$
such that the canonical images of $\sfB_i$ and $\sfC_i$ in $Y_i/Y_{i-1}$ are
identical for all $i$, $0\le i \le \ell$. Then $\sfB_0\cup\cdots\cup
\sfB_\ell$
is a basis of $Y$. }

\subsection{Lemma}\label{basis-new}{\it
Let $\{1_A\}\cup \sfB$ be a basis of the vector space $A$, and let $\le $ be
a total order on $\sfB$. Then
\[ \TB = \{ 1_A^{\ot \ell} \} \cup
 \{ \sym_\ell(b_1)\sym_\ell(b_2) \cdots \sym_\ell(b_j) : b_i \in \sfB, b_1 \le b_2 \le \cdots \le b_j, 1 \le j \le \ell\}\]
is a basis of the vector space $\TS^\ell(A)$. } 

\begin{proof} Recall the linear map $\pi \co A^{\ot \ell} \to \sfS^\ell(A)$,
$\pi(a_1 \ot \cdots \ot a_\ell) = a_1 \circ \cdots \circ a_\ell$. Since its
restriction to $\TS^\ell(A)$ is a vector space isomorphism, it suffices to
show that $\pi(\TB)$ is a basis of $\sfS^\ell(A)$. To this end, put
\[
  \sfC_j = \{ 1 \circ \cdots \circ 1 \circ b_1 \circ b_2 \circ \cdots
   \circ b_j : b_1 \le b_2 \le \cdots \le b_j\}, \quad (0 \le j \le \ell). \]
It is well-known that $\sfC=\bigcup_{j=0}^\ell \sfC_j$ is a basis of
$\sfS^\ell(A)$. For $b_1 \le \cdots \le b_j$ we get
\[ \sym_\ell (b_1) \cdots \sym_\ell (b_j) \equiv  \textstyle \frac{l!}{(l-j)!} \,
 1 \circ \cdots \circ 1 \circ b_1 \circ \cdots \circ b_j \, \mod
  S^\ell_{j-1}(A) := \Span_\kk(\sfC_0 \cup \cdots \cup \sfC_j). \]
Applying Lemma~\ref{lem:basis-filtration} (modulo scalars) shows that
$\pi(\TB)$ is indeed a vector space basis of $\sfS^\ell(A)$.
\end{proof}

\subsection{Seligman's symmetric identity}\label{sel-sym}
To a partition $p=(p_1 , \ldots, p_\ell)$ of $\ell$, i.e., $p_i \in \NN$ and
$p_1 + 2 p_2\cdots + \ell p_\ell = \ell$, we associate the conjugacy class
$\scC(p)\subset \frS_\ell$ of permutations whose cycle decomposition consists
of $p_i$ $i$-cycles. For example, $(\ell, 0, \ldots, 0)$ corresponds to (the
conjugacy class of) the identity $1_{\frS_\ell}$ and $(\ell -2, 1, 0, \ldots,
0)$ to the conjugacy class of any transposition. The map $p \mapsto \scC(p)$
is a bijection between the partitions of $\ell$ and the set of conjugacy
classes of $\frS_\ell$. Since the sign of a partition is constant on a
conjugacy class, $\sgn\big( \scC(p)\big)$ is well-defined as the sign of any
partition in $\scC(p)$. Following \cite{Sel80} we say that a linear map $\rho
\co A \to B$ into a unital associative $\kk$-algebra {\em satisfies the
$\ell^{\rm th}$-symmetric identity\/} if
\begin{equation}\label{sel-sym1}
       \sum_{p=(p_1, \ldots, p_\ell)} \, \sgn\big( \scC(p)\big) \, |\big( \scC(p)\big) |\, \rho(a)^{p_1} \, \rho(a^2)^{p_2} \, \cdots \rho(a^\ell)^{p_\ell} = 0
\end{equation}
for all $a\in A$, the sum being taken over all partitions $p$ of $\ell$.\sm

For example, it is easily seen that $\rho = \sym_2$, i.e., $\rho(a) = a \ot 1
+ 1 \ot a \in \TS^2(A)$, satisfies the third symmetric identity:
\[ \sym_2(a)^3 - 3 \, \sym_2(a) \, \sym_2(a^2) + 2\,  \sym_2(a^3) = 0.\]
But even more is true:

\subsection{Proposition (\cite[\S2]{Sel80})}\label{sel-Prop}
{\em The Lie homomorphism $\sym_\ell \co A \to \TS^\ell(A)$ satisfies the
$(\ell +1)^{\rm st}$-symmetric identity and $\sym_\ell(1_A) = \ell \,
1_{\TS^\ell(A)}$.}

In fact, we will show in Theorem~\ref{sel-thm}, $(\TS^\ell(A), \sym_\ell)$ is
universal with respect to these two properties. Our proof is a generalization
of \cite[Prop.~3.1]{Sel80} where this result is proven for finite-dimensional
$A$.

\subsection{Theorem}\label{sel-thm} {\em 
Let $\ell \in \NN_+$. Then for every unital associative $\kk$-algebra $B$ and
every Lie homomorphism $\rho\co A \to B$ satisfying the $(\ell + 1)^{\rm
st}$-symmetric identity and $\rho(1_A) = \ell \, 1_B$ there exists a unique
homomorphism $\vphi \co \TS^\ell(A) \to B$ of unital associative algebras
such that the diagram below is commutative.
\[ \xymatrix{
    A \ar[rr]^{\sym_\ell} \ar[dr]_\rho && \TS^\ell(A) \ar[dl]^\vphi \\ & B}\]
}

\begin{proof}
Let $\U=\U(A^-)$ be the universal enveloping algebra of the Lie algebra
$A^-$, let $\ga \co A^- \to U^-$ be the canonical embedding, and denote by $I
\subset U$ the ideal of the associative algebra $U$ generated by
\begin{enumerate}[(i)]
    \item the elements $\sum_{p=(p_1, \ldots, p_{\ell + 1})} \, \sgn\big(
        \scC(p)\big) \, |\big( \scC(p)\big) |\, \ga(a)^{p_1} \,
        \ga(a^2)^{p_2} \, \cdots \ga(a^{\ell +1})^{p_{\ell+1}}$ used in the
        definition of the $(\ell + 1)^{\rm st}$-symmetric identity, and
    \item $\ga(1_A) = \ell \, 1_U$.
    \end{enumerate}
Let $\can \co \U \to \U/I$ be the canonical quotient map. Then $\psi = \can
\circ \ga \co A \to \U/I$ is a Lie homomorphism satisfying the $(\ell +
1)^{\rm st}$-symmetric identity and $\psi(1_A) = \ell \, 1_{\U}$. Moreover,
since by Proposition~\ref{sel-Prop} the same holds for $(\TS^\ell(A),
\sym_\ell)$, there exists a unique homomorphism $\eta \co \U/I \to
\TS^\ell(A)$ of unital associative algebras such that $\sym_\ell = \eta \circ
\psi$:
\[ \xymatrix{
    A \ar[r]^>>>>>\ga \ar[dr]_{\sym_\ell} \ar@/^1pc/[rr]^\psi  & \U \ar[r]^<<<<<{\can} & \U/I \ar[dl]^\eta\\  & \TS^\ell(A)
}\]
We claim that $\eta$ is an isomorphism. Let $\TB\subset \TS^\ell(A)$ be the
basis of Lemma~\ref{basis-new}. Since $\eta$ is a unital algebra
homomorphism, $\TB$ is the image under $\eta$ of
\[ \sfC = \{ 1_{\U/I} \} \cup \{ \psi(b_1) \cdots \psi(b_j) : 1 \le j \le \ell , b_1 \le \cdots \le b_j\}. \]
Thus, $\eta$ is an isomorphism as soon as we show that $\sfC$ is a spanning
set of $\U/I$. To do so, we use the canonical filtration of $\U$ determined
by the generating set $\ga(A)$ of $\U$. Let $\U_{(t)}$ be the $t^{\rm
th}$-term of this filtration. By the PBW Theorem, $\U_{(\ell+1)}/\U_{(\ell)}
\cong \sfS^{\ell + 1}(A)$. Hence, by Newton's identities (see for example
\cite[III, \S6.1, Rem.~3]{Bou70}), $\U_{(\ell + 1)}/\U_{(\ell)}$ is spanned
by $\ga(a)^{\ell +1 } + \U_{(\ell)}$, $a\in A$. Since
\[
  \sum_{p=(p_1, \ldots, p_{\ell + 1})} \, \sgn\big( \scC(p)\big) \, |\big( \scC(p)\big) |\, \ga(a)^{p_1} \, \ga(a^2)^{p_2} \, \cdots \ga(a^{\ell +1})^{p_{\ell+1}} = \ga(a)^{\ell + 1} + u_{(l)}
\]
for some $u_{(\ell)} \in \U_{(\ell)}$ (the term $\ga(a)^{\ell + 1}$ occurs
for $p=(\ell+1, 0, \ldots, 0)$), it follows that every element in $\U_{(\ell
+ 1)}$ is congruent to $\U_{(\ell)}$ modulo $I$. Since $\U_{(t)} = \U_{(\ell
+ 1)} \U_{(t-\ell -1)}$ for $t> \ell$, this implies that $\U/I$ is spanned by
the image of $\U_{(\ell)}$ under the canonical map $\U\to \U/I$.

We extend the total order $\le$ of $\sfB$ to a total order on $\sfB'=\{1_A\}
\cup \sfB$ by $1_A \le b$ for all $b\in \sfB$. Again by the PBW Theorem (or
one of its corollaries), we then know that $\U/I$ is spanned by $1_{\U/I}$
and elements $\psi(b'_1) \cdots \psi(b'_j)$, $1 \le j \le \ell$, $b'_i \in
\sfB'$, $b'_1 \le \cdots \le b'_j$. Finally, replacing $1_A$ in such an
element by $1_{\U/I}$ using the relation (ii) above, we get that $\sfC$ spans
$\U/I$. \sm

Now let $\rho \co A \to B$ be a map as in the statement of the theorem. Using
the universal property of $\U$, it is immediate that there exists a unique
unital algebra homomorphism $\vphi' \co \U/I \to B$ such that $\rho = \vphi'
\circ \psi$:
\[  \xymatrix{
 A \ar[r]^\psi \ar[dr]_\rho & \U/I \ar[r]^\eta \ar[d]^{\vphi'} & \TS^\ell(A) \ar@{-->}[dl]^{\vphi} \\ & B }
\]
Putting $\vphi = \vphi' \circ \eta^{-1}$ shows $\vphi \circ \sym_\ell =
\vphi' \circ \eta^{-1} \circ \sym_\ell = \vphi' \circ \psi = \rho$.
\end{proof} \sm

A natural question arising at this point is how to check that a given map
$\rho \co A \to B$ satisfies the $k^{\rm th}$-symmetric identity for some
$k\in \NN_+$. Part of this problem is to get a good approach to the left hand
side of \eqref{sel-sym1}. In this paper we will use the following recursion.

\subsection{Recursion \cite[\S III.1, in particular (6)]{Sel80} }\label{rec}
{\em 
Let $\rho \co A \to B$ be a $\kk$-linear map into a unital associative
$\kk$-algebra $B$ such that $[a_1, a_2]=0 \Rightarrow [\rho(a_1), \rho(a_2)]
= 0$ for all $a_i \in A$. Assume that there exists a family $(g_t)_{t\in
\NN}$ of multilinear functions $g_t \co A^t \to B$ satisfying
\begin{enumerate}[\rm (i)]
  \item \label{rec1} $g_1 = \rho$, and

  \item \label{rec2} if $(a_1, \ldots, a_{t+1})$ is a family of commuting
      elements of $A$ then
  \begin{equation}
    \label{rec2a}
    \begin{split}
      g_{t+1}(a_1, \ldots, a_{t+1}) & =  \sum_{j=1}^{t+1} \,
          \rho(a_j) g_t(a, \ldots, \widehat{a_j}, \dots, a_{t+1})
    \\ & \qquad   - 2 \sum_{j<m} g_t(a_j a_m, a_1, \ldots, \widehat{a_j}, \ldots, \widehat{a_m}, \ldots, a_{t+1})
              \end{split}
    \end{equation}
    where $\widehat{\phantom{0}}$ indicates an omitted argument.
\end{enumerate}
Then $g_t(a, a, \ldots, a) $ equals the left hand side of \eqref{sel-sym1}$:$
\begin{equation} \label{rec3}    g_t(a, a, \ldots , a) =   \sum_{p=(p_1, \ldots, p_t)} \, \sgn\big( \scC(p)\big) \, |\big( \scC(p)\big) |\, \rho(a)^{p_1} \, \rho(a^2)^{p_2} \, \cdots \rho(a^t)^{p_t}
\end{equation} }
\sm

Seligman has in fact obtained a closed formula for $g_t(a_1, \ldots, a_t)$ (\cite[Lem.~1.1]{Sel80}).

In this paper we will use a two-pronged approach to verifying that a certain
function satisfies the $(\ell + 1)^{\rm st}$-symmetric identity: We will
define families $(g_t)$ of functions satisfying the recursion \eqref{rec2a}
and we will then know from the context that $g_{\ell + 1} \equiv 0$ (see for
example Proposition~\ref{Lema}).

\subsection{More properties of $\TS^\ell(A)$.}\label{mopro}
In light of the importance of the algebra $\TS^\ell(A)$ for this paper, it is
appropriate to establish more properties of this algebra, besides the ones
given in Lemma~\ref{lem:symm-folklore}. \sm

\begin{inparaenum}[(a)]
\item \label{mopro-a} ({\em Functoriality}) It is immediate that any
    homomorphism $\vphi \co A \to B$ of unital associative $\kk$-algebras
    induces a homomorphism $\TS^\ell(\vphi) \co \TS^\ell(A) \to \TS^\ell(B)$
    and that the assignments $A \mapsto \TS^\ell(A)$ and $\vphi \mapsto
    \TS^\ell(\vphi)$ define a functor $\TS^\ell$ on the category of unital
    associative $\kk$-algebras. It commutes with base field
    extensions,
\[ \TS^\ell(A) \ot K \cong \TS^\ell(A \ot K), \]
and moreover, is an exact functor. \sm

\item\label{mopro-cc} Let $A$ be a finite-dimensional
    semisimple $\kk$-algebra. Then $A^{\ot \ell}$ is semisimple by
    \cite[\S12.7, Cor.~1]{bou:A8} and so is $\TS^\ell(A)$ by
    \cite[Cor.~6]{Monty}. That $\TS^\ell(A)$ is semisimple in case $A$ is
    finite-dimensional central-simple, is also proven in \cite[IV.4]{Sel81}.  \sm

\item \label{mopro-b} ({\em Schur-Weyl duality\/} \cite[p.~467]{Sel80}) In
    particular, for $A= \Mat_d(\kk)$, the associative algebra of $d\times d$ matrices over
    $\kk$, we get
\[ \TS^\ell(\Mat_d(\kk)\big) \cong \big( \End_\kk(V^{\ot \ell})\big)^{\frS_\ell},\]
where the right hand side is the centralizer algebra of $\frS_\ell$ and is
thus described by the classical Schur-Weyl duality. That is, it is a
direct product of $p_d(\ell)$ matrix algebras where $p_d(\ell)$ is the number
of partitions of $\ell$ with at most $d$ parts, see for example
\cite[\S9.1.1]{GoWa}. \sm

\item \label{mopro-d} Assume $A\in \kalg$ is finitely generated. Then so is
    $A^{\ot \ell}$ and hence also $\TS^\ell(A)$, \cite[I, \S1.9, Th.~2]{Bou89}. In
    fact, if $\kk$ is algebraically closed and $X$ is an affine variety with
    coordinate ring $A$, then $\TS^\ell(A)$ is the coordinate ring of the
    symmetric product $X^{(\ell)}$, \cite[Ex.~10.23]{Ha-AG}.

    As a specific example, let $A= k[x]$, the polynomial ring in the
    variable $x$. Then $\TS^\ell(k[x])$ is the ring of symmetric polynomials in $
    \ell$ variables, hence a polynomial ring in the
    elementary symmetric polynomials.
\end{inparaenum}


\vfe

\section{Categories of integrable modules for pairs $(L,\g)$} \label{sec:gen}

Let $\g$ be a finite-dimensional split semisimple Lie algebra over $\kk$. In
this section we consider pairs $(L,\g)$ where
\begin{equation} \label{seti}
 \begin{split}
   & \text{$L$ is a Lie algebra containing $\g$ as a subalgebra such that}
 \\ & \quad \text{$L$ is an integrable module under the adjoint representation of $\g$, cf.\ \ref{int-mod}.} \end{split}
 \end{equation}

\subsection{Notation (standard).} \label{g-not} Let $\h \subset \g$ be a splitting Cartan subalgebra. We denote by

\quad $\De\subset \h^*$ the root system of $(\g, \h)$,

\quad $\g_\al$, $\al \in \De$, the corresponding root space,

\quad $h_\al\in \h$ the unique element in $[\g_\al, \g_{-\al}]$ satisfying
$\al(h_\al) = 2$,

\quad $\euQ = \sum_{\al \in \De} \ZZ \al$ the root lattice of $(\g,\h)$ and

\quad $\euP = \{ \la \in \h^* : \la(h_\al) \in \ZZ \text{ for all } \al \in
\De\}$ the (integral) weight lattice of $(\g, \h)$.\sm

\noindent We choose a root basis $\Pi = \{ \al_1, \ldots, \al_r\}$ of $\De$.
Associated to $\Pi$ are the following data:

\quad $I=\{1, \ldots, r\}$,

\quad $h_i \defeq h_{\al_i}\in \h$, $1\le i \le r$,

\quad $\euQ_+ = \bigoplus_{i\in I} \NN \al_i$ and $\De_+ = \De \cap \euQ_+$,

\quad $\euP_+ = \{ \la \in \h^* : \la(h_i) \in \NN, 1\le i \le r\}\subset
\euP$,

\quad $\varpi_i\in \euP_+$, $i\in I$, fundamental weights defined by
$\varpi_i(h_j) = \de_{ij}$,

\quad $\la \le \mu \iff \mu - \la \in \euQ_+$ for $\la, \mu \in \euP$,

\quad $\frn_+ = \bigoplus_{\al \in \De_+} \g_\al$.\sm

\noindent We choose non-zero $e_i \in \g_{\al_i}$, $1 \le i \le r$, and
denote by $f_i$, $1\le i \le r$, the unique element in $\g_{-\al_i}$
satisfying $[e_i, f_i] = h_i$. \sm

In the following we fix the data and notation introduced above, and use them
without further explanation.

\subsection{Integrable $\g$-modules.} \label{int-mod}
It is standard that for a $\g$-module $V$ the following are equivalent:
\begin{enumerate}[(i)]
  \item $V$ is a sum of finite-dimensional $\g$-modules;

  \item $V$ is a direct sum of finite-dimensional simple $\g$-modules;

   \item for each $v\in V$, the submodule $\U(\g)v$ is finite-dimensional;

   \item for all $i\in I$, the elements $e_i$, $f_i$ act locally nilpotently
       on $V$.
  \end{enumerate}
  The implication (iv)$\Rightarrow$(i) follows from
\cite[Prop.~3.8]{Kac90}.
A $\g$-module $V$ which satisfies the above properties is called an {\em integrable $\g$-module} -- it is sometimes
also called a {\em locally finite\/} $\g$-module.
 We will use the following properties of integrable $\g$-modules:

\begin{inparaenum}
  \item \label{int-mod1}In a short exact sequence $0 \to M' \to M \to M'' \to 0$ of $\g$-modules, $M$ is integrable if and only if $M'$ and $M''$ are integrable. In this case the sequence splits in the category of $\g$-modules.

 \item \label{int-mod2} The tensor product and the direct sum of integrable
  modules is again integrable.

\item \label{int-mod3} The tensor algebra ${\rm T}(L)$ and its quotient
    $\U(L)$ are  integrable $\g$-modules with respect to the canonical
    extension of the adjoint action of $\g$ on $L$ to ${\rm T}(L)$ and
    $\U(L)$ respectively. Denoting the latter action by $\rho_{\U}$, we have
    for $x\in \g$ and $u\in \U(L)$ the formula $\rho_{\U}(x) (u) = xu-ux$
    where the right hand side is given by the multiplication in $\U(L)$.
    \end{inparaenum} \sm

\subsection{Category $\cI(L,\g)$} \label{cat-int-def}
We let $\cI(L,\g)$ denote the category whose objects are $L$-modules which
are integrable as $\g$-modules, and whose morphisms are $L$-module
homomorphisms. Note that in our setting $L$ is an object of $\cI(L,\g)$ with
respect to the adjoint representation.

\subsection{Proposition}\label{indu-prop}
{\it Let $(V, \rho)$ be an integrable $\g$-module. Then
\[ P(V, \rho)\defeq \U(L)\otimes_{\U(\g)} V,\]
considered as $L$-module under left multiplication on the left factor, is a
projective object of $\cI(L,\g)$.
}
\sm

This is a classical result that has been proven in different contexts (e.g.,
\cite{CFK10}, but certainly goes back as far as \cite{Hoc56}). We include a
proof for the sake of completeness.

\begin{proof}
We start by showing that $\U(L) \ot_{\U(\g)} V$ is an object of $\scI(L,\g)$.
To do so, it suffices to prove that $\U(L) \ot_{\U(\g)} V$ is an integrable
$\g$-module with respect to the indicated action. We know from
\ref{int-mod}\eqref{int-mod3} and \ref{int-mod}\eqref{int-mod2} that $\U(L)
\ot_\kk V$ is an integrable $\g$-module with respect to the tensor product
action of $\g$, i.e.,  $x\in \g$ acts on $u\ot v\in \U(L) \ot_\kk V$ by
$\rho_{\U}(x)(u) \ot_\kk v + u \ot_\kk \rho(x)(v) = (xu-ux) \ot_\kk v + u
\ot_\kk \rho(v)$. This action descends to $\U(L) \ot_{\U(\g)} V$, and the
obvious surjective linear map
\[
\pi \co \U(L)\otimes_\kk V\to \U(L)\otimes_{\U(\g)} V\ \ ,\ \
u \otimes_\kk v\mapsto u\otimes_{\U(\g)}v
\]
is a $\g$-module map with respect to the indicated $\g$-actions. It then
follows from \ref{int-mod}\eqref{int-mod1} that $P(V,\rho)$ is indeed an
integrable $\g$-module.

To prove that $P(V,\rho)$ is projective, we must show that
$\Hom_{\U(L)}(P(V), -)$ is an exact functor on $\cI(L,\g)$. Using the
standard adjunction formula we have isomorphisms of abelian groups
\[ \Hom_{\U(L)}(\U(L)\otimes_{\U(\g)} V, M)\cong\Hom_{\U(\g)}(V,\Hom_{\U(L)}(\U(L),M))\cong\Hom_{\U(\g)}(V,M)\]
for any $M\in \Ob\scI(L,\g)$. Let $0\to K\to M\to N\to 0$  be an exact
sequence of objects in $\cI(L,\g)$. By \ref{int-mod}\eqref{int-mod1} this
sequence splits in the category of $\g$-modules. Hence the associated
sequence
 \[ 0\to \Hom_{\U(\g)}(V,K)\to\Hom_{\U(\g)}(V,M)\to\Hom_{\U(\g)}(V,N)\to 0\]
is also exact, in fact split-exact, which implies that $P(V, \rho)$ is
projective.
\end{proof}

\subsection{Corollary}{\it
 The category $\cI(L,\g)$ has enough projectives.}

\begin{proof}
  Indeed, for each $V\in\Ob\cI(L,\g)$ the surjective homomorphism $\epsilon_V \co  P(V)\to V$, $u\otimes v\mapsto u\cdot v$, is an $L$-module map.
\end{proof}

\subsection{Weight modules} \label{wei-mod}
A $\g$-module $M$ is called a \textit{weight module} (with respect to $\h$)
if
\[\ts M=\bigoplus_{\mu\in\h^*} M_\mu,\quad M_\mu=\{m\in M\ :\ h\cdot m=\mu(h)m,\ h\in\h\}.\]
It is equivalent to require $M = \sum_{\mu \in \h^*} M_\mu$. The {\it
weights\/}  of $M$ are those $\mu$ such that $M_\mu\ne 0$.

An integrable $\g$-module is necessarily a weight module with respect to any
splitting Cartan subalgebra of $\g$, hence in particular with respect to our
chosen $\h$. Furthermore, a weight module with a finite set of weights is
integrable.

A \textit{highest weight module of $\g$} of highest weight $\lambda\in\h^*$
is a module $V$ generated by a vector $v_\lambda$ which satisfies $\n_+\cdot
v_\lambda=0$ and $h\cdot v_\lambda=\lambda(h)v_\lambda$ for $h\in\h$. It
follows that $V$ is a weight module, the weights of $V$ lie in
$\lambda-\euQ_+$ and $\dim_\kk V_\mu<\infty$. Moreover, $V_\lambda=\kk
v_\lambda$.

For every $\la\in \euP_+$ there exists a unique (up to isomorphism) highest
weight module $V(\la)$ with highest weight $\la$. It is finite-dimensional
and given as $V(\lambda)=\U(\g)\big/J(\lambda)$, where $J(\lambda)$ is the
left ideal of $\U(\g)$ generated by $I(\lambda)$ and $f_i^{\lambda(h_i)+1}$.
Every irreducible finite-dimensional $\g$-module is isomorphic to a unique
$V(\la)$, $\la \in \euP_+$.

The following elementary lemma is a straightforward generalization
of~\cite[Prop.~3.2]{CFK10}.

\subsection{Lemma}\label{lem:proj-ann}
{\it Let $V$ be a cyclic $\U(\g)$-module, i.e., $V\cong \U(\g)/\Ann_{\U(\g)}
(v_0)$, where $v_0$ is a generator of $V$. Then $\U(L)\otimes_{\U(\g)} V$ is
a cyclic $\U(L)$-module under left multiplication with generator $1\otimes
v_0$, and $\Ann_{\U(L)}(1\otimes v_0)=\U(L)\Ann_{\U(\g)}(v_0)$. In
particular, if $V=V(\lambda)$, then $\U(L)\otimes_{\U(\g)} V(\lambda)$ is
presented, as a $\U(L)$-module, by the relations
	\[ \n_+(1\otimes v_\lambda)=0,\quad h(1\otimes v_\lambda)=\lambda(h)(1\otimes v_\lambda),\quad (f_i)^{\lambda(h_i)+1}(1\otimes v_\lambda)=0,\quad i\in I,\ \ h\in\h.\]
}

\begin{proof}[Proof (sketch).] $\U(L) \ot_{\U(\g)} V$ is clearly cyclic with
generator $1\otimes v_0$, and it is easy to see that
$\Ann_{\U(\g)}(v_0)\subseteq\Ann_{\U(L)}(1\otimes v_0)$, which is enough for
one containment. For the reverse containment, use the PBW Theorem to
decompose $\U(L)$ as a right $\U(\g)$-module.
\end{proof}
	
\vfe


\section{Integrable modules for root-graded Lie algebras}\label{sec:rg}

In this section we will specialize the pairs $(L,\g)$ of the previous section
\ref{sec:gen} assuming that $L$ is a root-graded Lie algebra with grading
subalgebra $\g$. This will allow us to introduce and study global Weyl
modules.

\subsection{Definition (Root-graded Lie algebras)} \label{def:rg} Let $\Theta$ be a
finite root system in the sense of \cite{bou:rac}, hence not necessarily reduced or irreducible, and let $\De\subset \Theta$ be
a subsystem. We will say that $(\Theta, \De)$ is an {\em admissible pair\/} if 
for every connected component $\Theta_c$ of $\Theta$ the following hold:
\begin{enumerate}

\item \label{def:rg0a} if $\Theta_c$ is reduced, then $\Theta_c \subset \De$, and

\item\label{def:rg0b} if $\Theta_c$ is not reduced, say $\Theta_c= \Theta_{c,1} \cup \Theta_{c,2} \cup
\Theta_{c,4}$ is of type $\rmbc_n$ where $\Theta_{c,i}$ is the set of roots of square
length $i$, then  \begin{enumerate}[(I)]
    \item \label{def:rg0bI} either $\De \cap \Theta_c = \Theta_{c,1} \cup \Theta_{c,2}$, hence is of type
$\rmb_n$, or

   \item \label{def:rg0bII} $\De \cap \Theta_c = \Theta_{c,2} \cup \Theta_{c,4}$ and thus is of type
$\rmc_n$. \end{enumerate}

\end{enumerate}
In case \eqref{def:rg0b}  with $n=1$ we have $\Theta_{c,2} = \emptyset$ and $\De \cap \Theta_c$ is of type
$\rmb_1 = \rma_1$ or $\rmc_1=\rma_1$. Observe that $\De$ is reduced and that $\Theta_c \cap \De$ is a connected component of $\De$.
We can and will view $\Theta$ as a set of weights of $\De$. \sm

Let $(\Theta, \De)$ be an admissible pair. A {\em $(\Theta,\De)$-graded Lie algebra\/}
is a pair $(L,\g)$ consisting of a Lie algebra $L$ over $\kk$ and a split semisimple
subalgebra $\g$ satisfying the following conditions:
\begin{inparaenum}[(i)]

\item \label{def:reg1} $\g$ contains a splitting Cartan subalgebra $\h$
    whose root system is (isomorphic to) $\De$ such that
\begin{equation} \label{def:rg1a}
   L = \tsum_{\al  \in \Theta \cup \{0\} } L_\al
   \end{equation}
\qquad where $L_\al = \{l\in L : [h,l] = \al(h) l \text{ for all } h\in \h\}$
for $\al \in \h^*$, and

\item \label{def:reg2} $L_0 = \tsum_{\al \in \Theta} [L_\al, L_{-\al}]$.
 \end{inparaenum}
\sm

\noindent As splitting Cartan subalgebras of $\g$ are conjugate under
automorphisms of $L$ which extend the elementary automorphisms of $\g$
\cite[VIII, \S3.3, Cor.~de la Prop.~10]{Bou75}, condition \eqref{def:reg1}
holds for one splitting Cartan subalgebra if and only if it holds for all. It
is therefore not necessary to include $\h$ in the  notation of a
$(\Theta,\De)$-graded Lie algebra. We will say that
$(L,\g)$ is \emph{$\De$-graded} if $\Theta=\De$ is reduced and that $(L,\g)$ is {\em root-graded\/}
if $(L,\g)$ is $(\Theta, \De)$-graded for some admissible pair $(\Theta, \De)$. The subalgebra $\g$,
called the {\em grading subalgebra\/}, will usually be omitted from the notation.
\sm

Condition \eqref{def:rg1a} says that $L$ is a weight module for $\g$ under
the adjoint action. As its set of weights is finite, it follows that $L$ is
an integrable $\g$-module. Thus $(L,\g)$ satisfies the assumption
\eqref{seti} of the previous section. 

The condition \eqref{def:reg2} serves as a normalization: If $L$ satisfies \eqref{def:reg1}
then $\big( \tsum_{\al \in \De} [L_\al, L_{-\al}]\big) \oplus
\big(\bigoplus_{\al \in \De} L_\al \big) $ is a $\De$-graded ideal of $L$.
From the definition of $L_\al$ it follows that
the decomposition \eqref{def:rg1a} is a grading by the weight lattice $\euP$ of
$\De$, even by the root lattice if $(\Theta,\De)$ does not have a connected component
of type \eqref{def:rg0bII}.
\sm

Root-graded Lie algebras for irreducible $\Th$ are considered in the papers
\cite{abg,abg2,BeSm,beze, BM92,n:3g}. But we caution the reader that in
\cite{abg2,BeSm} the case $(\Theta, \De) = (\rmbc_n, \rmd_n)$ for $n\ge 1$ is
allowed -- a situation not considered here.

\subsection{Examples}\label{rg-exa} \begin{inparaenum}[(a)]
\item \label{reg-exa1} A prime example of a root-graded Lie algebra is the
    Lie algebra $L=\g \ot A$ for any $A\in \kkalg$ and $(\g,\h)$ as in the
    definition above. The Lie algebra product of $L$ is $[l_1 \ot a_1, l_2
    \ot a_2] = [l_1, l_2] \ot a_1a_2$ where $l_i \in L$ and $a_i \in A$. If
    $\g = \bigoplus_{\mu \in \De \cup \{0\}} \g_\mu$ denotes the root space
    decomposition of $\g$ with respect to $\h$, whence in particular $\g_0 =
    \h$ and $\dim \g_\al = 1$ for $\al \in \De$, then the root spaces $L_\al$
    appearing in \eqref{def:rg1a} are the subspaces $L_\al = \g_\al \ot A$.
    We will refer to Lie algebras $\g\ot_\kk A$ as {\em map algebras}, since they can be
    interpreted as regular maps from the affine scheme $\Spec(A)$ to the affine scheme $\g$.
     (Sometimes they are also called \emph{generalized current algebras}).
\sm

\item\label{reg-exa2} A second example, relevant for this paper, is the Lie
    algebra $L=\lsl_n(A)$ for $n\in \NN_+$ and $A$ a unital associative
    $\kk$-algebra, see \ref{subsec:rev-liealg}. This is a $\De$-graded Lie
    algebra for $\De=\rma_{n-1}$. \sm

The two examples \eqref{reg-exa1} and \eqref{reg-exa2} above essentially
describe the structure of $\De$-graded Lie algebras $L$ up to central
extensions, where $\De$ is an irreducible reduced root system of ADE-type.
More precisely, if $\De$ is of type $\rmd$ or $\rme$ then $L$ is a central
extension of the example \eqref{reg-exa1} and if $\De$ is of type $\rma_l$,
$l\ge 3$, then there exists a unital associative $\kk$-algebra $A$ such that
$L/\scC(L) \cong \lsl_{l+1}(A) /\scC\big(\lsl_{l+1}(A)\big)$ for
$\scC(\cdot)$ the centre of the Lie algebra in question \cite{BM92}. We note
that for $\De=\rma_1$ and $\De=\rma_2$ more general Lie algebras than
$\lsl_{l+1}(A)$ occur \cite{beze, BM92, n:3g}, see \eqref{rg-exa-jor}
for $\De=\rma_1$. \sm

\item \label{reg-exa3} A finite-dimensional semisimple Lie algebra $L$, which is {\em isotropic\/}
in the sense that $L$ contains a non-zero split toral subalgebra, is root-graded with
respect to a maximal split toral subalgebra $\h$ and a suitable $\g$ (\cite{Sel76}). We point out that in this
example $\Theta$ need not be reduced. \sm

 \item\label{rg-exa-jor} The Tits-Kantor-Koecher algebra of a unital Jordan
    algebra $J$ and its central coverings are examples of root-graded Lie
    algebras with $\Theta=\De = \rma_1$. In fact, the $A_1$-graded Lie
    algebras are precisely the central coverings of the Tits-Kantor-Koecher
    algebra of a unital Jordan algebra $J$ (\cite[0.8]{beze}, \cite[2.7 and
    3.2(1)]{n:3g}).

\item\label{emas} Let $\frs$ be a finite-dimensional simple Lie algebra  over
    $\kk$, which we assume to be algebraically closed in this example.
    Further, let $\mu$ be a non-trivial diagram automorphism of $\frs$,
    whence $\mu$ has order $|\mu|=2$ or $|\mu|=3$. Let $\ze\in \kk$ be a primitive $|\mu|$th-root
    of unity, put $\frs^{(i)} = \{ x\in \frs : \mu(x) = \ze^i x\}$, so that $\frs$ has
    an eigenspace decomposition
\[ \frs = \frs^{(0)} \oplus \frs^{(1)}\quad \hbox{respectively}\quad
 \frs = \frs^{(0)} \oplus
\frs^{(1)} \oplus \frs^{(2)},
\]
which is a $(\ZZ/|\mu|\ZZ)$-grading of $\frs$. In particular,
$\frs^{(0)}=:\g$ is a subalgebra of $\frs$ and $\frs^{(i)}$, $i>0$, are
$\g$--modules under the adjoint action of $\frs^{(0)}$.

One knows (see for example  \cite[Prop.~7.9 and 7.10]{Kac90}) that $\g$ is a
simple subalgebra. Let $\h$ be a Cartan subalgebra of $\g$, and denote by
$\De$ the root system of $(\g,\h)$. It is further known that $\De$ is not
simply-laced and that the $\g$--modules $\frs^{(i)}$, $i>0$, are irreducible
of highest weight equal to the highest short root of $\De$ or  twice that
root. For $\al \in \h^*$ let $\frs_\al$ be the weight space of $\frs$ with
respect to $\frh$. We then have a weight decomposition
\begin{equation*} \label{rg-exa1}
   \frs = \textstyle \bigoplus_{\al \in \Theta\cup \{0\}} \, \frs_\al,
\end{equation*}
for $\Theta=\De$ or $\Theta = \rmbc_r\supset \De=\rmb_r$. Thus $(\Theta,
\De)$ is an admissible pair. In fact, {\em the pair $(\frs, \g)$ is a
$(\Theta, \De)$-graded Lie algebra\/} since the condition \eqref{def:reg2} of
\ref{def:rg} holds by simplicity of $\frs$. We note that every weight space
$\frs_\al$ respects the eigenspace decomposition of $\frs$ with respect to
$\mu$:
$\frs_\al = \textstyle \bigoplus_{0 \le i \le |\mu|}\,
               \frs^{(i)}_\al$ with $\frs^{(i)}_\al = \frs^{(i)} \cap \frs_\al$.
\sm

Generalizing this example, let $A\in \kalg$ and assume that $\mu$ also acts
on $A$ by an algebra automorphism. Let $L$ be the associated equivariant map
algebra, i.e., the fixed point subalgebra of the Lie algebra $\frs \ot_\kk A$
under the diagonal action of $\mu$ (\cite{NSS12}).
For $A$
finitely generated, these Lie algebras  are the main object of investigation
in \cite{FMSb} (the setting in this example is assumed from section 5 on in
\cite{FMSb}).
Using obvious notation, we
therefore get
\[
L = \textstyle \bigoplus_{0 \le i < |\mu|} \; \frs^{(i)} \ot A^{(-i)}.
\]
We set $\g$ equal to the subalgebra $\frs^{(0)} \ot \kk 1_A$ of $L$. Since $\h
\subset \g$ acts on the factor $\frs$ of $\frs \ot_\kk A$, it follows that
$L$ has a weight space decomposition
\begin{equation*}
  \label{rg-exa3} L = \textstyle \bigoplus_{\al \in \Theta \cup \{0\}} \, L_\al,
 \qquad L_\al = \bigoplus_{0 \le i < |\mu|} \, \frs_\al^{(i)} \ot A^{(-i)}.
\end{equation*}
We now sketch the argument that
condition
\eqref{def:reg2} of \ref{def:rg} holds, and therefore
{\em $(L,\g)$ is a root-graded Lie algebra.}
To this end, it is enough to verify that
\[
\frs_0^{(i)}\
\otimes A^{(-i)}=\sum_{\alpha\in\Theta}
\left[\frs_\alpha^{(0)}\otimes A^{(0)},\frs_{-\alpha}^{(i)}\otimes A^{(-i)}\right]
\text{ for }0\leq i<|\mu|.\]
The latter relation will follow immediately if we can show that
\begin{equation}
\label{s0(i)-Hadi}
\frs_0^{(i)}=\sum_{\alpha\in\Theta}\left[
\frs_\alpha^{(0)},\frs_{-\alpha}^{(i)}\right]
\text{ for }0\leq i<|\mu|.
\end{equation}
To see why \eqref{s0(i)-Hadi} is true, first note
that
by simplicity of $\frs$,
and the fact that $\frs_0^{(0)}$ is the Cartan subalgebra of $\frs^{(0)}$,
we have
\[\frs_0^{(0)}=\sum_{\alpha\in\Delta}
\left[\frs_\alpha^{(0)},\frs_{-\alpha}^{(0)}\right].
\]
This immediately implies \eqref{s0(i)-Hadi} for $i=0$.
 Next assume that $0<i<|\mu|$.
We can also assume that $\frs_0^{(i)}\neq\{0\}$.
Then $\frs^{(i)}$ is  a nontrivial irreducible $\frs^{(0)}$-module, and therefore it is generated by a nonzero vector $v_\beta\in \frs_{\beta}^{(i)}$, for some $\beta\neq 0$. Thus every nonzero vector in
$\frs_0^{(i)}$ is a linear combination of vectors of the form
$x_{\beta_1}\cdots x_{\beta_j}v_\beta$, for some $j>0$ and $\beta_1,\ldots,\beta_j\in\Delta\cup\{0\}$. Since $\frs_0^{(0)}$ is the Cartan subalgebra of $\frs^{(0)}$, we can assume that $ \beta_1\neq 0$. Consequently, $x_{\beta_1}\cdots x_{\beta_j}v_\beta
\in\left[
\frs_{\beta_1}^{(0)},\frs_{-\beta_1}^{(i)}\right]$.
This completes the proof of
\eqref{s0(i)-Hadi}.
\ms

\item\label{rg-exa-tori} {\em  Lie tori  are    examples    of
    root-graded     Lie    algebras} \cite{n:tori, n:eala} (these are the  cores and
    centreless     cores    of    extended     affine    Lie    algebras). In
    particular, if $\hat \g$ denotes an affine Kac-Moody Lie algebra then
    $\tilde \g = [ \hat \g, \hat \g]$ and the (twisted or untwisted) loop algebra
    $\tilde \g / \scC(\tilde \g)$ is a root-graded Lie algebra.
\ms

\item\label{rg-exa-hyper} Let $L$ be a $(\Theta, \De)$-graded Lie algebra
    with grading subalgebra $\g$ and let $L'$ be a subalgebra of $L$
    containing $\g$. Then $L'$ has a weight space decomposition $L'=\sum_{\al
    \in \Theta \cup \{0\}}\, L'_\al$ with $L'_\al = L' \cap L_\al$. It
    follows that $L'$ is root-graded as soon as condition \eqref{def:reg2} of
    \ref{def:rg} is fulfilled. We give examples of this situation: \sm

\begin{inparaenum}[(1)]

 \item\label{rg-exa-hyphyp}  The hyperspecial current algebra $\frC\g$ of \cite{CIK} is  $(\Theta,\De)$-graded for $(\Theta,\De)$ of type
    \ref{def:rg}\eqref{def:rg0bII}. It is a subalgebra of the twisted loop
    algebra of type $\rma_{2r}^{(2)}$ which is root-graded by \eqref{rg-exa-tori}.
Using the notation of \cite{CIK}, there exists an obvious epimorphism $\pi\co
\hat \h^* \to \h^*$ with
    $\delta\mapsto 0$. It follows that $\pi\big( \hat R_{\rm re}(+)\cup \NN
    \delta \cup \hat R_{\rm re}(-)\big) = R \cup \{0\} \cup \frac{1}{2}
    R_\ell \cong \rmbc_r$. The subalgebra $\g \subset \mathfrak{C} \g = L$ is
    simple of type $R=\rmc_r$. The Cartan subalgebra $\h$ of
    $\g$ acts diagonalizably on $L$ with weights $\al$ in $\rmbc_n\cup
    \{0\}$, namely $L_\al = \bigoplus_{\pi(\be) = \al} \hat \g_\be$.
    \sm

 \item The standard current algebras considered in \cite[\S9]{CIK} are
 root-graded. The details are similar to \eqref{rg-exa-hyphyp}. As a special case,
 the twisted current algebras in the sense of \cite{FK12} are $(\Theta,
    \De)$-graded with $(\Theta,\De)$ of type \ref{def:rg}\eqref{def:rg0bI}.
 \end{inparaenum}
\end{inparaenum}

The following facts about root-graded Lie algebras will be useful later.

\subsection{Lemma}\label{rogrfa}{\em
Let $(L,\g)$ be a $(\Theta,\De)$-graded Lie algebra. We use the notation of\/ {\rm \ref{g-not}}
and {\rm \ref{def:rg}}. \sm

\begin{inparaenum}[\rm (a)]

\item \label{rogrfa-a} For $\al \in \De$ let $(e_\al, h_\al, f_\al) \in
 \g_\al \times \h \times \g_{-\al}$ be an $\lsl_2$-triple. Then $[e_\al, L_\be] = L_{\al + \be}$
 whenever $\beta \in \Theta\cup \{0\}$ satisfies $\lan \be, \al\ch\ran \ge -1$. In particular,
\begin{equation}
    [e_\al, L_0] = L_\al \quad \hbox{and} \quad [e_\al, L_\al] = L_{2\al}. \label{rogrfa1}
\end{equation}

\item \label{rogrfa-b}  Let $\frC$ be the set of connected components
    $\Th_c$ of $\Th$ such that $(\Th_c, \De_c)\cong (\rmbc_n, \rmc_n)$ for some $n=n(c)$ depending on
    $c\in \frC$.
   For each $c\in \frC$ let $\be_c\in \Pi= \{\al_1, \ldots, \al_r\}$ be the unique root with $\be_c/2 \in \Th$.
   Then the $0$-weight space $L_0$ satisfies
\begin{equation} \label{eq-hadi01} \textstyle
L_0 = \sum_{i=1}^r \, [L_{\alpha_i},L_{-\alpha_i}] + \sum_{c\in \frC} \,
      [L_{\be_c/2}, L_{-\be_c/2}].
\end{equation}
 In particular, if $(\Theta, \De)$ does not have a connected component
of type {\rm \ref{def:rg}}\eqref{def:rg0bII} then
\begin{equation} \label{eq-hadi0} \textstyle
L_0 = \sum_{i=1}^r [L_{\alpha_i},L_{-\alpha_i}].
\end{equation}

\item \label{rogrfa-c} If $I \ideal L$ is an ideal of $L$ with $\g \subset
    I$, then $I=L$.
\end{inparaenum}
}

\begin{proof} \eqref{rogrfa-a} The subspace $\bigoplus_{n\in \ZZ} L_{\be + n \al}$ is
a weight module for the adjoint action of the subalgebra $\kk e_\al \oplus \kk h_\al \oplus f_\al
\cong \lsl_2(\kk)$. The result therefore follows from standard $\lsl_2$-representation theory.

\sm

\eqref{rogrfa-b}
 Let $\ga \in \Th \cap \euQ_+$ and assume that $\ga$ is not a simple root.  It is then
well-known that there exists a simple root $\al_i\in \Pi$ satisfying $\lan
\ga, \al_i \ch\ran >0$ and therefore $\beta = \ga - \al_i \in \Theta$.

 Let $(e_i, h_i,
f_i) \in \g_{\al_i} \times \h\times \g_{-\al_i}$ be an $\lsl_2$-triple. Then
$L_\ga = [e_i, L_\be]$ by \eqref{rogrfa-a}. Hence
\begin{equation}\label{eq-hadiLaLbLc} [L_\gamma,L_{-\gamma}]=
[[e_i,L_\beta],L_{-\gamma}] \subset
[e_i,[L_\beta,L_{-\gamma}]]+[[e_i,L_{-\gamma}],L_\beta] \subseteq
  [L_{\al_i}, L_{-\al_i}] + [L_\be, L_{-\be}].
\end{equation}
By induction on the height we get $[L_\ga, L_{-\ga}] \subset \sum_{i=1}^r
[L_{\alpha_i},L_{-\alpha_i}]$, in particular \eqref{eq-hadi0} follows when
$\frC=\varnothing$. In general,
\[ 
\textstyle
L_0 = \sum_{i=1}^r [L_{\alpha_i},L_{-\alpha_i}] +
         \sum_{c\in \frC} \sum_{\ga \in \Theta_c \setminus \De_c}
                  [L_\ga, L_{-\ga}].
\]
Fix $c\in \frC$. To simplify notation, we write $\Pi \cap \De_c = \{ \eps_1 -
\eps_2, \ldots, \eps_{n-1} - \eps_n, 2\eps_n\}$ so that $\be_c = 2\eps_n$. To
finish the proof of \eqref{eq-hadi01} we need to show $[L_{\eps_i},
L_{-\eps_i}] \subset \sum_{i=1}^r [L_{\alpha_i},L_{-\alpha_i}] + [L_{\eps_n},
L_{-\eps_n}]$. But this follows from \eqref{eq-hadiLaLbLc} for $\ga =
\eps_i$, $\al_i= \eps_i - \eps_{i+1}\in \Pi$, $\be = \eps_{i+1}= \ga - \al_i$
together with a downward induction starting at $n$.  \sm

\eqref{rogrfa-c} It follows from \eqref{rogrfa1} that $L_\al \oplus L_{2\al} \subset I$
for all $\al \in \De$. If there exists $\be \in \Theta$ with $\al = 2\be \in \De$ (case
\ref{def:rg}\eqref{def:rg0bII}) then $[e_\al, L_{-\be}] = L_\be \subset I$ by \eqref{rogrfa-a}, whence
$L_\ga \subset I$ for all $\ga \in \Theta$. Now
the relation
$L_0=\sum_{\alpha\in\Theta}[L_\alpha,L_{-\alpha}]$
proves $L_0 \subset I$, thus $L=I$.
\end{proof}

\subsection{The category $\cI(L,\g)^\lambda$}
\label{def:cat-I-bd} This category has been introduced in different settings
in several recent papers, e.g.\ \cite{CFK10,CIK,FK12,FMS11,FMSb}. We do this
here for a $(\Theta,\De)$-graded Lie algebra $(L,\g)$ which  covers all
previous papers apart from some preliminary results in the first four
sections of \cite{FMSb}. But note that the examples in \ref{rg-exa} show that
our setting applies to many more types of Lie algebras not covered in
\cite{FMSb}.  We use the notation introduced in
 \ref{g-not}. \sm

Let $\Theta_+\subset \Theta$ be a positive system, i.e., $(\Th_+ + \Th_+)\cap \Th \subset \Th_+$, $\Th_+ \cap (-\Th_+)
= \emptyset$ and $\Th = \Th_+ \cup (-\Th_+)$. We can and will assume that $\Theta_+ \cap \De = \De_+$.
We denote by $\NN[\Th_+]$ the monoid spanned by
$\Th_+$ in the weight lattice $\euP$ of $\De$. Then $\NN[\Th_+] \supset \euQ_+$ and
$\NN[\Th_+] / \euQ_+$ is an abelian group of order $2^s$ where $s$ is the number of irreducible
components of $(\Th, \De)$ of type \ref{def:rg}\eqref{def:rg0bII}. For later use we note
\begin{equation}
  \label{def:cat-I-bd-0} (-\Th_+) \cap \NN[\Th_+] = \emptyset.
\end{equation}
 It will be convenient to put
\[ \mu \preceq \la \quad \iff \quad \la - \mu \in \NN[\Th_+].\]
Since $\euQ_+ \subset \NN[\Th_+]$ we have $\mu \le \la \implies \mu \preceq \la$.
We define subalgebras
\[ L_{\pm}  \defeq \textstyle \bigoplus_{\al \in \pm \Theta_+} L_\al\]
of $L$. Thus
\begin{equation} \label{def:cat-I-bd-00}
    L=L_- \oplus L_0 \oplus L_+
\end{equation} is a decomposition
of $L$ as a direct sum of subalgebras, but we point out that this is in
general not a ``triangular decomposition'' in the sense of \cite{MP95} because for instance $L_0$ is not abelian in general.
\medskip

For $\lambda\in \euP_+$, let $\cI(L,\g)^\lambda$ denote
the full subcategory of $\cI(L,\g)$ whose objects $V$ satisfy
$\wt(V)\subseteq \lambda-\NN[\Th_+]$, where $\wt(V)$ denotes the weights of the
$\g$-module $V$. We note some elementary properties of the category
$\scI(L,\g)^\la$. \sm
\begin{inparaenum}[(a)]

\item\label{def:cat-I-bd1} Every $V\in \Ob\scI(L,\g)$ has a unique decomposition
$V=\bigoplus_{\mu \in \euP_+} V_{(\mu)}$ as $\g$-module where $V_{(\mu)}$
is the isotypic component of $V$ associated to the simple $\g$-module
$V(\mu)$. We have $V\in \Ob\scI(L,\g)^\la \Iff V = \bigoplus_{\mu \in \euP_+, \, \mu\preceq\la }
V_{(\mu)}$. By standard $\lsl_2$-theory, if $\lambda=\sum_{i\in I}\ell_i\varpi_i$, then
\begin{equation}\label{def:cat-I-bd1a}
    \la - (\ell_i +1) \al_i \not\in \wt V \text{ for every } V\in \Ob\scI(L,\g)^\la
   \text{ and every }i\in I.
\end{equation}
Since there are only finitely many $\mu \in \euP_+$ with $\mu \in \la - \NN[\Th_+]$,
the set of weights of any $V\in \Ob\scI(L,\g)^\la$ is finite. \sm

\item There is a natural truncation functor $\bfT^\la \co
    \cI(L,\g)\to\cI(L,\g)^\lambda$ given on objects by
\[\ts V\mapsto V^\lambda\defeq V\big/ \sum_{\mu\not\preceq \la} \U(L)\cdot V_\mu \]
and on morphisms by $\vphi \mapsto \vphi^\la$, where for an $L$-module
homomorphism $\vphi \co V \to N$ the map $\vphi^\la$ is the induced quotient
map keeping in mind that $\vphi(  \sum_{\mu \not\preceq \la} \U(L)\cdot V_\mu)
\subset  \sum_{\mu \not\preceq \la} \U(L)\cdot N_\mu$. An object $V\in
\scI(L,\g)$ lies in $\scI(L,\g)^\la$ if and only if $V=V^\la$.

The functor $\bfT^\la$ is exact and preserves projective objects. In
particular, in view of Proposition~\ref{indu-prop} we have (cf.~\cite[3.2,
Cor.~1]{CFK10})
\begin{equation}
  \label{def-cat-I-bd-p}
     \text{\it $P(V,\rho)^\la$ is projective for any $(V,\rho) \in \scI(L,\g)$.}
\end{equation}

\item\label{def:cat-I-bd2} If $\dim V_\la = 1$ and $V= \U(L) V_\la$, then $V$
    has a unique maximal $L$-submodule and hence a unique irreducible
    quotient.
\end{inparaenum}

\subsection{Remark}\label{supp-lem}{\em
For $V\in \scI(L,\g)^\la$ denote by $\wt(V)$ the set of weights of $V$. Then
$\wt(V) \subset \wt(V(\la))$ whenever one of the following holds,
\begin{enumerate}[\rm (i)]
  \item \label{supp-lem1} $\NN[\Th_+] = \euQ_+$, i.e, only the cases {\rm \ref{def:rg}}\eqref{def:rg0a} and
    {\rm \ref{def:rg}}\eqref{def:rg0bI} occur,

   \item \label{supp-lem2} $(\Th, \De) \cong (\rmbc_r, \rmc_r)$, $\la = \sum_{i\in I} \ell_i \varpi_i$ with
   $\ell_r >0$, where $\varpi_r$ is the fundamental weight corresponding to the unique long root  in $\Pi$.
\end{enumerate}
}

\subsection{Examples}\label{ex:int-fd} \begin{inparaenum}[(a)]
\item\label{def:cat-I-d} ($\la=0$)  By \ref{def:cat-I-bd}\eqref{def:cat-I-bd1}, an
    object $V$ in $\scI(L,\g)^0$ has the property that $\g\cdot V = 0$. It
    then follows from Lemma~\ref{rogrfa}\eqref{rogrfa-c} that $\Ann_L(V) =
    L$, i.e., $L \cdot V = 0$. \sm

\item\label{ex:int-fd-b} (\cite[I, \S2, Prop.~1.1]{Sel88} for
 Example~\ref{rg-exa}\eqref{reg-exa3}) {\em     For     every     finite-dimensional
 irreducible     $L$-module    $V$     there
 exists a unique $\la \in \euP_+$ such that $V\in \Ob\scI(L,\g)^\la$.}
Indeed, such a $V$ is an integrable $\g$-module by \ref{int-mod} with
finitely many weights. Since $L_\al \cdot V_\mu \subset V_{\al + \mu}$ there
exists a weight $\la \in \euP_+$ of $V$ such that $L_\al \cdot V_\la = 0$ for
all $\al \in \Theta_+$. One checks that $V' = \Span \{l_1 \cdots l_k \cdot v :
l_i \in L_{-\be_i}, \be_i \in \Th_+, v\in V_\la\}$ is an $L$-submodule, whence
$V'=V$. By construction the weights $\mu$ of $V'$ lie in $\la- \NN[\Th_+]$. This
also shows uniqueness of $\la$.
\sm

\item If $\Th $ is irreducible and $\theta$ the highest root with respect to $\Theta_+$,
we have $L \in \scI(L,\g)^\theta$.
\end{inparaenum}

\subsection{Global Weyl modules}\label{def:global-weyl}
For $\lambda\in \euP_+$ the \textit{global Weyl module (of highest weight
$\la$)\/} is defined as
\[ W(\lambda) \defeq  P(V(\lambda))^\lambda = \big(\U(L)\ot_{\U(\g)} V(\la)\big)^\la  \in\Ob\cI(L,\g)^\lambda.\]
By Lemma~\ref{lem:proj-ann} we see that $W(\lambda)$ is a cyclic
$\U(L)$-module, generated by the image of $1\otimes v_\lambda \in \U(L)
\ot_{U(\g)} V(\la)$ in $W(\la)$, which we will denote by $w_\lambda$, and
from \eqref{def-cat-I-bd-p} we know that
$W(\la)$ is projective
in $\scI(L,\g)^\la$.\sm

The terminology ``global Weyl module'' first appeared in this form in
\cite{CFK10} in the context $L=\g\otimes A$ for $A\in \kkalg$, but these
modules were first defined in \cite{CP01}. The name is justified in view of
their universal property \ref{prop:global-weyl-def-rel}\eqref{prop:gwdrb} and
the analogous universal property of Weyl modules in the theory of algebraic
groups \cite[II, 2.13(b)]{Jan:AG}.
For more background information on the significance of global Weyl modules and the choice of nomenclature we refer the reader to the papers by Chari and her collaborators that are listed in our bibliography.

\subsection{Proposition (Characterization of global Weyl modules, \cite[Prop.~4]{CFK10})}
\label{prop:global-weyl-def-rel}
{\it 
\begin{inparaenum}[\rm (a)]

\item \label{prop:global-weyl-def-rel-a} The global Weyl module $W(\lambda)$ is
    generated by $w_\lambda$ with defining relations 	\begin{equation}
    \label{prop:gwdr1}
  L_+\cdot w_\lambda=0,\qquad (h -\lambda(h)1)\cdot w_\lambda=0
   \quad (h\in \h), \qquad f_i^{\lambda(h_i)+1}\cdot w_\lambda=0, \ \ i\in I.\end{equation}
In other words, $\Ann_{\U(L)}(w_\lambda)$ is the left ideal of $\U(L)$
generated by the set
\begin{equation}\label{ann-w-lambda}
 L_+\cup\left\{f_i^{\lambda(h_i)+1}\ :\ i\in I\right\}\cup \left\{h-\lambda(h)\ :\ h\in\h\right\}.\end{equation}

\item \label{prop:gwdrb}  For every cyclic module $V \in \Ob \scI(L,\g)^\la$
    with generator $v_\la \in V_\la$ there exists a unique surjective
    $L$-module homomorphism $W(\la) \to V$ sending $w_\la$ to $v_\la$.
\end{inparaenum}}

\begin{proof} The proof is essentially the same as in \cite[Prop.~4]{CFK10}.
We give a sketch for the convenience of the reader.

(I) Let $V\in \Ob \scI(L,\g)^\la$ be a cyclic module as in
\eqref{prop:gwdrb}. Then $L_+ \cdot V_\la \subset \sum_{\mu \not\preceq \la}
V_\mu$ by \eqref{def:cat-I-bd-0}, which forces $L_+ \cdot v_\la = 0$.
Obviously, $h \cdot v_\la = \la(h) v_\la $ for all $h\in \h$. Finally,
\eqref{def:cat-I-bd1a} implies $f_i^{\la(h_i) + 1} \cdot v_\la = 0$.

(II) Proof of \eqref{prop:global-weyl-def-rel-a}: We have already noted that $W(\la)$ is a cyclic module. It
thus follows from (I) that $(W(\la), w_\la)$ satisfies the relations
\eqref{prop:gwdr1}. Let $W$ be the cyclic $L$-module generated by a vector
$w$ with the given relations. Since the relations hold in $W(\lambda)$, we
get that $W(\lambda)$ is a quotient of $W$. Since $\n_+\subseteq L_+$, it
follows from the presentation given in Lemma~\ref{lem:proj-ann} that $W$ is a
quotient of $P(V(\lambda))$. The relations imply that $W$ is a weight module,
and from $L_+\cdot w=0$, we see that the weights of $W$ lie in
$\lambda-\NN[\Th_+]$. Hence, $W$ is a quotient of
$W(\lambda)=P(V(\lambda))\big/\sum_{\mu\not\in \la - \NN[\Th_+]}\U(L)\cdot
P(V(\lambda))_\mu,$ so that $W\cong W(\lambda)$.

(III) Proof of \eqref{prop:gwdrb}: This follows from (I) and \eqref{prop:global-weyl-def-rel-a}. \end{proof}

\medskip

In the remainder of this section we relate the representation theory of $L$ with
that of a certain associative algebra. We start with a generalization
of~\cite[\S3.4]{CFK10} using essentially the same proof.

\subsection{Lemma} \label{right-action-L0}
{\it The formula
	\[ (u\cdot w_\lambda)\cdot a\defeq ua\cdot w_\lambda,\quad u\in \U(L),\quad a\in \U(L_0),\]
	endows $W(\lambda)$ with the structure of a $(\U(L),\U(L_0))$-bimodule.
}

\begin{proof}
	Let $J=\Ann_{\U(L)}(w_\lambda)$, the left ideal of $\U(L)$ described in~\eqref{ann-w-lambda}. Then the given action is well-defined if and only if $JL_0\subseteq J$. To see that this holds, it suffices to show that $ga\in J$ for $g$ one of the generators of $J$ described in Proposition~\ref{prop:global-weyl-def-rel} and $a\in L_0$.
	\begin{enumerate}[(i)]
	\item To see that this is true for $g=\ell_+\in L_+$, note that $\ell_+ a=[\ell_+, a]+a\ell_+\in L_++L_0L_+\subseteq J$, since $[L_\alpha, L_\beta]\subseteq L_{\alpha+\beta}$, and thus $[L_0,L_+]\subseteq L_+.$
	\item $(h-\lambda(h)1)a=a(h-\lambda(h)1)\in J$ since $[L_0,\h]=0$.
	\item In $W(\lambda)$ we have $\n_+a\cdot w_\lambda=[\n_+,a]\cdot w_\lambda+a\n_+\cdot w_\lambda=0$, because $\n_+$ and $[\n_+,a]$ lie in $L_+$. Since $(h-\lambda(h)1)a\cdot w_\lambda=0$ for $h\in\h$, either $a\cdot w_\lambda=0$, or $a\cdot w_\lambda$ is a primitive vector for the $\g$-action on $W(\lambda)$. In the latter case, since $W(\lambda)$ is a completely reducible $\g$-module, we have $U(\g)(a\cdot w_\lambda)\cong V(\lambda)$. It follows that for $1\le i\le n$ we have $f_i^{\lambda(\alpha_i^\vee)+1} a\cdot w_\lambda=0$, so that $f_i^{\lambda(\alpha_i^\vee)+1}a\in J$.
	\end{enumerate}
	This shows that the given action is well-defined. Since it clearly commutes with the left action of $\U(L)$, the proof is complete. \end{proof}


\subsection{The Chari--Pressley algebra $\mathbf A_\la$} \label{def:cha}
We have
 \begin{align*} \Ann_{\U(L_0)}(w_\lambda)&=\{x\in\U(L_0)\ :\ w_\lambda\cdot x=0\}\\
 &=\{x\in\U(L_0)\ :\ W(\lambda)\cdot x=0\}=\Ann_{\U(L_0)}(W(\lambda)),\end{align*}
 a two-sided ideal of $\U(L_0)$. Hence
 \[ \mathbf A_\lambda\defeq \U(L_0)\big/\Ann_{\U(L_0)}(w_\lambda)\]
 is an associative $\kk$-algebra,  called the Chari-Pressley algebra for now, but see Theorem~\ref{thm:chari=seli}.

 Let $\mathbf
A_\lambda\text{-mod}$ denote the category of left $\mathbf
A_\lambda$-modules.

\subsection{Lemma (\cite[Lem.~4]{CFK10})}
{\it For $\lambda\in \euP_+$ and $V\in\Ob\cI(L,\g)^\lambda$ we have
$\Ann_{\U(L_0)}(w_\lambda) \cdot V_\lambda=0.$ In particular, $V_\la$ is a left
$\mathbf A_\lambda$-module.
}

\begin{proof}
 Let $v\in V_\lambda$. Then $v$ satisfies the relations of Proposition~\ref{prop:global-weyl-def-rel}, so there is a homomorphism of $L$-modules $\pi: W(\lambda)\to \U(L)\cdot v$ which maps $w_\lambda\mapsto v$. In particular, if $u\in\Ann_{\U(L_0)}(w_\lambda)$ then $0=\pi(u\cdot w_\lambda)=u\cdot v$.
 \end{proof}
\vfe

\section{Admissible modules, Seligman algebras}\label{sec:admiss}

In this section we continue with the assumptions of  Section
\ref{sec:rg}: $L$ is a $(\Th,\De)$-graded Lie algebra with grading subalgebra
$\g$ and spliting Cartan subalgebra $\h$.

In this section we study in more detail the highest weight space $V_\lambda$
of modules $V\in\Ob\cI(L,\g)^\lambda$.
Our main goal is to relate the module categories $\cI(L,\g)^\lambda$ and $\mathbf A_\lambda\text{-mod}$ to the work of Seligman (\cite{Sel81,Sel88}) on rational modules of simple Lie algebras.

\subsection{The Harish-Chandra homomorphism} \label{lem:HC-hom}
We know from \ref{int-mod}\eqref{int-mod3} that $\U(L)$ is an integrable
$\g$-module with respect to the action $\rho_{\U}$, whence in particular a
weight module: $\U(L)=\bigoplus_{\mu\in \euQ}\U(L)_\mu$ with weight spaces
$\U(L)_\mu$.  Clearly $\U(L)_0$ is a subalgebra of $\U(L)$ containing the
centre of $\U(L)$. By a standard consequence of the PBW Theorem we have \sm

\begin{inparaenum}
\item $K=\big( L_-\U(L)\big)\cap \U(L)_0=\big(\U(L)L_+\big)\cap \U(L)_0$ is a
    two-sided     ideal     of     $\U(L)_0$.

\item\label{lem:HC_hom2} $\U(L)_0=\U(L_0)\oplus K$ and the associated
    projection $\pi_0: \U(L)_0\to \U(L_0)$ with $\Ker(\pi_0) = K$ is a homomorphism of algebras.
\end{inparaenum}
\noindent Following the standard terminology for the case $L=\g$, we call
$\pi_0$ the {\em Harish-Chandra homomorphism}.

\subsection{Definition (admissible modules, Seligman algebra)}\label{def:adm}
Fix $\lambda=\sum_{i\in I} \ell_i \omega_i\in \euP_+$. \sm

\begin{inparaenum}[(a)]

\item We say an $L_0$-module $M\ne 0$ has \textit{weight} $\lambda$ if $h\cdot
    m=\lambda(h)m$ for all $h\in\h$ and $m\in M$. \sm

\item An $L_0$-module $M$ of weight $\lambda$ is called
    $\lambda$-\textit{admissible} if for all $i\in I $ and for any pair of
    sequences
	\begin{align}
     \label{def:adm1}
(x_1,\ldots,x_{m_i}), \quad x_j\in L_{\be_j}, \be_j\in \Th_+,
\quad (y_1,\ldots,y_{\ell_i + 1}) ,\quad
    y_j\in L_{-\al_i}, \quad \textstyle \sum_{j=1}^{m_i} \be_j = (\ell_i + 1) \al_i,
    \end{align}
we have
	\[ \pi_0(x_1x_2\cdots x_{m_i}y_1y_2\cdots y_{\ell_i + 1})\cdot M=0.\]
	An $L_0$-module will be called \textit{admissible} if it is
$\lambda$-admissible for some $\lambda\in \euP_+$. If $\Th = \De$ is reduced
the condition $\sum_{j=1}^{m_i} \be_j = (\ell_i + 1) \al_i$ is equivalent to
$m_i = \ell_i + 1$ and all $\be_j = \al_i$. If $\Th$ is not reduced and
$s\al_i \in \Th$ with $s=2$ or $s=\frac{1}{2}$, then  $\be_j \in \{\al_i,
s\al_i\}$ follows, but not necessarily $\be_j = \al_i$. \sm

\item \label{def:adm3} Let $J^\la$ denote the two-sided ideal of the
    associative algebra $\U(L_0)$ generated by
    $\textstyle \bigcup_{i\in I } G_i$ where $G_i$, $i\in I$, is the set composed of

    \quad (i) $\pi_0(x_1x_2\cdots x_{m_i}y_1y_2\cdots
y_{\ell_i + 1})$   for   any pair of sequences  $(x_1, \ldots x_{m_i})$ and
$(y_1, \ldots, y_{\ell_i + 1})$ as in \eqref{def:adm1}

\qquad and

    \quad (ii) $h_i -\ell_i\idu$ for $\idu$ the identity in $\U(L_0)$.
 \sm

\item \label{def:admd} We call \[ \Se^\la(L, \g)= \U(L_0)/J^\la,\]
 often abbreviated  by $\Se^\la$, the \textit{Seligman algebra\/} since it seems to have
    appeared for the first time in \cite{Sel81},  but see Theorem~\ref{thm:chari=seli}.
    By construction, it is a unital associative $\kk$-algebra. Denoting by $u \mapsto \can (u)$  the
    canonical map $\U(L_0) \to \Se^\la$, a $\la$-admissible module $M$
    becomes an $\Se^\la$-module under the action $\can( u) \cdot m = u \cdot
    m$. In this way we obtain an {\em isomorphism between the category of
    $\la$-admissible $L_0$-modules whose morphisms are $L_0$-module maps
    and the category $\Se^\la\MOD$ of left $\Se^\la$-modules.} In the future
    we will take this isomorphism as an identification.
\end{inparaenum}


The motivation for studying admissible $L_0$-modules comes from the close
connection  between the categories $\scI(L,\g)^\la$ and $\Se(L,\g)^\la \MOD$
which we investigate now. We note that any
weight space of a weight module of $L$ is invariant under the action of
$L_0$.

\subsection{Proposition (Restriction)}\label{prop:res-1}
{\it Let $\lambda\in \euP_+$ and let $V\in\Ob\cI(L,\g)^\lambda$.
	\begin{enumerate}[\rm (a)]

	\item  \label{prop-res-1-a} The weight space $V_\lambda$ is a $\lambda$-admissible $L_0$-module.\sm
	\item  \label{prop-res-1-b} If $V$ is irreducible as an $L$-module, then $V_\lambda$ is irreducible as an $L_0$-module. \sm

    \item \label{prop:res1c} If $V$ is a cyclic $L$-module with generator
        $v_\la \in V_\la$, then $V_\la$ is a cyclic $L_0$-module.
\end{enumerate}
}

\begin{proof} \eqref{prop-res-1-a} Assume $\lambda=\sum_{i\in I}\ell_i\varpi_i$. It follows from
\eqref{def:cat-I-bd1a} that for any $i\in I$ and any sequence $(y_1,\ldots,
y_{\ell_i + 1})$ with $y_j \in L_{-\al_i}$ we have $ y_1\cdot y_2\cdots
y_{\ell_i + 1}\cdot V_\lambda \in V_{\la - (\ell_i + 1) \al_i} =0$. Now let
$(x_1,\ldots,x_{m_i})$ be a sequence as in \eqref{def:adm1} and consider the
monomial $\mathbf{m}\defeq x_1x_2\cdots x_{m_i}y_1y_2\cdots y_{\ell_i
+1}\in\U(L)_0$. By what we just observed, $\mathbf m\cdot V_\lambda=0$. On
the other hand, converting $\mathbf m$ to the PBW order given by the
decomposition $L=L_-\oplus L_0\oplus L_+$ of \eqref{def:cat-I-bd-00} yields
the equation $\mathbf m = \pi_0(\mathbf m) + \mathbf m'$, $\mathbf
m'\in\U(L)L_+$. Since $L_+\cdot V_\lambda=0$ by hypothesis, this shows that
$\pi_0(\mathbf m)\cdot V_\lambda=0$, which proves \eqref{prop-res-1-a}.

\eqref{prop-res-1-b} To see that $V_\lambda$ is irreducible as an
$L_0$-module, suppose $U\subseteq V_\lambda$ is a non-zero submodule. By the
irreducibility of $V$ and the PBW Theorem, we have $V=\U(L)\cdot
U=\U(L_-)\cdot U$. It follows that $V_\lambda=U$, which proves the second
statement.

\eqref{prop:res1c} Applying again the PBW Theorem we have $V= \U(L) v_\la =
\U(L_-)_+\U(L_0) v_\la$ with $\U(L_0) v_\la \subset V_\la$ and $\U(L_-)_+
V_\la \subset \bigoplus_{\mu \preceq \la, \mu \ne \la} V_\mu$ for $\U(L_-)_+$
the augmentation ideal of $\U(L_-)$. It follows that $\U(L_0) v_\la = V_\la$
(and $\U(L_-)_+ V_\la = \bigoplus_{\mu \preceq \la, \mu \ne \la} V_\mu$).
\end{proof}

The next result concerns the opposite direction: Associating to an admissible
$L_0$-module an integrable $L$-module.

\subsection{Proposition (Integrable induction)}\label{lem:corresp-L0}{\it 
\begin{inparaenum}[\rm (a)]
Fix $\lambda=\sum_{i\in I} \ell_i\varpi_i\in\euP_+$.

\item\label{corres-a}
Let $M$ be a
$\la$-admissible $L_0$-module.  We give $M$ the structure of an $(L_0\oplus
L_+)$-module by requiring $L_+\cdot M=0$ and define the induced $L$-module
 \begin{align*}
      \widetilde{M} &\defeq\U(L)\otimes_{\U(L_0\oplus L_+)} M \qquad \text{and its submodule}\\
      \wid{N} &\defeq \tsum_{i\in I} \tsum_{m\in M} \U(L) \cdot (f_i^{\ell_i +1} \ot m).
  \end{align*}
Then the $L$-module $\bfI(M)\defeq \wid{M}/\wid{N} \in \Ob \scI(L,\g)^\la$.
It has the following properties. \sm
\begin{inparaenum}[\rm (i)]

\quad \item \label{corres-ai} The map $\xi_M \co M \to \bfI(M)_\la$,
$m\mapsto 1 \ot_{\U(L_0 \oplus L_+)} m + \wid{N}$, is an isomorphism of
$L_0$-modules. \sm

\quad \item \label{corres-aii} If $M$ is a cyclic $L_0$-module, then
$\bfI(M)$ is a cyclic $L$-module.

\end{inparaenum}
 \sm

  \item\label{corres-b}
  For an $L_0$-module map $f\co M_1 \to M_2$ of $\la$-admissible $L_0$-modules the map
  $\wid{f} = \Id_{\U(L)} \ot f \co \wid{M_1} \to \wid{M_2}$ satisfies $\wid{f}(\wid{N_1}) \subset \wid{N_2}$ and hence induces an $L$-module map $\bfI(f) \co \bfI(M_1) \to \bfI(M_2)$. It has the following properties.

  \begin{inparaenum}[\rm (i)]
    \item \label{corres-bi} The diagram below commutes.
\begin{equation}\label{corres-bi1}  \xymatrix@C=50pt{ M \ar[r]^f \ar[d]_{\xi_M}
                   & N \ar[d]^{\xi_N} \\
                \bfI(M)_\la \ar[r]_{\bfI(f)_\la} &\bfI(N)_\la } \end{equation}

   \item\label{corres-bii} If $f$ is surjective, then so is\/ $\bfI(f)$.
  \end{inparaenum}
\end{inparaenum}
}

\begin{proof}
(a) Since $\U(L)$ is a weight module by \ref{int-mod}\eqref{int-mod3} and
\ref{wei-mod}, it follows that $\wid{M}$ too is a weight module. By the
PBW Theorem, $\wid{M} = \U(L_-) \ot_{\C} M$ as vector spaces, whence all
weights of $\wid{M}$ are contained in $\la - \NN[\Th_+]$. The same then holds
for the quotient $\bfI(M)$. To see that this $L$-module is integrable,
observe that by construction all $e_i$ and $f_i$ act locally nilpotently on
the generating set $(1\ot M) + \wid{N}$ of the $L$-module $\bfI(M)$. Since
$\ad_L e_i$ and $\ad_L f_i$ are also locally nilpotent (in fact nilpotent),
\cite[Lem.~3.4(b)]{Kac90} says that $e_i$ and $f_i$ act locally nilpotently
on $\bfI(M)$, whence $\bfI(M) \in \scI(L,\g)^\la$ by \ref{int-mod}. \sm

For the proof of \eqref{corres-ai}, we first consider a fixed $i\in I$ and
$m\in M$,  and set $N=N(i,m) = \U(L) \cdot (f^{\ell_i + 1} \ot m)$.  As a
submodule of a weight module (or since $\U(L)$ is a weight module), $N$ is a
weight module with weights contained in $\la-\NN[\Th_+]$ since this is so for
$\wid{M}$. The crucial point now is to show that the weights of $N$ are
strictly below $\la$. Suppose, to the contrary, that $N_\lambda\ne 0$. That
is, there is some
	\[ u=\sum_{j=1}^p u_j^-u_j^0u_j^+\in\U(L),\qquad u_j^\pm\in\U(L_\pm),\quad u_j^0\in\U(L_0),\]
with $0 \ne uf_i^{\ell_i+1}\otimes m\in N(i,m)_\lambda$, where we may assume
without loss of generality that $0 \ne u_j^+f_i^{\ell_i+1}\otimes m\in N(i,m)$
for each $1\le j \le p$ and that the elements $u_j^\pm$ are ``monomials'' in
$\U(L)$ in the following sense
\[ u_j^+=\prod_{s=1}^{r_j} b_{\eta_{s,j}},\quad b_{\eta_{s,j}}\in L_{\eta_{s,j}},\quad
     \eta_{s,j}\in\Th_+,\qquad
\text{and}\qquad u_j^-=\prod_{t=1}^{q_j} c_{\nu_{t,j}},\quad c_{\nu_{t,j}}\in L_{\nu_{t,j}},\quad
    \nu_{t,j}\in- \Th_+.
\]
Since $uf_i^{\ell_i+1}\otimes m\in N(i,m)_\lambda$ and $m\in \wid{M}_\la$, we
have
\[ \sum_{t=1}^{q_j} \nu_{t,j}+\sum_{s=1}^{r_j} \eta_{s,j}=(\ell_i+1)\alpha_i,\qquad 1\le j\le p.\]	
If $\nu_{t,j}\ne 0$ for some pair $(t,j)$, then $\sum_{s=1}^{r_j}
\eta_{s,j}=(\ell_i+1)\alpha_i +\mu_j$ for some $\mu_j\in \NN[\Th_+]\backslash\{0\}$, which
implies that $u_j^+f_i^{\ell_i+1}\otimes m=0$. Thus, we have in fact
\[ u=\sum_{j=1}^p u_j^0u_j^+ \quad \text{and} \quad
\sum_{s=1}^{r_j} \eta_{s,j}=(\ell_i+1)\alpha_i,\qquad 1\le j\le p.\] It then
follows that each $\eta_{s,j}$ is a positive multiple of $\al_i$. Thus
\begin{align*}
uf_i^{\ell_i+1}\otimes m&=\sum_{j=1}^p u_j^0u_j^+f_i^{\ell_i+1}\otimes m=
\sum_{j=1}^p 1\otimes u_j^0\pi_0(u_j^+f_i^{\ell_i+1})\cdot m=
\sum_{j=1}^p 1\otimes u_j^0\pi_0(b_{\eta_{1,j}}b_{\eta_{2,j}}\cdots
b_{\eta_{r,j}}f_i^{\ell_i+1})\cdot m.
\end{align*} By the $\lambda$-admissibility of
$M$, each term on the right is zero, letting us conclude that $N_\lambda=0$.

Since $\wid{N} = \sum_{i\in I} \sum_{m\in M} N(i,m)$ is a sum of weight
modules, it now follows that the weights of $\wid{N}$ are all strictly below
$\la$. Hence the composition of $L_0$-module maps
\[ M \ni m \mapsto \wid{m} = 1 \ot_{\U(L_0\oplus L_+)} m \mapsto \wid{m} + \wid N \in (\wid{M}/\wid{N})_\la
\]
is bijective. This proves \eqref{corres-ai}. For \eqref{corres-aii} suppose
$M = \U(L_0) m_0$. Then $\wid{M}$ is cyclic as $L$-module: $M = \U(L)
\ot_{\U(L_0 \oplus  L_+)} \U(L_0) m_0 = \U(L) \ot \kk m_0$. This implies
\eqref{corres-aii}. (b) is straightforward from the definitions, observing
that $f$ is an $L_0\oplus L_+$-module homomorphism.  \end{proof}

\subsection{The left-regular representation.}
Let $\la\in\euP_+$. Then the map $L_0\to \mathrm{End}_\kk(\U(L_0))$ given by
$u\mapsto \mathrm{L}_u$, where $\mathrm L_u(x)=ux$ for every $x\in \U(L_0)$,
is a representation of the Lie algebra $L_0$ on $\U(L_0)$. The ideal $J^\la$
is $L_0$-invariant, and it is clear from the definition that the quotient
$L_0$-module $\Se^\la=\U(L_0)/J^\la$ is $\la$-admissible, or equivalently, an
$\Se^\la$-module. We will call the latter $\Se^\la$-module the
\emph{left-regular $\Se^\la$-module}, and denote it by $_\reg \Se^\la$.

The next lemma shows that the ideal $J^\la$
defined in \ref{def:adm}\eqref{def:adm3}
 contains elements that are obtained by a relaxation of the definition of the generators.

\subsection{Lemma.}
\label{gener-Jla-Hadi}
Let $\la=\sum_{i\in I}\ell_i\varpi_i\in\euP_+$, and fix $i\in I$. Let $(x_1,\ldots,x_{m_i})$ and
$(y_1,\ldots,y_{n_i})$
be a pair of sequences
satisfying
\[
x_j\in L_{\beta_j},\, \beta_j\in \Theta_+,\,1\leq j\leq m_i,
\ \ \
y_j\in L_{-\gamma_j},\,\gamma_j\in\Theta_+,\,
1\leq j\leq n_i, \
\text{ and }\ \sum_{j=1}^{m_i}\beta_j=
\sum_{j=1}^{n_i}\gamma_j=r_i\alpha_i,
\]
where $r_i>\ell_i$.
Then $\pi_0(x_1\cdots x_{m_i}y_1\cdots y_{n_i})\in J^\la$.

\begin{proof}
By the universal property of $\U(L_0)$ , the $L_0$-action on $_\reg\Se^\la$
extends to a representation of the associative algebra $\U(L_0)$
on$_\reg\Se^\la$. The annihilator of this representation is $J^\la$.
Therefore it is enough to prove that $\pi_0(x_1\cdots x_{m_i}y_1\cdots
y_{n_i})({_\reg\Se^\la})=0$. By Proposition
\ref{lem:corresp-L0}\eqref{corres-ai}, it is enough to verify that $x_1\cdots
x_{m_i}y_1\cdots y_{n_i}$ acts on $\bfI(M)$ trivially. But indeed already
$y_1\cdots y_{n_i}$ acts  on $\bfI(M)$ trivially. The latter statement
follows from an $\lsl_2$-theory argument similar to the one given for
\eqref{def:cat-I-bd1a}.
\end{proof}

\subsection{Functors $\bfR$ and $\bfI$.} \label{func}
Let $\la \in \euP_+$. We abbreviate $\scI(L,\g)^\la =\scI^\la $ and
$\Se(L,\g)^\la = \Se^\la$.

For any morphism $\vphi \co V \to V'$ in $\scI^\la$ the restriction
$\vphi_\la \co V_\la \to V'_\la$ is an $L_0$-module homomorphism, whence, by
Proposition~\ref{prop:res-1}, a homomorphism of $\Se^\la$-modules. The
assignments $V \mapsto V_\la$ and $\vphi \mapsto \vphi_\la$ define a
restriction functor
\[\bfR^\la_{(L,\g)} \co \scI^\la \to \Se^\la \MOD, \]
which is exact. Using Proposition~\ref{lem:corresp-L0}, we have a functor
\[ \bfI^\la_{(L,\g)} \co \Se^\la \MOD \to \scI^\la\]
in the other direction, defined by $M \mapsto \bfI(M)$ and $f\mapsto
\bfI(f)$. Unless there is a danger of confusion we will abbreviate
\[ \bfR = \bfR^\la_{(L,\g)} \quad \hbox{and} \quad \bfI = \bfI^\la_{(L,\g)}.\]
A  re-interpretation of Proposition~\ref{lem:corresp-L0}\eqref{corres-bi} is
that the family of isomorphisms $(\xi_M)$ is an isomorphism of functors
\begin{equation}\label{func1}
 \Id_{\Se^\la\MOD} \xrightarrow{\cong} \bfR \circ \bfI .
 \end{equation}
We will deal with the composition of $\bfR$ and $\bfI$ in the other order in
the following result. It is in spirit closely related to \cite[Prop.~5 and
Cor.~2]{CFK10} dealing with the special case of map algebras, although the
functor $\bW$ used in loc.\ cit.\ is not the same as our functor $\bfI$.

\subsection{Proposition (The functor $\bfI \circ \bfR$)} \label{poth}
{\it 
We fix $\la \in \euP_+$ and use the abbreviations of {\rm \ref{func}}. \sm

\begin{enumerate}[\rm (a)]
  \item \label{potha} For $V\in \scI^\la$ the structure map $\U(L) \ot_\kk
      V \to V$ induces an $L$-module map $\eta_V \co \bfI(V_\la) \to V$
      with image $\U(L) \cdot V_\la$.
\sm

  \item\label{pothb}  The collection of maps $(\eta_V)$ is a natural
      transformation $\bfI \circ \bfR \Rightarrow \Id_{\scI^\la}$. \sm

   \item\label{pothc} $\bfR$ is a right adjoint of\/ $\bfI$.
\sm

  \item\label{pothd}  The functor $\bfI$ maps projective objects to projective objects.
\end{enumerate}
}

\begin{proof} \eqref{potha} The structure map $\U(L) \times V_\la \to V$, $(u,v_\la) \mapsto u \cdot v_\la$ is $\U(L_0 \oplus L_+)$-balanced, whence gives rise to an $\U(L)$-linear map $\U(L) \ot_{\U(L_0 \oplus L_+)} V_\la = \wid{V_\la} \to V$ with $\U(L)$ acting on the left factor of $\wid{V_\la}$. Because of \eqref{def:cat-I-bd1a}, this map annihilates the submodule $\wid{N}$ of Proposition~\ref{lem:corresp-L0}\eqref{corres-a}. The quotient map is the map $\eta_V$ of our claim.

Regarding \eqref{pothb} we need to show that for any morphism $f\co V\to V'$
in $\scI^\la$ the diagram
\[\xymatrix@C=50pt{
     \bfI(V_\la) \ar[d]_{\bfI(f_\la)} \ar[r]^{\eta_V} & V \ar[d]^f \\
         \bfI(V'_\la) \ar[r]_{\eta_{V'}} & V'
}\]
commutes, which is immediate from the definitions.

\eqref{pothc} We need to prove that for all $M\in \Se^\la\MOD$ and $V\in
\scI^\la$ there exists an isomorphism of abelian groups
\[
  \ta_{M,V} \co \Hom_{\scI^\la}\big( \bfI(M), V \big) \to \Hom_{\Se^\la}\big(M, \bfR(V)\big)
\]
which is natural in $M$ and $V$. There is no other reasonable choice but to
define $\ta_{M,V}(\vphi) = \vphi_\la \circ \xi_M= \bfR(\vphi)\circ \xi_M$ for
$\vphi \in \Hom_{\scI^\la}\big( \bfI(M), V \big)$. This map is clearly
additive in $\vphi$. It is injective since $\bfI(M)$ is generated by
$\xi_M(M) \cong M$ as $L$-module. It is surjective since
$\tau_{M,V}\big(\eta_V \circ \bfI(f)\big) = f$ for any $f\in
\Hom_{\Se^\la}\big(M, \bfR(V)\big)$. Indeed, the commutative diagram
\eqref{corres-bi1} specializes to
\[
\xymatrix@C=60pt{ M \ar[r]^f \ar[d]_{\xi_M}
                   & V_\la \ar[d]^{\xi_{V_\la}} \\
                \bfI(M)_\la \ar[r]_{(\bfR \circ \bfI)(f)} &\big( \bfI(V_\la)\big)_\la }
\]
so that $\ta\big( \eta_V \circ \bfI(f)\big) = \bfR\big( \eta_V \circ
\bfI(f)\big)\circ \xi_M = \bfR(\eta_V) \circ (\bfR \circ\bfI)(f) \circ \xi_M
= \bfR(\eta_V) \circ \xi_{V_\la} \circ f$. Our claim then follows from
$\bfR(\eta_V) \circ \xi_{V_\la} = \Id_{V_\la}$ which easily follows from the
definitions. We leave it to the reader to check that $\ta_{M,V}$ is natural
in $M$ and $V$.

\eqref{pothd} is a standard property of adjoint functors. \end{proof}

\subsection{Corollary (The Weyl module $W(\la)$)
} \label{cor:we}
{\it 
The map $\eta_{W(\la)} \co \bfI(W(\la)_\la) \to W(\la)$ of Proposition~{\rm
\ref{poth}\eqref{potha}} is an isomorphism of $L$-modules. Moreover, let
$_\reg \Se^\la$ be the left-regular $\Se^\la$-module and  $f\co {_\reg
\Se^\la} \to W(\la)_\la$ be the $\Se^\la$-module map sending $1_{\Se^\la}$ to
the generator $w_\la$ of the cyclic $L_0$-module $W(\la)_\la$. Then the map
$\bfI(f) \co \bfI({_\reg \Se^\la}) \to \bfI(W(\la)_\la)$ is an isomorphsim of
$L$-modules. In sum,
\[
      \bfI(W(\la)_\la) \cong W(\la) \cong \bfI({_\reg \Se^\la})
\]
as $L$-modules.}

\begin{proof} By Proposition~\ref{prop:res-1}\eqref{prop:res1c},
$W(\la)_\la$ is a cyclic $L_0$-module, whence $\bfI(W(\la)_\la)$ is a cyclic
 $L$-module in $\scI^\la$ by Proposition~\ref{lem:corresp-L0}\eqref{corres-aii}.
The map $\eta_{W(\la)}$ sends the generator $\bfI(w_\la)$ onto the canonical
generator $w_\la$ of $W(\la)$. It then follows from the universal property of
$W(\la)$, cf.\ Proposition~\ref{prop:global-weyl-def-rel}\eqref{prop:gwdrb},
that $\eta_{W(\la)}$ is an isomorphism.

Since $W(\la)_\la$ is a cyclic $\Se^\la$-module we have a surjective
$\Se^\la$-module map $f \co {_\reg \Se^\la} \to W(\la)_\la$ and then by
Proposition~\ref{lem:corresp-L0}\eqref{corres-bii} a surjective $L$-module
map $\bfI(f) \co \bfI({_\reg \Se^\la}) \to \bfI(W(\la)_\la)$. Furthermore, by
Proposition \ref{lem:corresp-L0}\eqref{corres-aii} the module
 $\bfI({_\reg \Se^\la})$ is cyclic, and therefore
 the canonically induced map $W(\la)\to
\bfI({_\reg \Se^\la})$ is a surjection. It then follows as in the first part
that $\bfI(f)$ is an isomorphism.
\end{proof}

\subsection{Theorem} \label{thm:chari=seli}{\em
For any $\lambda\in\EuScript P_+$ we have $\mathbb S^\lambda=\mathbf
A_\lambda$ as associative $\kk$-algebras.} \sm

 Because of this result, we will refer to $\mathbf A_\la$ and $\Se^\la$ as the
{\em Seligman-Chari-Pressley algebra\/}  in the future. For simpler notation
we will  denote it by $\Se^\la$.

\begin{proof}
Recall that $\Se^\la=\U(L_0)/J^\la$ and $\mathbf
A_\lambda=\U(L_0)/\Ann_{\U(L_0)}(w_\lambda)$. From Proposition
\ref{prop:res-1}(a) it follows that $J^\lambda\subseteq
\Ann_{\U(L_0)}(w_\lambda)$. To complete the proof, we need to show the
reverse inclusion. Consider the $L_0$-module map
\[
f:\Se^\lambda\to W(\lambda)_\lambda\ ,\ u\mapsto u\cdot w_\lambda.
\]
By Proposition \ref{lem:corresp-L0}(b) and Corollary \ref{cor:we}, $f$ gives
rise to an isomorphism of $L$-modules
\[
\mathbf{Int}(f):\mathbf{Int}({_\reg \Se^\la})\to
\mathbf{Int}(W(\lambda)_\lambda)\cong W(\lambda).
\]

Write $\lambda=\sum_{i=1}^r\ell_i\varpi_i$ and fix $u_\circ\in\Ann_{\mathrm
U(L_0)}(w_\lambda)$. Our goal is to prove that $u_\circ\in J^\lambda$. Let
$\overline u_\circ$ be the projection of $u_\circ$ in $\Se^\lambda$. First
note that $\mathbf{Int}(f)(1\otimes \overline u_\circ)=1\otimes (u_\circ\cdot
w_\lambda)=0$, and since $\mathbf{Int}(f)$ is an isomorphism, we should have
$1\otimes \overline u_\circ=0$ as an element of $\mathbf{Int}({_\reg
\Se^\la})$. It follows that
\[
1\otimes \overline u_\circ\in \sum_{i=1}^r\sum_{s\in \Se^\lambda}\U(L)\cdot (f_i^{\ell_i+1}\otimes s).
\]
Fix a basis of $L$ consisting of $\mathfrak h$-root vectors. The basis of $L$
yields a PBW basis of $\mathrm  U(L)$.  We can now write
\begin{equation}
\label{eq-hadi-1ucrc}
1\otimes\overline u_\circ=\sum_{i=1}^r\sum_{j=1}^\infty
q_{i,j}f_i^{\ell_i+1}\otimes s_{i,j},
\end{equation}
where every $q_{i,j}$ is a monomial in the PBW basis of $\U(L)$, and all but
finitely many of the $s_{i,j}\in \Se^\lambda$ are zero. Next fix $ i,j$ and
write  \[ q_{i,j}=y_1\cdots y_a h_1\cdots h_b x_1\cdots x_c,
\] where $y_s\in L_{-\beta_s}$ for $1\leq s\leq a$,
$h_s\in L_0$ for $1\leq s\leq b$, and $x_s\in L_{\gamma_s}$ for $1\leq s\leq
c$, with $\beta_1,\ldots, \be_a,\gamma_1,\ldots,\gamma_c\in \Th_+$. The
equality of $\mathfrak h$-weights of both sides of \eqref{eq-hadi-1ucrc}
implies that
\[
\gamma_1+\cdots +\gamma_c=\beta_1+\cdots+\beta_a+(\ell_i+1)\alpha_i.
\]
If $a\geq 1$, then the $\mathfrak h$-weight of $x_1\cdots x_c f_i^{\ell_i+1}$
is in $\NN[\Th_+]\setminus\{0\}$, from which it follows that $x_1\cdots x_c
f_i^{\ell_i+1}\in \U(L)L_+$, and therefore $q_{i,j}f_i^{\ell_i+1}\otimes
s_{i,j}=0$. We can therefore assume $a=0$. Then $x_1\cdots
x_cf_i^{\ell_i+1}\in \U(L)_0$ and
\[
   x_1\cdots x_cf_i^{\ell_i+1}=\pi_0(x_1\cdots x_cf_i^{\ell_i+1})+u_+ \ \text{ for some }\
u_+\in \U(L)L_+,
\] so that
\begin{align*}
q_{i,j}f_i^{\ell_i+1}&\otimes s_{i,j}=
\left(h_1\cdots h_b\pi_0(x_1\cdots x_cf_i^{\ell_i+1})+
h_1\cdots h_bu_+\right)
\otimes s_{i,j}\\
&=
h_1\cdots h_b\pi_0(x_1\cdots x_cf_i^{\ell_i+1})
\otimes s_{i,j}=1\otimes
h_1\cdots h_b\pi_0(x_1\cdots x_cf_i^{\ell_i+1})
s_{i,j}.
\end{align*}
Note that $ h_1\cdots h_b\pi_0(x_1\cdots x_cf_i^{\ell_i+1}) s_{i,j}=0$ in
$\Se^\lambda $ since already $\pi_0(x_1\cdots x_cf_i^{\ell_i+1})=0$ in
$\Se^\la$. Therefore the above arguments show that $1\otimes \overline
u_\circ=0$  as an element of $\U(L)\otimes_{\U(L_0\oplus L_+)}\Se^\lambda$.
Since $\U(L)$ is a free $\U(L_0\oplus L_+)$-module, it follows that
$\overline u_\circ=0$ in $\Se^\lambda$, so that $u_\circ\in J^\lambda$.
\end{proof}

\ms

In the remainder of this section we will derive results concerning the
structure of $\Se^\la$.

\subsection{Example $\U(L_0)_+ \subseteq J^\la$:} \label{ex:seli}
Let $\U(L_0)_+$ be the augmentation ideal of $\U(L_0)$. Since $\U(L_0)_+$ is
a maximal ideal there are exactly the cases \eqref{ex:seli1} and
\eqref{ex:seli2} below whenever $\U(L_0)_+\subseteq J^\la$: \sm

\begin{inparaenum}\item\label{ex:seli1} $\U(L_0)_+ = J^\la:$  We have
\[
    \U(L_0)_+ = J^\la \quad \iff \quad \la = 0,
    \]
and in this case (obviously) $\Se^\la(L,\g) = \kk 1$ is $1$-dimensional.

Indeed, if $\U(L_0)_+ = J^\la$, then $h - \la(h)\idu \in J^\la$ and $h\in
J^\la$ forces $\la=0$. Conversely, assume $\la = 0$. Then $\pi_0(L_{\al_i}
L_{-\al_i})=[L_{\al_i}, L_{-\al_i}] \subset J^0$ by definition  of the ideal
$J^\la$ in \ref{def:adm}\eqref{def:adm3}. Using the notation of
Lemma~\ref{rogrfa}\eqref{rogrfa-b}, we also have $\pi_0(L_{\be_c/2}
L_{\be_c/2}L_{-\be_c}) = [L_{\be_c/2}, [L_{\be_c/2}, L_{\be_c}] =
[L_{\be_c/2}, L_{\be_c/2}] \subset J^0$. It then follows from
\eqref{eq-hadi01} and that $J^0$ is generated by $L_0$, i.e., $J^0 =
\U(L_0)_+$.

We point out that $\dim \Se^\la= 1$ does not imply $\la = 0$. For example, as
already shown in \cite[p.~59]{Sel81}, for $L=\lsl_n(A)$, $A$ central-simple,
$\la = \ell \varpi_i$ for $1< i <n-1$, we always have $\dim \Se^\la \le 1$
and $\dim = 1 \iff d|\ell$ where $d$ is the degree of $A$. We will see in
\ref{ex:disc-th} and \ref{ex:csan} that examples with $\dim \Se^\la \le 1$
arise naturally.\sm

 \item\label{ex:seli2} $\U(L_0) = J^\la$, i.e., $\Se^\la(L,\g) = \{0\} :$ Then $\la \ne 0 $,
  and the only $\la$-admissible $L_0$-module is $\{0\}$. Hence any module $V$ in
   $\scI^\la(L,\g)$ has $V_\la = \{0\}$.
\end{inparaenum}

\subsection{The subalgebras $\Se^\la_i < \Se^\la$}\label{subi}
One approach to the analysis of the structure of
$\Se^\la(L, \g)
= \Se^\la$ is through certain subalgebras $\Se^\la_i$, which we now define using
the notation of \ref{g-not} and \ref{def:adm}.
For $i\in I$ let
\begin{equation}\label{subi1}
 \begin{split}  L_i &= \Big(\textstyle  \sum_{\beta\in\R\alpha_i\cap \Theta} \, L_{\beta} \Big)
   \bigoplus \sum_{\beta\in\R_+\alpha_i\cap\Theta}\, [ L_{\beta}, L_{-\beta} ], \\
     \g_i &= \g_{\al_i} \oplus \kk h_i \oplus \g_{-\al_i} \simeq \lsl_2(\kk). \end{split}
\end{equation}
Then $L_i$ is a subalgebra of $L$ and $(L_i,\g_i)$ is either an $\rma_1$-graded or a $\mathrm{BC}_1$-graded
Lie algebra with $0$-part $L_{i0}  = L_i \cap L_0$.
It is immediate from the Jacobi identity  that
$L_{i0}$ is an ideal of $L_0$. Let $\xi_i $ be the restriction of the
fundamental weight $\varpi_i$ to $\kk h_i$, the fundamental weight of $\g_i$.
Under the canonical injection $\io_i \co \U(L_{i0}) \hookrightarrow \U(L_0)$
the ideal $J^{\ell_i \xi_i}$ is mapped to $J^\la$, giving rise to a
commutative diagram
\[\xymatrix{
    0 \ar[r] & J^{\ell_i \xi_i} \ar[r] \ar@{^{(}->}[d] &\U(L_{i0})
         \ar[r]^>>>>>{\can}
             \ar@{^{(}->}[d]^{\io_i}
           &\Se^{\ell_i\xi_i}(L_i, \g_i) \ar[r] \ar@{-->}[d]^{\phi_i} & 0 \\
 0 \ar[r] & J^\la  \ar[r] & \U(L_0) \ar[r]^{\can}
           &\Se^\la (L, \g) \ar[r] & 0
}\] where $\Se^{\ell_i\xi_i}(L_i,\g_i)$ is
defined for $(L_i,\g_i)$ as in Definition \ref{def:adm}, the maps ``$\can$" are the canonical quotient
maps, and $\phi_i$ is the unital algebra homomorphism induced by $\io_i$. We
denote by
\[   \Se^\la_i = \phi_i \left( \Se^{\ell_i \xi_i}(L_i, \g_i) \right)
      < \Se^\la
\]
the image of $\phi_i$. Since $L_0 = \sum_{i\in I} L_{i0}$ by
\eqref{eq-hadi01}, the algebra $\Se^\la$ is generated as associative algebra
by the subalgebras $\Se^\la_i$, $i\in I$. Note $\Se^\la_i = \kk \idu $ if
$\ell_i = 0$ since then $J^{0\xi_i}$ is generated by $L_{i0}$ and thus equals
the augmentation ideal of $\U(L_{i0})$.

Using the subalgebras $\Se^\la_i$ to determine the structure of $\Se^\la$
leads to two problems: First, $\Se^\la_i$ need not be isomorphic to
$\Se^{\ell_i\xi_i}(L_i, \g_i)$ and, second, the interplay between the
subalgebras $\Se_i^\la$ of $\Se^\la$ seems to be complicated in general. But
it can be understood in the following situation.

\subsection{Proposition} \label{seli-str}{\em
Let $L$ be a $\De$-graded Lie algebra and let $M_0 \ideal L_0$ be an ideal of
$L_0$ satisfying
\begin{enumerate}[\rm (i)]

\item $M_0 \subset J^\la$, viewed in $\U(L_0)$,

\item $M_0 = \sum_{i\in I} M_0 \cap L_{i0}$, and

\item denoting by $\overline{\phantom{0}}$ the canonical quotient map $L_0
    \to     L_0/M_0$    we     have $L_0 /M_0 = \boxplus_{i\in I} \overline{L_{i0}}$,
    a direct sum of ideals    $\overline{L_{i0}}$ of $L_0/M_0$.
\end{enumerate}

\noindent Then multiplication induces a unital algebra isomorphism
between the tensor product algebra of the $\Se^\la_i$, $i\in I$, and $\Se^\la :$
\[
    \textstyle \bigotimes_{i\in I} \Se^\la_i \xrightarrow{\simeq} \Se^\la.
\]}

\begin{proof} By \cite[\S2.7, Cor.~6]{Bou71}, multiplying the factors is an isomorphism of vector
spaces $\textstyle \bigotimes_{i\in I} \U(\overline{L_{i0}}) \to
\U(\overline{L_0})$. But since the factors $ \U(\overline{L_{i0}})$ commute
in $\U(\overline{L_0})$, this map is in fact a unital algebra isomorphism. To
simplify the notation, we will view it as an identification in the following.

Let $\id(M_0) \ideal \U(L_0)$ be the ideal of $\U(L_0)$ generated by $M_0$.
By \cite[\S2.3, Prop.~3]{Bou71} this is the kernel of the quotient map $\bar
p \co \U(L_0) \to \bigotimes_{i\in I} \U(\overline{L_{i0}})$. Observe that
$\bar p$ maps $\U(L_{j0}) \subset \U(L_0)$ onto $\U(\overline{L_{j0}})
\subset \bigotimes_{i\in I} \U(\overline{L_{i0}})$. We have the following
commutative diagram, defining the unital surjective algebra homomorphism
$\vphi :$
\[\xymatrix{
    0 \ar[r] & \id(M_0) \ar[r] \ar@{^{(}->}[d] &\U(L_0)
         \ar[r]^>>>>>{\bar p}        \ar@{=}[d]
           &   \textstyle \bigotimes_{i\in I} \U(\overline{L_{i0}})
            \ar[r] \ar@{-->}[d]^{\vphi} & 0 \\
 0 \ar[r] & J^\la  \ar[r] & \U(L_0) \ar[r]^{\can}
           &\Se^\la (L, \g) \ar[r] & 0
}\] The kernel of $\vphi$ is $\bar p(J^\la)$, the ideal of $\bigotimes_{i\in
I} \U(\overline{L_{i0}})$ generated by $\bigcup_{i\in I} \bar p(G_i)$. From
the algebra structure of $\bigotimes_{i\in I} \U(\overline{L_{i0}})$ it
follows that $\bar p(J^\la) = \sum_{i\in I} \U(\overline{L_{10}}) \ot \cdots
\U(\overline{L_{i-1,0}}) \ot   X_i \ot \U(\overline{L_{i+1, 0}}) \ot \cdots
\U(\overline{L_{r0}})$ where $X_i$ is the ideal of $\U(\overline{L_{i0}})$
generated by $\bar p(G_i)$. Hence,

\[ \Se^\la \simeq \big( \textstyle \bigotimes_{i\in I} \U(\overline{L_{i0}}) \big)
       \big/ \Ker(\vphi) \simeq
    \bigotimes_{i\in I} \U(\overline{L_{i0}})/X_i
 \]
But $\Se^\la_i = \U(L_{i0})/ (J^\la \cap \U(L_{i0})) = \bar p(\U(L_{i0}))\big
/ \bar p( \U(L_{i0}) \cap J^\la) = \U(\overline{L_{i0}})/X_i$.\end{proof}

\subsection{Examples}\label{seli-se}
\begin{inparaenum}
\item\label{seli-se-a} Let $(L,\g)= (\g\ot A, \g)$ be a map algebra as in \ref{rg-exa}\eqref{reg-exa1}.
Then $L_0 = \h \ot A$ is abelian and a direct sum of the ideals $\g_i \ot A$. Thus,
Proposition~\ref{seli-str} applies with $M_0 = \{0\}$. Observe that in this case $\bar p$
is an isomorphism and $J^{\ell_i \xi_i} = X^i$ in the notation of the proof above. It follows that
the maps $\phi_i$ of \ref{subi} are isomorphisms. We will use this in the proof of Theorem~\ref{thm:comm}.
\sm

\item\label{seli-se-b} Suppose $\supp \la = \{ i\in I : \ell_i >0 \}$ is
{\em totally disconnected\/} in the sense that $i\in \supp \la$ implies that both $i-1\not\in
\supp\la$ and $i+ 1\not\in \supp \la$. Then $M_0 = \sum_{k\not\in \supp \la} L_{k0}$ is
an ideal of $L_0$ satisfying the conditions (i)--(iii) of Proposition~\ref{seli-str}.
\end{inparaenum}

\sm

We also note the following structural result which generalizes \cite[Prop. I.10]{Sel81}
proven there for finite-dimensional semisimple Lie algebras as in Example~\ref{rg-exa}\eqref{reg-exa3}.

\subsection{Proposition} \label{fd-Se} {\em
Let $(L, \g, \h)$ be a $(\Theta,\De)$-graded Lie algebra, and let
$\lambda=\sum_{i\in I}\ell_i \varpi_i\in\EuScript P_+$. We set
$\ell_{\max}:=\max\{\ell_i\ :\ i\in I\}$. If\/ $\dim L_0<\infty$, then \[
\dim\Se^\lambda(L,\g)\leq (2\ell_{\max}+1)^{\dim L_0-\dim\h}.
\] If there are no connected components $\Theta_c$ of $\Theta$ satisfying $(\Theta_c,\De_c)
\cong (\mathrm{BC}_n,\mathrm{C}_n)$, then the stronger inequality
$\dim\Se^\lambda(L,\g)\leq (\ell_{\max}+1)^{\dim L_0-\dim\h}$ holds. }

\begin{proof}
We can essentially follow Seligman's proof in loc.\ cit. But we give a sketch of
the argument for the convenience of the reader.

(I) For $\al \in \Theta_+$, $x\in L_\al, y\in L_{-\al}$ and $j\in \NN_+$
there exist $a_{jk}\in \NN_+$ such that
\[ x^j y^j \equiv \textstyle \sum_{k=0}^j a_{jk}
y^k [x,y]^{j-k} x^k \mod \sum_{k=0}^{j-1} L_{-\al}^k \U(L_0)_{(j-k-1)}
L_\al^k.\]
This is proven in \cite[Lemma~I.2]{Sel81} and \cite[4.5, Lem.~5]{CFK10} for examples
\ref{rg-exa}\eqref{reg-exa3} and \ref{rg-exa}\eqref{reg-exa1} respectively. Both proofs work in our
setting. In particular, for $\al = \al_i$ and $j=\ell_i + 1$ we get, after
application of $\pi_0$, that
\begin{equation} \label{fd-Se1}
    \pi_0 (x^{\ell_i + 1} y^{\ell_i + 1} ) \equiv a_{\ell_i + 1, 0}\, [x,y]^{\ell_i
           + 1} \mod \U(L_0)_{(\ell_i)}. \end{equation}
The left hand side lies in $J^\la$, hence $[x,y]^{\ell_i+1}\in J^\la$. If
$\alpha=\frac{1}{2}\alpha_i$ (this occurs when $\Theta$ has connected
components of the form $(\Theta_c,\De_c)\cong(\mathrm{BC}_n,\mathrm{C}_n)$),
then setting $j=2\ell_i+1$ and using Lemma \ref{gener-Jla-Hadi} implies that
$[x,y]^{2\ell_i+1}\in J^\la$. The algebra $L_0$ has a basis $\sfB$ that
consists of elements of the form $[x,y]$ with $x\in L_{\alpha_i}$ and $y\in
L_{-\alpha_i}$, or $x\in L_{\frac12\alpha_i}$ and $y\in L_{-\frac12\alpha_i}$
(where in the latter case $\alpha_i$ is a long simple root of a connected
component $\De_c$ of type $\mathrm{C}_n$). Furthermore, using the Chevalley
basis of $\g$ we can assume that $\sfB$ contains a basis of $\mathfrak
h\subset\g$. By the PBW Theorem, the relation (ii) of
\ref{def:adm}\eqref{def:adm3}, and the observations made above,
$\U(L_0)/J^\la$ is spanned by monomials in the elements of $\sfB\backslash
\mathfrak h$ in which none of the exponents exceed $2\ell_i$. It follows that
$\dim\Se^\lambda(L,\g)\leq (2\ell_{\max}+1)^{\dim L_0-\dim\h}$. From the
proof it is clear that when there are no connected components of
$(\mathrm{BC}_n,\mathrm{C}_n)$ type, the stronger inequality holds.
 \end{proof}
\vfe

\section{The Seligman--Chari--Pressley algebra of map algebras}
\label{sec:rep-gencur}

We continue with the notation of the previous sections: $\kk$ is a field of
characteristic $0$ and $\g$ is a finite-dimensional split semisimple Lie
algebra over $\kk$. The ultimate goal of this section is to determine  $\Se(L, \g)^\bullet$ for the root-graded Lie $\kk$-algebra $L=\g \ot
A$, $A\in \kkalg$, of \ref{rg-exa}\eqref{reg-exa1}.

We will start this section by defining and describing the structure for
$\lsl_n(A)$ for $A$ unital associative (\ref{subsec:rev-liealg}--\ref{Lnull}).
We will then specialize to $n=2$ with
$A$ as before, and at the end consider $L=\g\ot A$ for $A \in \kkalg$.
Proceeding in this way will allow us to use some formulas needed for $\g \ot
A$ also in the case $\lsl_n(A)$ of \S\ref{sec:slna}. \sm

Thus, unless stated otherwise, in this section $n\in \NN$, $n\ge 2$, and $A$
is a unital associative, but not necessarily commutative $\kk$-algebra.

\subsection{The Lie algebra $\lsl_n(A)$.} \label{subsec:rev-liealg}
By definition $\lsl_n(A)$ is the derived algebra of the Lie algebra
$\Mat_n(A)^-= \gl_n(A)$ of $n\times n$ matrices over $A$:
\[ L = \lsl_n(A) = [\gl_n(A), \, \gl_n(A)] \]
(see \eqref{sln-des} for another description of $\lsl_n(A)$).
For general $A$ the structure map $\kk \to \kk\cdot1_A$ gives rise to an
injective homomorphism of $\lsl_n(\kk)= [\gl_n(\kk), \gl_n(\kk)]$ into
$\lsl_n(A)$. Its image is the subalgebra of $\lsl_n(A)$ consisting of the
trace-less matrices over $\kk$. In the future we will take this homomorphism
as an identification:
\[ \g = \lsl_n(\kk) \subset L.\]
We claim that $L$ is a $\De$-graded Lie algebra for the root system
$\De=\rma_{n-1}$ in the sense of \ref{def:rg}. Indeed, $\g$ is split simple
and contains
\[ \h = \{h\in \lsl_n(\kk): h \text{ diagonal}\} \]
as a splitting Cartan subalgebra. Denoting by $E_{ij}$ the usual matrix units
in $\Mat_n(A)$ and abbreviating $E_{ij}(a) = a E_{ij}$ for $a\in A$, the Lie
algebra $\lsl_n(A)$ has a weight space decomposition
\begin{align*} 
 L=  L_0 \oplus \ts \bigoplus_{i\ne j} L_{\eps_i - \eps_j},
\quad    L_{\eps_i - \eps_j} = \{ E_{ij}(a): a\in A\} =\co E_{ij}(A)
 \end{align*}
where $L_0$ consists of the diagonal matrices in $\lsl_n(A)$ and
$L_{\eps_i-\eps_j}$ is the weight space for the root $\eps_i-\eps_j \in \De$
realized in the standard way: $\eps_i$ is the projection onto the $i^{\rm
th}$-component of $h\in \h$. Hence \ref{def:rg}\eqref{def:reg1} holds. In our
setting condition \ref{def:rg}\eqref{def:reg2} says $L_0 \subseteq \sum_{i\ne
j} [E_{ij}(A), E_{ji}(A)]$. Since $\gl_n(A) = \gl_n(A)_0 \oplus
\big(\bigoplus_{i\ne j} E_{ij}(A)\big)$ is a grading by the root lattice of
$\De$, it follows from the definition of $\lsl_n(A)$ that we need to show
$[\gl_n(A)_0, \gl_n(A)_0] \subset\sum_{i\ne j} [E_{ij}(A), E_{ji}(A)]$, or
equivalently
\begin{equation} \label{subsec:rev-liealg2}
[a,b]E_{ii}\in  \sum_{i\ne j} [E_{ij}(A), E_{ji}(A)].
\end{equation}
To do so we will use the matrices
\begin{align}  e_i(a)& \defeq E_{i,i+1}(a),\quad f_i(a)\defeq E_{i+1,i}(a), \nonumber
\\  H_i(a,b) & \defeq[e_i(a),f_i(b)] = ab E_{ii} - ba E_{i+1, i+1} \in L_0,\label{rev-liealg1}
 \\ h_i(a) & = H_i(a,1_A)=H_i(1_A,a) =a(E_{ii} - E_{i+1, i+1}) \nonumber
 \end{align}
defined for $a, b\in A$ and $1\le i < n$. It is straightforward to verify
\begin{equation}\label{Lnull00}
  [a,b] E_{ii} = \begin{cases} H_i(a,b) - h_i(ba), &(1\le i < n), \\
         H_{i-1}(a,b)-h_{i-1}(ab) &(1< i \le n), \end{cases}
 \end{equation}
which clearly implies
\eqref{subsec:rev-liealg2}. \sm

We will use the standard root basis of the root system $\De$, for which the
simple roots are the $\al_i = \eps_i - \eps_{i+1}$, cf.~\ref{g-not}. Hence
the subalgebra $L_+$ (resp. $L_-$) is the subalgebra of strictly upper
(resp.\ lower) triangular matrices in $\lsl_n(A)$.

\subsection{Structure of $L_0=\lsl_n(A)_0$} \label{Lnull}
Since $L_0$ is ${\rm A}_{n-1}$ graded, we know that
\begin{equation} \label{Lnull1}
    L_0 = \textstyle\sum_{i\ne j} [E_{ij}(A), E_{ji}(A)] = \Span_\kk \{ H_i(a,b): 1\le i < n\}.
\end{equation}
The second equality follows from the general formula \eqref{eq-hadi0}, which
in our setting can of course be easily verified:
\[
  abE_{ii} - baE_{jj} =\textstyle  \big( \sum_{k=i}^{j-2} h_k(ab) \big)
       - H_{j-1}(a,b) \qquad  (1\le i < j-1 < n).
  \]
Observe that \eqref{Lnull1} implies that every $x=\sum_{i=1}^n x_i E_{ii}\in
L_0$ has trace $\tr(x) = \sum_{i=1}^n x_i \in [A,A]$. Conversely,
since any diagonal $x=\sum_i x_i E_{ii}\in \gl_n(A)$ can be written in the
form
\begin{align}
x &=\tr(x) E_{11} + \tsum_{j=2}^n x_j (E_{jj} -E_{11})
      = \tr(x) E_{11} + \sum_{j=2}^n [x_j E_{j1}, E_{1j}]
 \nonumber  \\  \label{Lnull3} & =
 \tr(x) E_{11} - \tsum_{j=2}^n \big( h_1(x_j) + \cdots + h_{j-1}(x_j)\big)
\end{align}
it follows from \eqref{Lnull00} that any diagonal $x\in \gl_n(A)$ with
$\tr(x) \in [A,A]$ lies in $L_0$. Thus
\begin{equation}\label{sln-des}
   \lsl_n(A) =\{x \in\gl_n(A) :  \tr(x) \in [A,A]\}.
\end{equation}  This formula justifies the notation $\lsl_n$. Indeed, if $A$
is commutative then $\lsl_n(A)$ consists of the trace-less  matrices in
$\gl_n(A)$ and is in fact the base ring extension of $\lsl_n(\kk)$ by $A$,
i.e., $\lsl_n(A) \cong \lsl_n(\kk) \ot_\kk A$.

For later use we will establish some more formulas for $L_0$. We put
\[ H_j(A,A)= \Span \{ H_j(a,b) : a,b\in A\}
 \quad \hbox{and} \quad h_j(A) = \{ h_j(a) : a \in A\}.
\]
Obviously
\begin{equation} \label{Lnull4}
   L_0' \defeq \{ l\in L_0 : \tr(l) = 0 \} = \textstyle
       \bigoplus_{j=1}^{n-1} h_j(A).
\end{equation}
Moreover, \eqref{Lnull3} implies that $L_0 = [A,A]E_{11} \oplus L_0' =
H_1(A,A)+ L'_0$ where the second equality follows from \eqref{Lnull00}. By
symmetry we thus have for any $j$, $1\le j <n$,
\begin{equation}  \label{Lnull5}
L_0 = [A,A] E_{jj} \oplus L_0' = H_j(A,A) + L_0'.
\end{equation}

For $a,b,c\in A$ we abbreviate
\[ \{ a\ b\ c\} = abc + cba, \]
the Jordan triple product of $A$. It enters in the description of the algebra
$L$:
\begin{equation}\label{Lnull2}\begin{split}
  [H_i(a,b), \, H_i(c,d)] & = H_i(\{a\ b\ c\}, d) - H_i(c, \{b\ a\ d\}), \\
 [H_i(a,b), \, e_i(c) ] & = e_i(\{ a\ b\ c\}) \quad \hbox{and} \\
  [H_i(a,b), f_i(c)] &=  - f_i(\{b\ a\ c\}).
 \end{split}\end{equation}
A special case of \eqref{Lnull2} is
\begin{equation}
  \label{Lnull4a} [ h_i(a), h_i(b)] = [a,b] (E_{ii} + E_{i+1, i+1}).
\end{equation}

\subsection{Identities in $\U(L)$ and $\U(L_0)$.}\label{U0id}
In this subsection we will establish several identities involving products
with factors $e_i(a), f_i(b)$ for $a,b\in A$ and a fixed $i$, see
\eqref{rev-liealg1} for the definition of $e_i$ and $f_i$. For better
readability we put
\begin{equation}\label{U0id00}
 e(a) = e_i(a), \quad H(a,b) = H_i(a,b),
    \quad h(a) = h_i(a)\quad \hbox{and}  \quad f(a) = f_i(a)
\end{equation}
for $a,b\in A$. \sm

Applying the basic commutation rules $[e(a), f(b)] = H(a,b)$ and $[e(a),
H(b,c)] = e(-\{ a\, c\, b\})$ in $\U(L)$ repeatedly to the product $e(a_1)
\cdots e(a_t) \in \U(L)$, $t\in \NN_+$, yields
\begin{align*}
  [e(a_1) \cdots e(a_t), f(b)] &=
   \sum_{j=1}^t e(a_1) \cdots e(a_{j-1}) H(a_j,b) e(a_{j+1}) \cdots e(a_t) \\
 [ e(a_1) \cdots e(a_t), H(b,c)] &= \sum_{m=1}^t e(a_1) \cdots e(a_{m-1}) e(-\{ a_m\ c\ b\})e(a_{m+1}) \cdots e(a_t) \\
   &= - \sum_{m=1}^t  e(\{a_m\ c\ b\}) e(a_1) \cdots \widehat{e(a_m)} \cdots e(a_t),
  \end{align*}
where in the last equality we used $[e(A), e(A)] = 0$ in $\U(L)$ and where as
usual $\widehat{x}$ indicates that $x$ has been omitted. Combining these two
formulas we get for $t\in\NN_+$
\begin{align*}
 [ e(a_1) \cdots e(a_t), f(b)] &= \sum_{j=1}^t H(a_j,b) e(a_1) \cdots \widehat{e(a_j)} \cdots e(a_t) \\
   & \quad - \sum_{1\le j<m\le t} e(\{ a_j\ b\ a_m\}) e(a_1) \cdots
        \widehat{e(a_j)} \cdots \widehat{e(a_m)} \cdots e(a_t).
\end{align*}
Applying the Harish-Chandra homomorphism $\pi_0 \co \U(L)^0 \to \U(L_0)$
defined in  \ref{lem:HC-hom} and keeping in mind that $\pi_0\big(f(b_1)\cdots
f(b_t)e(a_1)\cdots e(a_t)\big)=0$ we now obtain, again for $t\in\NN_+$,
\begin{align} \label{longcalcc}\notag
\pi_0&\big(
e(a_1)\cdots e(a_t)f(b_1)\cdots f(b_t) \big) = \\
&= \sum_{j=1}^t H(a_j,b_1)\pi_0\big(e(a_1)\cdots \widehat{e(a_j)}\cdots e(a_t)f(b_2)\cdots f(b_t)\big)\\
 &\ \ \ \ \ \ \ \ \ \ - \sum_{1\leq j<m\le t}\pi_0\big(e( \{ a_j\ b_1\ a_m\})
 e(a_1)\cdots \widehat{e(a_j)} \cdots \widehat{e(a_m)}\cdots e(a_t)f(b_2)\cdots
   f(b_t)\big), \notag \end{align}
cf.\ \cite[p.~62]{Sel81}. Denoting by $\U(L_0)_{(t)}$ the $t^{\rm th}$-term
of the standard filtration of $\U(L_0)$ we now claim for any $t\in \NN_+$
\begin{equation}\label{U0id1}
 \pi_0\big( e(a_1)\cdots e(a_t)f(b_1) \cdots f(b_t)\big)\equiv \sum_{\si \in \frS_t} H(a_1, b_{\si(1)}) \cdots H(a_t, b_{\si(t)}) \mod \U(L_0)_{(t-1)}. \end{equation}
This is obvious for $t=1$, and follows for $t\ge 2$ by induction. Indeed,
assuming \eqref{U0id1} for $t-1$ the second term on the right hand side of
\eqref{longcalcc} lies in $\U(L_0)_{(t-1)}$ whence
\begin{align*}
  \pi_0&\big( e(a_1)\cdots e(a_t)f(b_1) \cdots f(b_t)\big) \equiv
      \sum_{j=1}^t H(a_j,b_1)\pi_0 \big( e(a_1)\cdots \widehat{e(a_j)}\cdots e(a_t)f(b_2)\cdots f(b_t)\big) \\
   & \equiv \sum_{j=1}^t H(a_j, b_1) \Big( \sum_{\si \in \frS_t, \si(j)=1} H(a_1, b_{\si(1)} )\cdots \widehat{H(a_j, b_1)} \cdots H(a_t, b_{\si(t)})\Big) \mod \U(L_0)_{(t-1)}.
   \end{align*}
Using the PBW Theorem we can move $H(a_j,b_1)$ to the $j^{\rm th}$ place
modulo terms in $\U(L_0)_{(t-1)}$, finishing the proof of \eqref{U0id1}. In
particular,
\begin{equation}\label{U0id1a}
 \pi_0\big( e(a_1)\cdots e(a_t)f(b)^t\big)\equiv t!\,  H(a_1, b)\cdots H(a_t, b) \mod \U(L_0)_{(t-1)}. \end{equation}

\subsection{Proposition}\label{Lema} {\it 
We fix $i$, $1\le i < n$, and assume that $\eta \co \U(L_0)\to B$ is a unital
algebra homomorphism
Let $\rho \co A \to B$ be the $\kk$-linear map defined by
$\rho(a) = \eta\big( h_i(a)\big)$. \sm

\begin{inparaenum}[\rm (a)] \item \label{Lema0}
The sequence $(g_t)_{t\in \NN_+}$ of
functions, \begin{equation} \label{Lema0gt}
       g_t \co A^t \to B, \qquad g_t(a_1, \ldots, a_t) = (\eta \circ \pi_0)\,
        \big(e_i(a_1) \cdots e_i(a_t ) f_i(1_A)^t \big)
\end{equation}
satisfies the recursion {\rm \ref{rec}} with respect to $\rho$. \sm

\item \label{Lema3} Let $\ell \in \NN_+$, and assume that $\rho$ is  a Lie homomorphism.
 Then $\rho$ satisfies  the
    $\ell^{\rm th}$-symmetric identity if and only if
\begin{equation} \label{Lema-ad}
  (\eta \circ \pi_0)\,
    \big(e_i(a_1) \cdots e_i(a_\ell) f_i(b_1) \cdots f_i(b_\ell)\big) = 0
\end{equation}
    for   arbitrary $a_j , b_j \in A$.
\end{inparaenum} }

\begin{proof} \eqref{Lema0}
We will employ the notation \eqref{U0id00}.
First note that
the condition
\ref{rec}\eqref{rec1} holds by definition: $g_1(a) = (\eta \circ \pi_0)
\big(e(a) f(1_A)\big) = \eta\big(h(a)\big) = \rho(a)$. To verify the
condition  \ref{rec}\eqref{rec2} we suppose that $(a_1, \ldots, a_{t+1})$
is a family of commuting elements of $A$. Then, using \eqref{longcalcc} with
$b_i=1_A$ we get
\begin{align*}
   g_{t+1}(a_1, \ldots, a_{t+1}) = \sum_{j=1}^{t+1} \rho(a_j) g_t(a_1, \ldots, \widehat{a_j}, \ldots a_{t+1})
    - 2 \sum_{1\le j < m \le t} g_t(a_ja_m, a_1, \ldots, \widehat{a_j}, \ldots, \widehat{a_m}, \ldots, a_{t+1})
\end{align*}
since $\{a_j\ 1_A\ a_m\} = a_j a_m + a_m a_j = 2a_ja_m$. \sm

\eqref{Lema3}
 We start by noting that
 \[
 [\eta(h(a)),\eta(h(b))]=\eta(h([a,b]))
 =\eta([a,b]E_{i,i}-[a,b]E_{i+1,i+1})
 \]
since $\rho$ is a Lie homomorphism. But $\eta$ is also a homomorphism of
associative algebras, so that the left hand side of the last equation is
equal to $\eta\big([a,b](E_{i,i}+E_{i+1,i+1})\big)$. It follows that
  $\eta([a,b]E_{i+1,i+1})=0$, and the latter relation can be written as
\begin{align}
  \label{etaHetah-Hadi}
  \eta(H(a,b))=\eta(h(ab)).
\end{align}
Next we show for $t\in \NN_+$ and $a_k, b\in A$ that
\begin{equation}\label{Lema1}
(\eta \circ \pi_0)\, \big(e(a_1) \cdots e(a_t) f(b)^t\big) =
   (\eta \circ \pi_0)\, \big(e(a_1b) \cdots e(a_tb) f(1_A)^t\big).
\end{equation}
Since $\pi_0\big( e(a)f(b) \big) = H(a,b)$, our assertion holds for $t=1$.
Let now $t>1$. Then, by  \eqref{etaHetah-Hadi}, \eqref{longcalcc} and
induction,
\begin{align*} (\eta \circ \pi_0) & \big(
e(a_1)\cdots e(a_t)f(b)^t \big)
 = \sum_{j=1}^t \eta\big(H(a_j,b)\big) \, (\eta \circ \pi_0) \big(e(a_1)\cdots
        \widehat{e(a_j)}\cdots e(a_t)f(b)^{t-1} \big)\\
 &\quad   - \sum_{1\leq j<m\le t}(\eta \circ \pi_0) \big(e( \{ a_j\ b\ a_m\}) e(a_1)
   \cdots \widehat{e(a_j)} \cdots \widehat{e(a_m)}\cdots e(a_t)f(b)^{t-1}\big)  \notag \\
 &= \sum_{j=1}^t  \eta\big( h(a_jb)\big) \, (\eta \circ \pi_0) \big(e(a_1b)\cdots \widehat{e(a_jb)}\cdots e(a_tb)f(1_A)^{t-1} \big)\\
 &\quad - \sum_{1\leq j<m\le t}(\eta \circ \pi_0) \big(e( \{ a_j\ b\ a_m\}b) e(a_1b)\cdots \widehat{e(a_jb)} \cdots \widehat{e(a_mb)}\cdots e(a_tb)f(1_A)^{t-1}\big)
 \end{align*}
The same calculation applied to the right hand side of \eqref{Lema1} yields
  \begin{align*} &(\eta \circ \pi_0) \big(
e(a_1b)\cdots e(a_tb)f(1_A)^t \big)
  = \sum_{j=1}^t \eta\big( h(a_jb)\big) \, (\eta \circ \pi_0) \big(e(a_1b)\cdots
      \widehat{e(a_jb)}\cdots e(a_tb)f(1_A)^{t-1} \big)\\
 &\ \ \ \ \ \ \ \ \ \ - \sum_{1\leq j<m\le t}(\eta \circ \pi_0)
      \big(e( \{ a_jb,\ 1_A,\ a_mb\}) e(a_1b)\cdots \widehat{e(a_jb)} \cdots
            \widehat{e(a_mb)}\cdots e(a_tb)f(1_A)^{t-1}\big).  \end{align*}
The equality \eqref{Lema1} now follows by comparing the two equations and
using $\{a_j\ b\ a_m\}b = a_jba_mb+ a_mba_jb = \{ a_jb,\ 1_A,\ a_mb\}$.

Since $[f(A), f(A)] = 0$ in $L$, the function $A^{\ell} \to \U(L_0)$, $(b_1,
\ldots ,b_{\ell}) \mapsto e(a_1) \cdots e(a_{\ell}) f(b_1) \cdots f(b_{\ell
})$, is symmetric. Hence \eqref{Lema-ad} is equivalent to $(\eta \circ
\pi_0)\, \big(e(a_1) \cdots e(a_{\ell} ) f(b)^{\ell}\big) = 0$ for arbitrary
$a_i, b \in A$. Using \eqref{Lema1}, this is in turn equivalent to
\begin{equation} \label{Lema2}
 g_\ell(a_1, \ldots, a_\ell) =
(\eta \circ \pi_0)\, \big(e(a_1) \cdots e(a_{\ell} ) f(1_A)^{\ell}\big) = 0,
\end{equation}
Since the functions $g_\ell$ are symmetric, \eqref{Lema-ad} is also
equivalent to $g_{\ell} (a, \ldots, a) = 0$. But by \eqref{Lema0} and
\eqref{rec3} this holds if and only if $\rho$ satisfies the $\ell^{\rm
th}$-symmetric identity.
\end{proof}

\subsection{Corollary}\label{cor:od} {\it Assume $\la = \sum_{j=1}^{n-1} \ell_j \varpi_j
\in \scP_+$ satisfies $\ell_{i+1} = 0$ for some $i$, $1\le i < n-1$. Denote
by $\can \co \U(L_0) \to \Se^\la (L,\g) = : \Se^\la$ the canonical
epimorphism and define $\rho_i \co A \to \Se^\la$, $\rho_i (a) = \can
\big(h_i(a)\big)$. \sm

Then $\rho_i$ satisfies the $(\ell_i + 1)^{\rm st}$-symmetric identity and
$\rho(1_A) = \ell_i 1_{\Se^\la}$. Hence there exists a unique homomorphism of
unital assocative algebras
\[ \vphi_i \co \TS^{\ell_i}(A) \to \Se^\la, \quad \vphi_i \big( \sym_{\ell_i}
(a)\big)  = \can \big( h_i(a) \big). \]
}

\begin{proof}
  From \ref{def:adm}\eqref{def:adm3} with $i$ replaced by $i+1$ we obtain
  $\pi_0\big(e_{i+1}(a) f_{i+1}(b)\big) = H_{i+1}(a,b) \in J^\la$, whence
  $\can \big(H_i(a,b) - h_i(ab) = \can([a,b]E_{i+1, i+1}) = 0$ since
  $[a,b]E_{i+1, i+1} \in H_{i+1}(A,A)$ by \eqref{Lnull00}. Thus $\eta = \can$
  satisfies the assumptions of Proposition~\ref{Lema}. The map $\rho\co A \to
  \Se^\la$ defined there is our map $\rho_i$. Again by \ref{def:adm}\eqref{def:adm3}
the function $g_{\ell_i + 1}$ of \eqref{Lema0gt} vanishes, whence by
\eqref{rec3} the function $\rho_i$ satisfies the $(\ell_i + 1)^{\rm
st}$-symmetric identity. We also have $\rho_i(1_A) = \can\big(h_i(1_A)\big) =
\can \big(\la(h_i) \idu) = \ell_i 1_{\Se^\la}$. The last part of our claim is
an application of Theorem~\ref{sel-thm}. \end{proof}

\subsection{Lemma}\label{LemBsl2} {\it
Let $L=\lsl_2(A)$ for $A\in \kalg$, and let $\la = \ell
\varpi_1\in \scP_+$ for some $\ell \in \NN$. We further suppose that\/
$\sfB'$ is a basis of $A$ containing $1_A$ and that $\le$ is a total order on
$\sfB=\sfB'\setminus \{1_A\}$. Then
\[
   \{h_1(b_1) \cdots h_1(b_k) + J^\la : b_i \in \sfB, b_1 \le b_2 \cdots
                 \le b_k, k\in \NN, k\le \ell,  \}
\]
is a spanning set of the vector space $\Se^\la(L, \g)$.
}

\begin{proof}
Since $A$ is commutative, we have $L_0 = h_1(A)$ by \eqref{Lnull5}, whence
$L_0$ is commutative too. It then follows from the PBW Theorem that
$\Se(L,\g)^\la$ is spanned by
\begin{equation}\label{LemBsl21}
   h_1(b'_1) \cdots h_1(b'_k) + J^\la, \quad b'_i \in \sfB', b'_1 \le b'_2 \le
       \cdots \le b'_k, \quad k\in \NN,
\end{equation}
with $k=0$ being the identity element $\idu$. In such an expression we can
eliminate a factor $b'_i = 1_A$ using the relation $h_1(1_A) - \ell \idu \in
J^\la$. Thus we can assume that all $b'_i = b_i \in \sfB$. We can then bound
the length $k$ of the products using \eqref{U0id1a}:
\[ \textstyle \frac{1}{(\ell + 1)!}
    \pi_0\big( e_1(b_1)\cdots e_1(b_{\ell +1})f_1(1)^{\ell +1} )
       \equiv h_1(b_1)\cdots h_1(b_{\ell +1}) \mod \U(L_0)_{(\ell)}.
\]
Since the left hand side lies in $J^\la$, any expression \eqref{LemBsl21}
with $k> \ell$ can be iteratively reduced to a linear combination of elements
\eqref{LemBsl21} with $k \le \ell$ (observe that one can achieve an
increasing order of the $b_i$'s using commutativity of $L_0$).
\end{proof}

We now have enough preparation to identify
 $\Se^\la(L,\g)$ when $L=\g\ot A$, defined in Example~\ref{rg-exa}\eqref{reg-exa1}. We use
the notation~\ref{g-not}.

\subsection{Theorem}\label{thm:comm}
{\it 
Let $L=\g\ot A$ for some $A\in \kkalg$ and let $\la = \sum_{i=1}^r \ell_i
\varpi_i \in \euP_+$. Then there exists an isomorphism
\[ \Se^\la(L,\g) \to
\TS^\la(A) \defeq \TS^{\ell_1}(A) \ot_\kk \cdots \ot_\kk \TS^{\ell_r}(A)\]
of unital associative $\kk$-algebras induced by
\[  h_i(a) \mapsto 1 \ot \cdots \ot 1 \ot \sym_{\ell_i}(a) \ot 1 \ot \cdots \ot 1 \]
with $\sym_{\ell_i}(a)$ put in the $i^{th}$ factor.
The algebra structure of\/ $\TS^\la(A)$ is that of the tensor product of the associative commutative algebras $\TS^{\ell_i}(A)$. }

\sm

For $\kk$ an algebraically closed field and $A$ a finitely generated
$\kk$-algebra with trivial Jacobson radical the isomorphism $\Se^\la(L,\g)
\cong \TS^\la(A)$ was shown in \cite[Thm.~4]{CFK10} with a different proof.
The proof of loc.\ cit.\ does not generalize to our setting.

\begin{proof}
For $1\le i \le r$ define $\g_i\simeq \lsl_2(\kk)$ as the
subalgebra generated by $e_i, f_i$ and put $L_i = \g_i \ot A = (\kk e_i \ot
A) \oplus (\kk h_i \oplus A) \oplus (\kk f_i \ot A)$. We have seen in Proposition~\ref{seli-str}
and Example~\ref{seli-se}\eqref{seli-se-a} that multiplication induces a unital algebra
homomorphism $\bigotimes_{i\in I} \Se^{\ell_i\xi_i}(L_i,\g_i) \simlgr \Se^\la(L,\g)$.
It is therefore sufficient to deal with the rank-$1$-case, i.e., $L=\lsl_2(A)$ and $\g=\lsl_2(\kk)$,
and show $\Se^\la:= \Se^\la(L,\g) \cong \TS^\ell(A)$
for $\la = \ell \om$, the isomorphism being induced by $h(a) \mapsto
\sym_\ell(a)$.

Since $A$ is commutative, hence so is $L_0 = h(A)$. By
\eqref{subsn:symm-alg1}, we have a homomorphism $L_0 \to \TS^\ell(A)^-$,
$h(a) \mapsto \sym_\ell(a)$, whence a homomorphism
\[ \eta \co \U(L_0) \to \TS^\ell(A), \quad h(a) \mapsto \sym_\ell(a)
\]
of unital associative algebras. We claim that $\eta$ annihilates the ideal
$J^\la \ideal \U(L_0)$. To prove this, recall from
\ref{def:adm}\eqref{def:adm3}  that $J^\la$ is generated by elements of type

\quad (i) $\pi_0\big( e(a_1) \ldots e(a_{\ell +1})f(b_1)\cdots f(b_{\ell
+1})\big)$, $a_i, b_i \in A$, and

\quad (ii) $h(1_A) - \ell \idu$.

\noindent It is clear that $\eta\big( h(1_A) - \ell \idu \big) = 0$. We need
to work more to deal with elements of type (i). First, note that we can apply
Proposition~\ref{Lema} with $B=\TS^\ell(A)$ because $A$ is commutative. The map
$\rho$ of loc.~cit.~ is $\rho = \sym_\ell$ and hence satisfies the
$(\ell+1)^{st}$ symmetric identity by Proposition~\ref{sel-Prop}. Thus, by
Proposition~\ref{Lema}(b) with $\ell$ replaced by $\ell + 1$, $\eta$ also
annihilates elements of type (i). Hence $\eta$ descends to a homomorphism
\[ \eta \co \U(L_0)/J^\la = \Se^\la \to \TS^\la(A), \quad u + J^\la \mapsto \eta(u)\]
of unital associative algebras. By Lemma~\ref{LemBsl2} and
Lemma~\ref{basis-new}, $\eta$ maps a spanning set of $\Se^\la(A)$ to a vector
space basis of $\TS^\ell(A)$, whence is an isomorphism.\end{proof}

\vfe

\section{$L=\lsl_n(A)$, $n\ge 3$,  and $A$ associative}\label{sec:slna}

In this section we consider the root-graded Lie algebra
\[ L=\lsl_n(A), \; \g= \lsl_n(\kk); \quad  n\ge 3, \] as defined in \ref{subsec:rev-liealg}.
We point out that $A$ is a unital associative $\kk$-algebra which need not be
commutative. Our goal in this section is to describe  $\Se^\la =
\Se^\la(L,\g)$ for special weights $\la \in \euP_+$. We do this for $\la$
totally disconnected in Theorem~\ref{disc-th}. Furthermore, in Proposition
\ref{lemmaASlaHadi} we show that for $(L,\g)=(\lsl_4(A),\lsl_4(k))$ and
$\la=\varpi_1+\varpi_2$, we have $\Se^\la\cong A$ for $A=\Mat_2(\kk)$,
whereas $\Se^\la(A)=\{0\}$ for $A=\Mat_{2n}(\kk)$ with $n>1$. In particular,
unlike the case when $A$ is commutative, for $\la=\varpi_1+\varpi_2$ the
algebra $\Se^\la$ is not necessarily isomorphic to $A\otimes A$.

\subsection{The isomorphism $\theta$}\label{the-iso}
Let $d\in \Mat_n(\kk)$ be the matrix with $1$'s on the second diagonal and
$0$'s elsewhere.
Note that $d^2 = \sum_{i=1}^nE_{i,i}$. The map
\[
   \theta = \theta_A \co \lsl_n(A) \to \lsl_n(A\op), \quad x \mapsto - d x^t d
\]
is an isomorphism of Lie algebras with $\theta_A^{-1} =\theta_{A\op}$. To
simplify notation in the following we will view $\g=\lsl_n(\kk)$ as a
subalgebra of both $\lsl_n(A)$ and $\lsl_n(A\op)$, and do the same for $\h$.
Since $\theta(aE_{ij}) = - a E_{n+1-j, n+1-i}$ we have $\theta(\g) = \g$ and
$\theta(\h) = \h$. We denote by $\theta^* \co \h^* \to \h^*$ the transpose of
$\theta|_\h = (\theta|_\h)^{-1}$. Then $\theta(\lsl_n(A)_\al) =
\lsl_n(A\op)_{\theta^*(\al)}$ holds for all $\al \in \De \cup \{0\}$, and
$\theta^*$ induces the non-trivial automorphism of the Dynkin diagram of
$(\g,\h)$. Hence
\begin{equation}
   \label{the-iso-1}
\theta^*(\varpi_i) = \varpi_{n-i} \end{equation}

Assigning to a representation $\rho \co \lsl_n(A)\to \gl(V)$ the representation $\rho
\circ \theta_{A\op} \co \lsl_n(A\op) \to \gl(V)$ gives rise to an isomorphism
between the representation categories of $\lsl_n(A)$ and $\lsl_n(A\op)$. It
preserves integrable representations and because of \eqref{the-iso-1} induces an isomorphism
\begin{equation}
  \label{the-iso-2} \scI(\lsl_n(A), \g)^\la \; \simlgr  \; \scI(\lsl_n(A\op),\g)^{\theta^*(\la)}.
   \end{equation}
The functors $\bfR$ and $\bfI$ of \ref{func} together with the isomorphism \eqref{the-iso-2} induce
isomorphisms of the module categories
$\Se^\la\MOD$ and $\Se^{\theta^*(\la)}\MOD$, that is,
\[ \xymatrix{
   \scI(\lsl_n(A), \g)^\la \ar[r]^{\simeq\ \ \ \ \ }
     & \scI(\lsl_n(A\op), \g)^{\theta^*(\la)} \ar[d]^{\bfR} \\
    \Se^\la\MOD \ar[u]_{\bfI} \ar@{-->}^\sim[r] & \Se^{\theta^*(\la)}\MOD
}\]
whence the isomorphisms
\begin{equation}
  \label{the-iso-3} \Se^\la(\lsl_n(A), \g) \; \simlgr \; \Se^{\theta^*(\la)}(\lsl_n(A\op), \g).
\end{equation}
of the associative algebras.
The isomorphism \eqref{the-iso-3} can of course also be derived directly from the definitions.
Indeed, the canonical extension of $\theta$ to an isomorphism $\U(\theta) \co \U(\lsl_n(A)) \to
\U(\lsl_n(A\op))$ maps $\U(\lsl_n(A)_0)$ to $\U(\lsl_n(A\op)_0)$, is equivariant with respect to the
Harish-Chandra homomorphisms, maps the ideal $J^\la\ideal \U(\lsl_n(A)_0)$ of \ref{def:adm}\eqref{def:adm3}
to the ideal $J^{\theta^*(\la)} \ideal \U(\lsl_n(A\op))$ and thus descends to the isomorphism
\eqref{the-iso-3}. \sm

Recall the subalgebras $L_i=L_i(A)$ of $L=\lsl_n(A)$ and $\Se^\la_i$ of $\Se^\la$
defined in \ref{subi}. In our setting, for $i\in I= \{1, \ldots, n-1\}$ we have
$L_i(A) = e_i(A) \oplus H_i(A,A) \oplus f_i(A)$. It is an $\rma_1$-graded subalgebra of $L$
with grading subalgebra $\g_i = e_i(1_A) \oplus h_i(1_A) \oplus f_i(1_A)$ and $0$-part $L_{i0} = H_i(A,A)$.
Restricting the isomorphism $\theta$ to $L_i(A)$ shows
\begin{equation}\label{the-iso-4}
   (L_i(A), \g_i) \simeq (L_{n-i}(A\op), \g_{n-i}).
\end{equation}
As before we let $\can \co \U(L_0) \to \Se^\la$ be the canonical
algebra homomorphism. Then
\begin{equation}\label{the-iso-4a}
 \Se^\la_i = \Se^\la_i(A) = \lan \can\big(H_i(A,A)\big)\ran
\end{equation}
is the unital subalgebra of $\Se^\la$ generated by $\can\big(H_i(A,A)\big)$. The isomorphism
\eqref{the-iso-4} induces an isomorphism
\begin{equation*} \label{the-iso-5}
   \Se^\la_i(A) \simeq \Se^{\theta^*(\la)}_{n-i}(A\op)
\end{equation*}
for all $i\in I$.

\subsection{The structure of $L_0 = \lsl_n(A)_0$ again.} \label{Lnua} We fix $k\in I=\{1, \ldots, n-1\}$ and
recall the decomposition \eqref{Lnull5}:
\begin{equation}   \label{Lnua1}
 L_0 = [A,A]E_{kk} \oplus \big( \oplus_{i\in I} h_i(A) \big).
\end{equation}
Observe for $a,b\in A$ and $i\in I$
\begin{align}
   aE_{ii} &=   \nonumber
    \begin{cases}
         h_i(a) + \cdots + h_{k-1}(a) + a E_{kk}, & i<k. \\
         aE_{kk} - \big( h_k(a) + \cdots + h_{i-1}(a) \big), & i>k.
    \end{cases} \\
  H_i(a,b) & = abE_{ii} - baE_{i+1, i+1} \nonumber  \\ &=
      \begin{cases}
         h_i(ab) + h_{i+1}([a,b]) + \cdots + h_{k-1}([a,b]) + [a,b]E_{kk} , & i+1 < k, \\
        h_i(ab) + [a,b]E_{kk}, & i+1 = k, \\
         [a,b]E_{kk} + h_k(ba), & i=k, \\
        [a,b]E_{kk} + h_k([b,a]) + \cdots + h_{i-1}([b,a]) + h_i(ba), &i>k.
      \end{cases} \label{Lnua-2}
\end{align}
We use this to describe the Lie algebra structure of $L_0$ in terms of the decomposition
\eqref{Lnua1}: \begin{align}
   [A,A] E_{kk} & \ideal L_0; \label{Lnua-3} \\
   [h_i(a), h_{i+1}(b)] &= [b,a]E_{i+1,i+1} \nonumber
     \\ &  = \begin{cases}
          h_{i+1}([b,a]) + \cdots + h_{k-1}([b,a]) + [b,a]E_{kk}, & i+1 < k, \\
         h_k([a,b]) + \cdots  + h_i([a,b])+ [b,a]E_{kk},  & i+1 > k;
           \end{cases}  \label{Lnua-4}\\
  [h_i(a), h_i(b)] &= [a,b] (E_{ii} + E_{i+1, i+1} ) \nonumber
     \\ &= \begin{cases}
           h_i([a,b]) + 2 h_{i+1}([a,b]) + \cdots 2 h_{k-1}([a,b]) + 2[a,b]E_{kk}, & i+1 < k, \\
           h_i([a,b]) + 2[a,b]E_{kk},   & i+1 = k, \\
             2[a,b]E_{kk} - h_k([a,b], & i=k, \\
           2[a,b]E_{kk} - 2 \big( h_k([a,b]) + \cdots + h_{i-1}([a,b]) - h_i([a,b]), &i>k;
     \end{cases}  \label{Lnua-5} \\
 [h_i(a), h_j(b) &= 0 \quad \hbox{if $|i-j| > 1$.} \label{Lnua-6}
\end{align}

\subsection{Lemma}\label{conA}
{\em
As in {\rm \ref{Lnua}} we fix $k\in I$. Let $\rho_i \co A \to B$, $i\in I$, be
linear maps into a unital associative $\kk$-algebra $B$ satisfying \sm

\quad {\rm (I)} $[\rho_i(A), \rho_j(A)]= 0 $ if $i\ne j$;

\quad {\rm (II)} $\rho_i([A,A]) = 0$ if $1<i<n-1$ and also for $i=1$ in case $k=1$;

\quad {\rm (III)} $\rho_i$ is a Lie homomorphism for $1\le i <k$, and a
     Lie anti-homomorphism for $k\le i\le n-1$.
\sm

Then the map
\[ \eta_{L_0} \co L_0 \to B, \quad \textstyle
      c E_{kk} + \sum_{i\in I} h_i(a_i) \mapsto \sum_{i\in I} \rho_i(a_i)
      \]
is a Lie homomorphism. }

\begin{proof}
This is immediate from the multiplication rules \eqref{Lnua-3}--\eqref{Lnua-6}.

\end{proof}

\subsection{Lemma}\label{abelian-generators-Hadi}
{\em Let $\la = \sum_{i\in I} \ell_i \varpi_i \in \euP_+$ be totally disconnected,
cf.\ {\rm \ref{seli-se}\eqref{seli-se-b}}.
Then the subalgebras $\Se^\la_i$ of $\Se^\la$, cf.\ {\rm \eqref{the-iso-4a}},
have the following properties.

\begin{inparaenum}[\rm (a)]

\item \label{agH-a} $\Se^\la_i$ is generated by $\can\big(h_i(A)\big)$ as associative algebra.

\item \label{agH-b} If $\ell_i = 0$, then $\Se^\la_i = \kk 1_{\Se^\la}$.

\item \label{agH-c} If $1< i < n-1$, then $\Se^\la_i$ is commutative and $h_i([A,A]) = 0$.

\item \label{agH-d} $[\Se^\la_i(A), \Se^\la_j(A)] = 0$ for $i\ne j$.
\end{inparaenum}}

\begin{proof} Recall from \ref{def:adm}\eqref{def:admd} that $H_k(A,A) \subset J^\la$
whenever $\ell_k = 0$, which implies \eqref{agH-b}. The identities \eqref{Lnull00},
\begin{equation*}\label{lemBn1n}
H_j(a,b) - h_j(ba) = [a,b]E_{jj} = H_{j-1}(a,b) - h_{j-1}(ab), \qquad j \in I, j>1,
\end{equation*}
imply that $\can (H_i(A,A)) = \can (h_i(A))$ since the image under $\can$ of the
left or right side of the equation vanishes. This proves \eqref{agH-a}. Moreover, for $1< i < n-1$
with $\ell_i >0$ the equation also implies that $[h_i(a), h_i(b)] = [a,b] (E_{ii} + E_{i+1, i+1} )
\in H_{i-1}(A,A) + H_{i+1}(A,A) \subset J^\la$, whence the first part of \eqref{agH-c}. The second follows
in the same way from $\can\big( h_i([a,b])\big) = \can \big( [a,b](E_{ii} - E_{i+1, i+1})$.
Finally, for the proof of
\eqref{agH-d} we can assume that $|i-j| >1$, in which case already $[H_i(A,A), H_j(A,A)]= 0$ in
$\U(L_0)$. \end{proof}

We are now in a position to prove the main result of this section. For simpler notation we put $\TS^0(A) =
\kk 1$ and let $\sym_0 \co A \to \TS^0(A)$ be the zero map.

\subsection{Theorem}\label{disc-th} {\em 
Let $\la = \sum_{i\in I} \ell_i \varpi_i \in \euP_+$ be totally disconnected. For $1< i < n-1$
denote by $\scC_i$ the ideal of $\TS^{\ell_i(A)}$ generated by the commutator space
$[\TS^{\ell_i}(A), \TS^{\ell_i}(A)]$ and define
\[
  B_i = \begin{cases}
     \TS^{\ell_1}(A), & i=1 \\ \TS^{\ell_i}(A)/\scC_i, & 1<i<n-1, \\
      \TS^{\ell_{n-1}}(A\op) & i = n-1.
  \end{cases}
\]
Then
\[ \Se^\la\big(\lsl_n(A)\big)\cong \textstyle \bigotimes_{i\in I} B_i \]
as unital associative $\kk$-algebras, where the algebra structure of\/ $\bigotimes_{i\in I} B_i$
is that of the tensor product algebra.}

\begin{proof} We will construct unital algebra homomorphisms in both directions and show
that they are inverses of each other. For easier notation we put $\Se^\la = \Se^\la\big(\lsl_n(A)\big)$
and $\scC_i=\{0\} \subset B_i$ for $i=1$ and $i=n-1$. \sm

(a) {\em The homomorphism $\bar \eta \co \Se^\la \to B = \bigotimes_{i\in I} B_i$.}
We fix $k\in I$ with $\ell_k = 0$, and for $i\in I$ define
\[
   \rho_i \co A \to B, \quad \rho_i(a) = 1_{B_1} \ot \cdots \ot 1_{B_{i-1}} \ot
     \big( \sym_{\ell_i}(a) + \scC_i\big) \ot 1_{B_{i+1}} \ot \cdots \ot 1_{B_{n-1}}.
\]
We claim that the map \[ \eta_{L_0} \co L_0 \to B, \quad \textstyle
cE_{kk} + \sum_{i\in I} h_i(a) \mapsto \sum_{i\in I} \rho_i(a_i)\]
of Lemma~\ref{conA} is a Lie homomorphism. To show this, we verify
the conditions (I)---(III) of loc.\ cit.

Identifying $B_i$ with the obvious unital subalgebra of $B$, we have $\rho_i(A) \subset B_i$,
so that (I) holds because of $[B_i, B_j] = 0$ for $i\ne j$. Recall from \eqref{subsn:symm-alg1} that
$a \mapsto \sym_\ell(a)$ is a Lie homomorphism, whence $\rho_i([A,A]) = [\sym_{\ell_i}(A),
\sum_{\ell_i}(A)] + \scC_i = \scC_i$ for $1<i<n-1$, and (II) follows in that case. Since
$\rho_1= 0$ in case $k=1$, we have now proven (II). The condition (III) is clear
from the proof of (II) if $1\le i \le k$, and also for $k \le i < n-1$ since then $\rho_i([a,b]) = 0 =
[\rho_i(b), \rho_i(a)]$. Finally, (III) for $i=n-1$ follows from the fact that
$\sym_{\ell_{n-1}} \co A\op \to \TS^{\ell_{n-1}}(A\op) = \TS^{\ell_{n-1}} (A)\op$
is a Lie homomorphism. \sm

Having proven that $\eta_{L_0}$ is a Lie homomorphism, we next claim that the unique extension of $\eta_{L_0}$ to a unital
associative algebra homomorphism $\eta \co \U(L_0) \to B$ annihilates the ideal $J^\la$
defining  $\Se^\la$. Thus, specializing \ref{def:adm}\eqref{def:admd} to our
setting, we need to show that $\eta$ vanishes on elements of the following types:
\begin{inparaenum}[(i)]

\qquad \item\label{conA-1}  $H_i(a,b)=\pi_0\big(e_i(a) f_i(b)\big)$ for
        $i\in I$ with $\ell_i= 0 $;

\qquad \item  \label{conA-2} $\pi_0\big(e_i(a_1)\cdots e_i(a_{\ell_i
    +1})f_i(b_1)\cdots f_i(b_{\ell_i+1})\big)$ for  $a_r,b_s \in A$ and $\ell_i >0$;

\qquad \item\label{conA-3} $h_i(1_A)-\ell_i \idu$ for $i\in I$ with $\ell_i >0$.
\end{inparaenum}
\sm

Re \eqref{conA-1}: We have $\eta\big( H_i(a,b)\big) = \eta_{L_0}\big( H_i(a,b)\big)$ which
can be calculated using the formulas in \eqref{Lnua-2}. For example, if $i+1< k$ then
$\eta_{L_0}\big(H_i(a,b)\big)= \rho_i(ab) + \rho_{i+1}([a,b]) + \cdots + \rho_{k-1}([a,b])=0$
using $\rho_i = 0$ because $\ell_i = 0$ and $\rho_j ([a,b]) = 0$ for $1<j<n-1$ as shown in
the proof of (II) above. The proof in the other cases is similar.

Re \eqref{conA-2}: By Proposition~\ref{sel-Prop} the map $A \to \TS^{\ell_i}(A)$ satisfies the
$(\ell_i + 1)^{\rm st}$-symmetric identity, whence so does $\rho_i$ for $1\le i < n-1$. Since
$\rho_i(a) = \eta\big( h_i(a)\big)$, Proposition~\ref{Lema} applies and yields our claim. We
can argue similarly in case $i=n-1$: The map $\sym_{\ell_{n-1}} \co A \op \to  \TS^{\ell_{n-1}}
(A\op) = B_{n-1}$ satisfies the $(\ell_{n-1} + 1)^{\rm st}$-symmetric identity. Since powers in
$A$ and $A\op$ coincide, it follows from Definition~\ref{sel-sym} that $\rho_{n-1} \co
A \to B_{n-1}$ satisfies the same symmetric identity. We can then conclude as before by applying
Proposition~\ref{Lema}.

Re \eqref{conA-3}: These elements are annihilated because $\rho_i(1_A)) = \ell_i 1_{B_i}$ and
$\eta(1_{\U(L_0)} = 1_B$. \sm

We now know that $\eta(J^\la) = 0$. The induced algebra homomorphism $\bar \eta \co
\Se^\la \to B$ is the map we were looking for. Observe that
\begin{equation}\label{disc-th-1}
  \bar \eta \big( \can \big(h_i(a)\big)\big) = \sym_{\ell_i}(a) + \scC_i \qquad (i\in I).
\end{equation}
\sm

(b) {\em The homomorphism $\bar \vphi \co B \to \Se^\la$.}
We will use Lemma~\ref{abelian-generators-Hadi} to construct this map.
First we claim that there exists a well-defined unital algebra homomorphism
\[
   \vphi_i \co B_i \to \Se^\la, \qquad \vphi_i(\sym_{\ell_i} (a) + \scC_i)
             = \can\big( h_i(a)\big) \in \Se^\la_i. \]
This is obvious in case $\ell_i = 0$ since then $B_i = \kk 1_B$ (after identification)
and also $\can\big(h_i(a)\big)= 0$. Let then $\ell_i > 0$. In case $1 \le i < n-1$ the existence
of $\vphi_i$ follows from Corollary~\ref{cor:od}, namely directly in case $i=1$ and from
commutativity of $\Se^\la_i$ in case $1< i < n-1$. It remains to consider $i=n-1$.
Here $\theta^*(\la) = \ell_{n-1} \varpi_1 + 0 \varpi_2 + \cdots$, whence Corollary~\ref{cor:od} for
$A\op$ yields the existence of a unital algebra homomorphism
\[ \vphi_{n-1}\op \co \TS^{\ell_{n-1}}(A\op) \to \Se^{\theta^*(\la)} \big(\lsl_n(A\op)\big),
\quad \sym_{\ell_{n-1}}(a) \mapsto \can \big(h_1(a)\big).
\]
Following this map with the isomorphism $\Se^{\theta^*(\la)}_1\big( \lsl_nA\op)\big) \cong
\Se^\la_i(A)$ established in \eqref{the-iso-3} provides a map $\vphi_{n-1}\co B_{n-1} \to \Se^\la$ as claimed.

We combine the maps $\vphi_i$, $i\in I$, to get a multilinear map
$\vphi \co \textstyle \prod_{i\in I} B_i \to \Se^\la$, $
            (b_i)_{i\in I} \mapsto \prod_{i\in I} \vphi_i(b_i)$, and denote by
$\bar \vphi \co \bigotimes_{i\in I} B_i \to \Se^\la$ the unique linear map defined by
$\bar\vphi(b_1 \ot \cdots \ot b_{n-1}) \mapsto \vphi_1(b_1) \cdots \vphi_{n-1}(b_{n-1})$.
It is immediate that $\bar \vphi$ is a unital algebra homomorphism.

By construction
\begin{equation}
  \label{disc-th-2}  \bar \vphi(\sym_{\ell_i}(a) + \scC_i) = \can \big( h_i(a)\big)
  \qquad (i\in I).
\end{equation}
\sm

(c) {\em $\bar \eta$ and $\bar \vphi$ are inverses of each other.} The associative algebra
$B$ is generated by $B_i \subset B$, $i\in I$, which, as shown in Lemma~\ref{lem:symm-folklore}\eqref{sym-folk-c},
are in turn generated by $\sym_{\ell_i}(A)$. Similarly, since $L_0 = \sum_{i\in I} H_i(A,A)$,
the associative algebra $\Se^\la$ is generated by $\can\big( H_i(A,A)\big) = \can(h_i(A)$, $i\in I$. Because of
\eqref{disc-th-1} and \eqref{disc-th-2}, $\bar \eta \circ \bar \vphi$ and $ \bar \vphi \circ \bar \eta$ are
the identity map on generating sets of $B$ and $\Se^\la$ respectively, proving that they
are inverses of each other. \end{proof}

\subsection{Example $\Se^\la = \{0\}$} \label{ex:disc-th}
Let $\la= \sum_i \ell_i \varpi_i\in \euP_+$ be totally disconnected with
$\ell_1 = 0 = \ell_{n-1}$. Then $\Se^\la $ is commutative by
Theorem~\ref{disc-th} (or Lemma~\ref{abelian-generators-Hadi}\eqref{agH-c}).
If $A = (\kk 1_A + [A,A]) \oplus X$ for some complementary space $X$, then
Lemma~\ref{abelian-generators-Hadi}\eqref{agH-c} shows that $\Se^\la_i$ is
generated by $\can \big( h_i(X)\big)$ as unital algebra. Moreover, if $1_A
\in [A,A]$ then $\Se^\la= \{0\}$ since all $\Se^\la_i$ vanish in that case.
We note that there exist natural algebras $A$ with $1\in [A,A]$, for example
the Weyl algebra. See \ref{ex:seli}\eqref{ex:seli2} for a discussion of the
consequences of $\Se^\la = \{0\}$. \sm

We will use the next lemma to construct representations.

\subsection{Lemma}\label{constr-n}{\em
Let $\rho \co A \to B$ be a linear map into a unital  associative
$\kk$-algebra $B$ satisfying
\begin{equation}  \label{construct1n}
   \rho([a,b]) = 0 = [\rho(a), \rho(b)]
\end{equation}
for $a,b\in A$. For $i$ with $1<i < n-1$ define $\eta\co L_0 \to B$ by
\begin{equation} \label{def:eta22n}
 [a,b]E_{11} \oplus \big(\textstyle \bigoplus_{j=1}^{n-1}
 h_j(a_j)\big)  \mapsto \rho(a_i).
\end{equation}
Then $\eta\co L_0 \to B^-$ is a Lie homomorphism. The unique extension of
$\eta$ to a  unital associative algebra homomorphism $\U(L_0) \to B$ has the
following properties: \sm

\begin{enumerate}[\rm (i)]

\item \label{constructin} $\eta\big(h_i(a)\big) = \rho(a)$ and $\eta\big(
    H_i(a,b)\big) = \eta\big( h_i(ab)\big)$ for all $a,b\in A$; \sm

\item \label{constructiin} $\eta$ annihilates the ideal $J^\la$ if and only
    if $\rho(1_A) = \ell 1_B$ and $\rho$ satisfies the $(\ell + 1)^{\rm
    st}$-symmetric identity.
\end{enumerate}
}

Part \eqref{constructin} is the special case $k=1$ of Lemma~\ref{conA}. The
proof of part \eqref{constructiin} is similar to part (a) in the proof of
Theorem~\ref{disc-th} and will be left to the reader.

\subsection{Example $A=\euZ(A) \oplus [A,A]$} \label{ex:csan}
We denote by $\euZ(A)$ the centre of $A$ and suppose that $A= \euZ(A) \oplus
[A,A]$. Some examples  for which this is the case are:
\begin{enumerate}[(i)]
\item \label{ex:csa0n} $A$ is commutative;

 \item\label{ex:csa1n} $A$ is finite-dimensional central-simple over $\kk$;

 \item\label{ex:csa2n} $A$ is an Azumaya over some $R\in \kalg$;

 \item\label{ex:csa3n} $A$ is a quantum torus over $\kk$.

\end{enumerate}

\noindent We denote by $\ta \co A \to \euZ(A)$ the projection with kernel
$[A,A]$, and by $L_z\in \End_\kk \big(\euZ(A)\big)$ the left multiplication by $z\in
\euZ(A)$. Let $\la = \ell \varpi_i$ for some $1 <i < n-1$ and $\ell \in
\NN_+$. The map
\[ \rho \co A \to \End_\kk \euZ(A), \quad a \mapsto \ell \, L_{\ta(a)}.
\]
satisfies \eqref{construct1n}.

The homomorphism $\eta \co \U(L_0) \to \End_\kk\big( \euZ(A)\big)$, defined in
\eqref{def:eta22n}, descends to a unital homomorphism $\Se^\la \to
\End_\kk\big(\euZ(A)\big)$ and thus {\em defines a $\la$-admissible
$L_0$-module if and only if $\rho$, equivalently $\ell\tau$, satisfies the
$(\ell + 1)^{\rm st}$-symmetric identity.}   \sm

In particular, assume $\euZ(A) = \kk 1_A$. Then
Lemma~\ref{abelian-generators-Hadi}\eqref{agH-c} shows $\dim \Se^\la \le 1$.
We claim in this case
\begin{equation} \label{ex:csa-eq2n}
 \dim \Se^\la = 1  \quad \iff \quad \hbox{$\ell \tau$ satisfies the
   $(\ell + 1)^{\rm st}$-symmetric identity.}
\end{equation}
Indeed, if $\Se^\la = \kk 1_{\Se^\la}\ne \{0\}$, the map $\rho_i \co A \to
\Se^\la$ of Corollary~\ref{cor:od} is here given by $\rho_i(a) = \ell \ta(a)
1_{\Se^\la}$,  and hence can be identified with $\ell \ta \co A \to \kk$. We
then know from Corollary~\ref{cor:od} and $1_{\Se^\la}\neq 0$ that $\ell \ta$ satisfies the $(\ell +
1)^{\rm st}$-symmetric identity. Conversely, if this is the case, we have
seen above that then $\Se^\la$ has a $1$-dimensional irreducible
representation. Since the zero algebra does not have such a representation we
must have $\dim \Se^\la = 1$. \sm

The assumption $A = \kk 1_A \oplus [A,A]$ holds in case \eqref{ex:csa1n}
above. For those algebras, \eqref{ex:csa-eq2n} is shown in \cite[V.2]{Sel81}.
We note that Seligman has also shown (\cite[Cor.~III.V]{Sel81}) that
\eqref{ex:csa-eq2n} holds if and only if $\ell$ is a multiple of the degree
of $A$.

\subsection{Remark}\label{rem:compi} Let $\scC_A = \id([A,A])$ be the ideal of $A$ generated by
$[A,A]$. By exactness of the functor $\TS^\ell$, \ref{mopro}\eqref{mopro-a},
the epimorphism $c\co A \to A/\scC_A$ leads to an epimorphism $\TS^\ell(c)
\co \TS^\ell(A)\twoheadrightarrow  \TS^\ell(A/\scC_A)$ which by commutativity
of $\TS^\ell(A/\scC_A)$ factors through an epimorphism
\begin{equation}
  \label{rem:compi1}
    \bar c \co \TS^\ell(A)/\scC \twoheadrightarrow \TS^\ell(A/\scC_A).
\end{equation}
The map $\bar c$ is in general not an isomorphism. For example, if $A$ is
finite-dimensional central-simple over $\kk$ and non-commutative, i.e., $\dim
A >1$, we have $A/\scC_A =\{0\}$, whence also $\TS^\ell(A/\scC_A) = \{0\}$.
But by Example~\ref{ex:csan} we know
$\dim \TS^\ell(A)/\scC = \dim \Se^{\ell \varpi_i} = 1$ if (and only if)
$\ell$ is a multiple of the degree of $A$.

\subsection{$\Se^\lambda$
for $\lambda=\varpi_1+\varpi_2$ and $L=\lsl_4(A)$.}
Recall that a \emph{generalized quaternion algebra}
is a $\kk$-algebra with generators $\mathbf i,\mathbf j$
 which satisfy the relations
\[
\mathbf i^2=a,\ \mathbf j^2=b,\ \mathbf i\mathbf j=-\mathbf j\mathbf i\text{ where }a,b\in\kk^\times.
\]
Setting $\mathbf k=\mathbf i\mathbf j$, we obtain a basis
$\{1,\mathbf i,\mathbf j,\mathbf k\}$ for $A$ with relations
$\mathbf k^2=-ab$, $\mathbf k\mathbf i=-a\mathbf j$, and $\mathbf j\mathbf k=-b\mathbf i$.
The following proposition provides an example that indicates that beyond the case of totally disconnected $\lambda$,  the structure of  $\Se^\lambda$ might depend on $A$ as well as on $\lambda$.

\subsection{Proposition}
\label{lemmaASlaHadi}
{\em Assume that $L=\lsl_4(A)$, and
set $\lambda=\varpi_1+\varpi_2$. If $A$ is a generalized quaternion algebra over $\kk$, then the map
\[
A\to \Se^\lambda\ ,\ a\mapsto h_2(a)+J^\la\] is an isomorphism of associative
algebras. Furthermore, if  $A=\Mat_{2n}(\kk)$ for an integer $n\geq 2$, then
$\Se^\la=\{0\}$. }

\ms
In the rest of the paper,  we give the proof of Proposition \ref{lemmaASlaHadi}, which appears in the  Journal of Algebra version of this paper without proof.

Let $\lambda$ and $L$ be as in
Proposition \ref{lemmaASlaHadi}. 
The Harish-Chandra homomorphism satisfies
\[
\pi_0(e_j(a_1)e_j(a_2)f_j(b_1)f_j(b_2))=H_j(a_1, b_1)H_j(a_2, b_2)+H_j(a_2, b_1)H_j(a_1, b_2)-H_j(\{a_1b_1a_2\}, b_2)
\]
for $a_1, a_2, b_1, b_2\in A$. Using the latter relation, we obtain that the ideal $J^\lambda\subset \U(L_0)$ is generated by
the following:
\begin{itemize}
\item[(i)] $H_1(a_1,b_1)H_1(a_2,b_2)+H_1(a_2,b_1)H_1(a_1,b_2)-
H_1(\{a_1\ b_1\ a_2\},b_2)$
for $a_1, a_2, b_1, b_2\in A$.
\item[(ii)] $H_2(a_1,b_1)H_2(a_2,b_2)+H_2(a_2,b_1)H_2(a_1,b_2)-
H_2(\{a_1\ b_1\ a_2\},b_2)$
for $a_1, a_2, b_1, b_2\in A$.
\item[(iii)] $H_3(a,b)$ for $a,b\in A$.
\item[(iv)] $h_1(1)-1$, $h_2(1)-1$, $h_3(1)$.

\end{itemize}
Also recall that $L_0=H_1(A,A)+H_2(A,A)+H_3(A,A)$. Next observe that
\begin{equation}
\label{eq-H3-hadiabfall}
H_2(a,b)-h_2(ab)=[a,b]E_{3,3}=H_3(a,b)-h_3(ba)
\in J^\lambda
\end{equation}
and
\begin{equation}
\label{eq-H4-hadiabfall}
H_1(a,b)-h_1(ab)=[a,b]E_{2,2}=H_2(a,b)-h_2(ba).
\end{equation}
The last two relations together with (iii) and $h_3(a)=H_3(a,1)$  show that the algebra $\U(L_0)/J^\lambda$ is generated by $h_1(A)$ and $h_2(A)$.

Set $a_1=a$, $b_1=a_2=1$, $b_2=b$ in (i) and (ii). Since $H_i(1,1)=h_i(1)$ and $h_i(1)-1\in J^\lambda$ for $i=1,2$, we obtain
\[
h_i(a)h_i(b)+H_i(a,b)-H_i(2a,b)\in J^\lambda\ \text{ for }i=1,2,
\]
from which it follows that
\begin{equation}
\label{hihiabababa-Hadi}
h_i(a)h_i(b)-H_i(a,b)\in J^\lambda.
\end{equation} The last relation together with \eqref{eq-H3-hadiabfall} implies that $h_2(a)h_2(b)-h_2(ab)\in J^\lambda$.
Furthermore,
\[
H_1(a,b)-h_1(ab)+J^\lambda=[a,b]E_{2,2}+J^\lambda=
H_2(a,b)-h_2(ba)+J^\lambda=
h_2(ab)-h_2(ba)+J^\lambda.
\]
The latter calculation together with \eqref{hihiabababa-Hadi} for $i=1$ yields
\[
h_1(a)h_1(b)-h_1(ab)-h_2([a,b])\in J^\la.
\]
Finally, $[h_1(a),h_2(b)]=[b,a]E_{2,2}$, so that
using
\eqref{eq-H4-hadiabfall} and
\eqref{eq-H3-hadiabfall} we obtain
\[
[h_1(a),h_2(b)]+J^\la=[b,a]E_{2,2}+J^\la=H_2(b,a)-h_2(ab)+J^\la
=h_2([b,a])+J^\la.
\]
Thus, from the above discussion it follows that $\Se^\la$ is generated by $h_1(A)$ an $h_2(A)$ which satisfy the following
relations mod $J^\la$:
\begin{equation}
\label{casesh1h2[]-hadi}
\begin{cases}
h_1(a)h_1(b)=h_1(ab)+h_2([a,b])\\
h_2(a)h_2(b)=h_2(ab)\\
h_2(b)h_1(a)=h_1(a)h_2(b)+h_2([a,b])
\end{cases}
\end{equation}

\subsection{Proposition}
\label{lemmaASlaHadi-1}
{\em Assume that $A$ is a generalized quaternion algebra over $\kk$. Then the map
\[
A\to \Se^\lambda\ ,\
a\mapsto h_2(a)+J^\la\]
is an isomorphism of associative algebras.
}\begin{proof}
 From (i) above for $a_1=\mathbf i$, $b_1=\mathbf j$,  $a_2=1$, and $b_2=\mathbf i$ we obtain $\{a_1\ b_1\ a_2\}=0$ and therefore
\begin{align}
\label{H1i,jh11+h11}
H_1(\mathbf i,\mathbf j)h_1(\mathbf i)+h_1(\mathbf j)h_1(a)\in J^\la
.
\end{align} Next we have
\[
H_1(\mathbf i,\mathbf j)=h_1(\mathbf k)+2h_2(\mathbf k)
+2\mathbf kE_{3,3},
\]
and since $\mathbf k\in [A,A]$, we obtain
from \eqref{eq-H3-hadiabfall} that $\mathbf k E_{3,3}\in J^\la$.
Therefore \eqref{H1i,jh11+h11} can be written as
\[
(h_1(\mathbf k)+2h_2(\mathbf k))h_1(\mathbf i)
+a h_1(\mathbf j)\in J^\la.
\]
Expanding the left hand side of above by means of \eqref{casesh1h2[]-hadi}, we obtain
\begin{align}
\label{Equ1-Hadi}
ah_2(\mathbf j)+h_1(\mathbf i)h_2(\mathbf k)\in J^\la.
\end{align}
It follows that
\[
h_1(\mathbf i)+J^\la=
h_1(\mathbf i)h_2(\mathbf k)h_2(-\frac{1}{ab}\mathbf k)+J^\la=\frac{1}{b}h_2(\mathbf j)h_2(\mathbf k)+J^\la=\frac{1}{b}h_2(-b\mathbf i)+J^\la=
-h_2(\mathbf i)+J^\la.
\]
A similar argument shows that $h_1(\mathbf j)+h_2(\mathbf j)\in J^\la$ and
$h_1(\mathbf k)+h_2(\mathbf k)\in J^\la$. Consequently, we obtain
\[
h_1(a)=h_2(a^*)\text{ for }a\in A,
\]
where $a\mapsto a^*$ is defined by
\[
c_1+c_2\mathbf i+c_3\mathbf j+c_4\mathbf k\mapsto
c_1-c_2\mathbf i-c_3\mathbf j-c_4\mathbf k
\]
where $c_1,\ldots,c_4\in\kk$. But we know that $\Se^\lambda$ is generated by
$h_1(A)$ and $h_2(A)$. It follows that the map $a\mapsto h_2(a)+J^\lambda$ is
surjective.

Using the facts that $a+a^*$ is always in the centre of $A$ and that $[a,b]+[a,b]^*=0$ for $a,b\in A$, it is easily verified that the
map $L_0\to A$ given by $\diag(a_1,a_2,a_3,a_4)\mapsto a_1^*+a_1^{}+a_2^{}$ is indeed a Lie homomorphism. The latter map induces a homomorphism of associative algebras $\mathsf e:\U(L_0)\to A$. We prove that
$J^\la\subset\ker(\mathsf e)$. To this end, we need to verify that the generators of $J^\la$ of type (i)--(iv) are in $\ker(\mathsf e)$. This is straightforward for the generators of type (ii)--(iv).
For generators of type (i), this amounts to
verifying that
 \begin{align*} ((a_1b_1)^*+[a_1,b_1]) ((a_2b_2)^*+
 [a_2,b_2]) &+ ((a_2b_1)^*+[a_2,b_1]) ((a_1b_2)^*+[a_1,b_2])\\ &- \big(
 (\{a_1\ b_1\ a_2\}b_2)^* +
[ \{a_1\ b_1\ a_2\},b_2]
  \big)=0.
  \end{align*}
To simplify the left hand side, we use the relation
$[a,b]=[b,a]^*$ for each of the commutators. Then we can write
$(a_1b_1)^*+[a_1,b_1]=a_1^*b_1^*$ and so on. The expression then simplifies easily.

Since $J^\la\subset\ker(\mathsf e)$, we obtain a homomorphism of associative
algebras $\overline{\mathsf e}:\Se^\la=\U(L_0)/J^\la\to A$. Clearly the maps
$\overline{\mathsf{e}}$ and $a\mapsto h_2(a)+J^\la$ are inverse to each
other. This completes the proof of the proposition.
\end{proof}

\subsection{Proposition}\label{prop-A=Mat2n-hadi}
{\em Let $A=\Mat_{2n}(\kk)$ for an integer $n\geq
1$, and let $\la=\varpi_1+\varpi_2$. Then $\Se^\la=A$ for $n=1$, and
$\Se^\la=\{0\}$ for $n>1$. }

\begin{proof}
For $n=1$ the statement follows from Proposition \ref{lemmaASlaHadi-1}, because it is well known that $\Mat_2(\kk)$ is a generalized quaternion algebra for parameters $a=b=1$. From now on we assume $n>1$.

\textbf{Step 1.}
Let $A_\circ\subset A$ be the subalgebra of $A$   generated by unital $\kk$-subalgebras of $A$ that are isomorphic to $\Mat_2(\kk)$. We prove that $A_\circ=A$. To this end, for $2\times 2$
invertible matrices $p_1,\ldots,p_n$, we consider the subalgebra $B=B(p_1,\ldots,p_n)$ of $A$ that consists of block-diagonal matrices of
the form \[
\diag(p_1xp_1^{-1},\ldots,p_nxp_n^{-1})\text{ where }x\in \Mat_2(\kk).
\]
Clearly $B(p_1,\ldots,p_n)\cong \Mat_2(\kk)$. Now
let $a,b,c,d\in\kk\backslash\{0\}$ be distinct. Set \[
x=\begin{bmatrix}a& 0\\0 & b
\end{bmatrix},\
y=\begin{bmatrix}c& 0\\0 & d
\end{bmatrix},\
\text{ and }
p_i=\begin{bmatrix}1 & 1\\ 1 & x_i
\end{bmatrix}\text{ where }x_i\in\kk\backslash\{1\}
\text{ for }1\leq i\leq n.
\]
By calculating the eigenvalues of $yp_ixp_i^{-1}$, it is straightforward to show that one can choose the $x_i$ such that the eigenvalues of the block-diagonal matrix
\[
\diag(yp_1xp_1^{-1},\ldots,yp_nxp_n^{-1})
\]
are distinct elements of $\kk$. But the latter element belongs to
$B(E_{1,1}+E_{2,2},\ldots,E_{1,1}+E_{2,2})B(p_1,\ldots,p_n)$, hence it belongs to $A_\circ$. Since $A_\circ$ is invariant under the action of $\mathrm{GL}_{2n}(\kk)$ by conjugation, it follows that $A_\circ $ contains a diagonal matrix $x_\circ$ with distinct eigenvalues. We now prove that every diagonal element of $A$ is in $A_\circ$.
This is because for every diagonal matrix $y_\circ$, the linear system
\[
c_0+c_1x_\circ+\cdots c_{2n}x^{2n}_\circ=y_\circ
\]
can be solved for $c_0,\ldots,c_{2n}\in\kk$.
Finally, from $\mathrm{GL}_{2n}(\kk)$-invariance of $A_\circ$ it follows that $A_\circ=A$.

\textbf{Step 2.}  Set
\[
U=\{a\in A\ :\ h_1(a)\in h_2(A)\}.
\]
We will prove that $U=A$. From \eqref{casesh1h2[]-hadi} it follows that $U$
is a unital $\kk$-subalgebra of $A$. Furthermore, if $B\subset A$ is any
unital $\kk$-subalgebra of $A$ which is isomorphic to $\Mat_2(\kk)$, then the
proof of Proposition \ref{lemmaASlaHadi-1}  shows that $B\subset U$. From Step 1 it
follows that $U=A$.

\textbf{Step 3.} From Step 2 it follows that the homomorphism of
$\kk$-algebras $h_2:A\to \Se^\la$ is surjective. Our goal is to prove that
$h_2(A)=\{0\}$ (from which it follows that $\Se^\la=\{0\}$). Suppose that
$h_2\neq 0$. Since $A$ is simple, $h_2$ should be injective. Thus for every
$a\in A$, there exists a unique $\phi(a)\in A$ such that
$h_1(a)=h_2(\phi(a))$. It is not difficult to verify that $\phi:A\to A$ is
$\kk$-linear. Furthermore, from the first relation in
\eqref{casesh1h2[]-hadi} it follows that
\begin{equation}\label{mystphi-hadi}
\phi(ab)=\phi(a)\phi(b)-[a,b].
\end{equation}
The above relation implies
 that $\phi:A\to A$ is a homomorphism of Jordan algebras (where $A$ is
 equipped with the standard Jordan product $\{a,b\}=ab+ba$). By a well-known
 theorem of Jacobson and Rickart \cite[Remarks after Th.~10]{JaRi}, one knows that the only Jordan homomorphisms of $A$ are of the form $x\mapsto pxp^{-1}$ or $x\mapsto px^tp^{-1}$. Substituting the map $\phi(x)=pxp^{-1}$ in
 \eqref{mystphi-hadi}, we reach a contradiction. Substituting the map
$\phi(x)=px^tp^{-1}$ in \eqref{mystphi-hadi}, we obtain the relation
$pyp^{-1}=-y^t$ for every $y\in[A,A]$, which is impossible for $n>1$.
\end{proof}

\bibliographystyle{alpha}
\bibliography{global-weyl-modules-biblist}

\vfill\eject
\end{document}